\documentclass[11pt,reqno,psamsfonts]{amsart}
\pdfoutput=1

\usepackage{amssymb}
\usepackage{amsthm}
\usepackage{amsmath}
\usepackage{latexsym}
\usepackage[T1]{fontenc}
\usepackage[utf8]{inputenc}
\usepackage[russian, french, 
english]{babel}

\usepackage{graphicx}
\usepackage{wrapfig}
\usepackage{mathtools}
\usepackage{amsbsy}
\usepackage[inline]{enumitem}
\usepackage{mathrsfs}
\usepackage{array}
\usepackage{multicol}
\usepackage{stmaryrd}
\usepackage{cancel}
\usepackage{lmodern}
\usepackage{mathabx}
\usepackage{upgreek}
\usepackage{titlesec}
\usepackage{titletoc}
\usepackage[spacing=true,kerning=true,babel=true,tracking=true]{microtype}
\usepackage[shortcuts]{extdash}
\usepackage[foot]{amsaddr}
\usepackage[left=1in,right=1in,top=1in,bottom=1in,bindingoffset=0cm]{geometry}
\usepackage{bm}
\usepackage{centernot}
\usepackage{mdframed}
\usepackage[hidelinks]{hyperref}
\usepackage{xspace}
\usepackage[most]{tcolorbox}
\usepackage{framed}
\usepackage[justification=centering, labelfont=bf]{caption}
\usepackage{tikz}
\usetikzlibrary{math}
\usetikzlibrary{shapes,snakes}
\usetikzlibrary{arrows.meta}
\usetikzlibrary{decorations.pathmorphing}
\usetikzlibrary{patterns}
\usetikzlibrary{positioning,calc,intersections}
\usepackage{tikzsymbols}
\usepackage{float}

\DeclareFontFamily{U}{skulls}{}
\DeclareFontShape{U}{skulls}{m}{n}{ <-> skull }{}

\usepackage[
    sortcites,
    backend=biber, style=alphabetic, sorting=nyt, maxnames=100,backref=true]{biblatex}
\newtheoremstyle{bfnote}%
{}{}%
{\slshape}{}%
{\bfseries}{\bfseries.}%
{ }%
{\thmname{#1}\thmnumber{ #2}\thmnote{ \ep{\normalfont{}#3}}}

\theoremstyle{bfnote}
\newtheorem{theo}{Theorem}[section]
\newtheorem{theorem}[theo]{Theorem}
\newtheorem*{theo*}{Theorem}
\newtheorem{prop}[theo]{Proposition}
\newtheorem{lemma}[theo]{Lemma}
\newtheorem{claim}[theo]{Claim}

\newtheorem*{corl*}{Corollary}

\theoremstyle{definition}
\newtheorem{definition}[theo]{Definition}
\newtheorem*{defn*}{Definition}
\newtheorem{exmp}[theo]{Example}

\newtheorem{remk}[theo]{Remark}
\newtheorem{question}[theo]{Question}
\newtheorem{remks}[theo]{Remarks}
\newtheorem*{remks*}{Remarks}
\newtheorem*{exmp*}{Example}

\theoremstyle{remark}
\newtheorem*{ques*}{Question}
\newtheorem*{remk*}{Remark}

\newcommand*{\myproofname}{Proof}
\newenvironment{claimproof}[1][\myproofname]{\begin{proof}[#1]}{\end{proof}}

\makeatletter
\newcommand{\neutralize}[1]{\expandafter\let\csname c@#1\endcsname\count@}
\makeatother

\newcommand{\0}{\varnothing}
\newcommand{\set}[1]{\{#1\}}

\newcommand{\Q}{\mathbb{Q}}
\renewcommand{\epsilon}{\varepsilon}

\renewcommand{\phi}{\varphi}
\renewcommand{\theta}{\vartheta}
\renewcommand{\leq}{\leqslant}
\renewcommand{\geq}{\geqslant}
\newcommand{\defeq}{\coloneqq}

\newcommand{\bemph}[1]{{\normalfont#1}} 
\newcommand{\ep}[1]{\bemph{(}#1\bemph{)}} 

\newcommand{\emphd}[1]{{\fontseries{b}\selectfont\textsf{#1}}}
\newcommand{\acts}{\mathrel{\reflectbox{$\righttoleftarrow$}}}
\newcommand{\G}{\Gamma}
\newcommand{\pto}{\dashrightarrow}
\newcommand{\Free}{\mathcal{F}}

\newcommand{\rest}[2]{{{#1}\vert_{#2}}}

\newcommand{\dist}{\mathsf{dist}}
\newcommand{\fs}[2]{{[{#1}]^{#2}}}
\newcommand{\finset}[1]{\fs{{#1}}{{<\infty}}}
\newcommand{\LOCAL}{$\mathsf{LOCAL}$\xspace}
\newcommand{\N}{\mathbb{N}}
\newcommand{\Z}{\mathbb{Z}}
\newcommand{\R}{\mathbb{R}}
\newcommand{\F}{\mathcal{F}}

\newcommand{\dom}{\mathsf{dom}}

\newcommand{\inv}{^{-1}}

\newcommand{\CONT}{\mathtt{CONTINUOUS}}
\newcommand{\BOREL}{\mathtt{BOREL}}
\newcommand{\MEAS}{\mathtt{MEASURE}}
\newcommand{\BAIRE}{\mathtt{BAIRE}}
\newcommand{\FIID}{\mathtt{FIID}}
\newcommand{\FFIID}{\mathtt{FFIID}}
\newcommand{\COMP}{\mathtt{COMPUTABLE}}
\newcommand{\BS}{\mathtt{BaireSHIFT}}
\newcommand{\Shift}{\mathcal{S}}
\newcommand{\Leb}{\mathcal{L}}
\newcommand{\Cay}{\mathsf{Cay}}
\newcommand{\RBG}{\set{\mathtt{R}, \mathtt{B}, \mathtt{G}}}
\newcommand{\Waffle}{\mathsf{RT}}
\newcommand{\blank}{\mathtt{blank}}
\newcommand{\CW}{\mathsf{CRT}}
\newcommand{\Sn}{\set{\pm \mathsf{e}_1,\ldots,\pm \mathsf{e}_n}}

\newcommand{\POS}{{1}}
\newcommand{\NEG}{{0}}
\newcommand{\CWplus}{\CW^+}

\numberwithin{equation}{section}

\newenvironment{scproof}[1][]{\begin{proof}[\textsc{\upshape{Proof}}#1]}{\end{proof}}

\titleformat{\section}[block]{\large\bfseries\sffamily}{\thesection.}{1ex}{}
\titleformat{\subsection}[block]{\bfseries\sffamily}{\thesubsection.}{1ex}{}
\titleformat{\subsubsection}[block]{\itshape}{\bfseries\upshape\sffamily\thesubsubsection.}{1ex}{}

\titlespacing*{\section}{0pt}{*3}{*1}
\titlespacing*{\subsection}{0pt}{*3}{*1}
\titlespacing*{\subsubsection}{0pt}{*2}{*1}

 \titlecontents{section}
 [1.5em] %
 {\smallskip}
 {\bfseries\thecontentslabel\hspace{1.02em}}
{\bfseries}
 {\,\,\titlerule*[0.77pc]{}\bfseries\contentspage}
\titlecontents{subsection}
 [4em] %
 {\smallskip}
 {\thecontentslabel\hspace{1.02em}}
{\hspace*{2.32em}}
 {\,\,\titlerule*[0.77pc]{.}\contentspage}

\renewbibmacro{in:}{}

\renewbibmacro*{volume+number+eid}{%
	\printfield{volume}%
	\setunit*{\addnbspace}
	\printfield{number}%
	\setunit{\addcomma\space}%
	\printfield{eid}}

\DeclareFieldFormat[article]{volume}{\textbf{#1}\space}
\DeclareFieldFormat[article]{number}{\mkbibparens{#1}}

\DeclareFieldFormat{journaltitle}{#1,}
\DeclareFieldFormat[thesis]{title}{\mkbibemph{#1}\addperiod}
\DeclareFieldFormat[article, unpublished, thesis]{title}{\mkbibemph{#1},}
\DeclareFieldFormat[book]{title}{\mkbibemph{#1}\addperiod}
\DeclareFieldFormat[unpublished]{howpublished}{#1, }

\DeclareFieldFormat{pages}{#1}

\DeclareFieldFormat[article]{series}{Ser.~#1\addcomma}

\setlist{topsep=3pt,itemsep=3pt}

\title{\sffamily Separating complexity classes of LCL problems on grids}
\date{}

\author{Katalin~Berlow}
\address{\normalfont (KB, FW) Department of Mathematics, University of California, Berkeley, CA, USA}
\email{katalin@berkeley.edu}

\author{Anton~Bernshteyn}
\address{\normalfont (AB) Department of Mathematics, University of California, Los Angeles, CA, USA}
\email{bernshteyn@math.ucla.edu}

\author{Clark~Lyons}
\address{\normalfont (CL) 
Alfréd Rényi Institute of Mathematics, Budapest, Hungary}
\email{lyons.clark@renyi.hu}

\author{Felix~Weilacher}
\email{weilacher@berkeley.edu}

\thanks{KB's research is partially supported by the NSF Graduate
Research Fellowship Program under grant DGE-2146752.  AB's research is partially supported by the NSF CAREER grant DMS-2528522. CL is supported by the ERC Advanced Grant “ERMiD”. FW is supported by the NSF under award number DMS-2402064.}

\begin{document}

\maketitle

\begin{abstract}
    We study the complexity of locally checkable labeling (LCL) problems on $\Z^n$ from the point of view of descriptive set theory, computability theory, and factors of i.i.d. Our results separate various complexity classes that were not previously known to be distinct and serve as counterexamples to a number of natural conjectures in the field.
\end{abstract}

\section{Introduction}\label{sec:intro}

    \subsection{Descriptive complexity of LCL problems}

    A common theme throughout mathematics and computer science
    is the desire to understand, in rigorous terms, the complexity inherent in various problems of interest. Depending on the subject area, different notions of complexity can be proposed, resulting in different complexity hierarchies. In this paper we investigate the complexity of so-called locally checkable labeling problems on $\Z^n$ from the perspective of descriptive set theory, probability theory, and computability theory. Our contributions separate certain complexity classes that were not previously known to be distinct and serve as counterexamples to an array of natural conjectures. 
    Further results along these lines will be presented in the forthcoming paper \cite{Sequel}.
    
    To start, let us identify the class of problems we will be concerned with. In the definition given below, the underlying structure is provided by 
    a group $\G$ (typically assumed to be finitely generated); in our main results, we will have
    $\G = \Z^n$ for $n \geq 2$. This definition can also be phrased, somewhat more generally, in graph-theoretic terms, but we will focus on groups for convenience. Extensions of our results beyond the group-theoretic framework are 
    discussed in \S\ref{subsec:unlabeled}.

    \begin{definition}[LCL problems on a group]
        A \emphd{locally checkable labeling} \ep{\emphd{LCL} for short} \emphd{problem} on a group $\G$ is a triple $\Pi = (W, \Lambda, \mathcal{A})$, where:
        \begin{itemize}
            \item $W \subseteq \G$ is a finite set, called a \emphd{window},
            \item $\Lambda$ is a finite set whose elements are called \emphd{labels}, and
            \item $\mathcal{A} \subseteq \Lambda^W$ is a family of mappings $W \to \Lambda$, called \emphd{allowed configurations}.
        \end{itemize}
        %
        A \emphd{solution} to $\Pi$, or a \emphd{$\Pi$-labeling} of $\G$, is a function $f \colon \G \to \Lambda$ such that for all $\gamma \in \G$, the mapping \[W \to \Lambda \ \colon \ w \mapsto f(w\gamma)\] is in $\mathcal{A}$. The set of all $\Pi$-labelings of $\G$ is called a \emphd{subshift of finite type} (an \emphd{SFT} for short).
    \end{definition}

    The term ``locally checkable labeling problem'' was introduced by Naor and Stockmeyer \cite{NaorStock} in a computer science context. At the same time, subshifts of finite type are studied 
    in topological and symbolic dynamics, where they form an important class of dynamical systems \cite{subshifts1,subshifts2}. 

    \begin{exmp}[Proper $k$-colorings of Cayley graphs]\label{exmp:coloring}
        Let $\G = \langle S \rangle$ be a group generated by a finite 
        set $S$ with $\bm{1} \notin S$. The 
        \emphd{Cayley graph} $\Cay(\G, S)$ has vertex set $\G$ and edge set \[E(\Cay(\G,S)) \,\defeq\, \set{\set{\gamma, \sigma \gamma} \,:\, \gamma \in \G,\, \sigma \in S}.\] 
        Following a standard set-theoretic convention, we identify each natural number $k \in \N$ with the $k$-element set $\set{0,1,\ldots, k-1}$. Given $k \in \N$, a \emphd{proper $k$-coloring} of $\Cay(\G, S)$ is a function $f \colon \G \to k$ 
        such that $f(\gamma) \neq f(\delta)$ whenever $\gamma$ and $\delta$ are adjacent in $\Cay(\G, S)$. Note that $f$ is a proper $k$-coloring of $\Cay(\G,S)$ if and only if $f$ is a solution to the LCL problem $\Pi = (W,\Lambda, \mathcal{A})$, where 
        \[
            W \,\defeq\, \set{\bm{1}} \cup S, \qquad \Lambda \,\defeq\, k, \quad \text{and} \quad \mathcal{A} \,\defeq\, \big\{\phi \in \Lambda^W \,:\, \phi(\bm{1}) \neq \phi(\sigma) \text{ for all } \sigma \in S\big\}.
        \]
        In particular, the set of all proper $k$-colorings of $\Cay(\G,S)$ is an SFT.
    \end{exmp} 

    How difficult is it to solve a given LCL problem $\Pi$? This question can be interpreted in a number of ways. For instance, it can be studied from the standpoint of \emph{algorithmic decidability}; that is, one may seek an algorithm that, given an LCL problem $\Pi$, determines whether $\Pi$ has a solution. Such an algorithm exists for LCL problems on $\Z$ and, more generally, on finitely generated free groups \cite[Prop.~9.3.35]{ABJ_survey}. On the other hand, refuting a conjecture of Wang \cite{Wang}, Berger showed that there is no such algorithm for LCL problems on $\Z^2$ \cite{Berger}; see also \cite{Robinson} for a simplified proof due to Robinson. Considering other finitely generated groups gives rise to a rich theory and an active research area; see the survey \cite{ABJ_survey} by Aubrun, Barbieri, and Jeandel for further references.

    Another approach to the complexity of LCL problems is provided by \emph{distributed computing}, the area of computer science concerned with decentralized computational networks. For an introduction to this field, see the book \cite{BE} by Barenboim and Elkin. Roughly speaking, a distributed algorithm for an LCL problem $\Pi$ is a communication protocol that allows autonomous processes placed at the nodes of an underlying graph $G$ to compute their corresponding labels forming a solution to $\Pi$ (in our setting $G$ is the Cayley graph of $\G$). Starting with the seminal work of Naor and Stockmeyer \cite{NaorStock}, there has been a concerted effort to obtain a classification of LCL problems according to their distributed complexity. In particular, a series of contributions by a large number of researchers has culminated in a complete characterization of possible complexities of LCL problems on $\Z^n$ \cite{NaorStock,CS2,CS3,CS8,CS7} and the free groups $\mathbb{F}_n$ \cite{NaorStock,CS1,CS2,CS3,CS4,CS5,CS6,CS7} in the so-called \LOCAL model of distributed computation.
    
    In this paper we study the complexity theory for LCL problems developed using concepts from \emph{descriptive set theory}. The idea of employing descriptive set-theoretic notions in a combinatorial context dates back to the pioneering work of Kechris, Solecki, and Todorcevic \cite{KST}. By now it has grown into a fruitful field---called \emph{descriptive combinatorics}---with many connections to other areas; see \cite{KechrisMarks} for a survey by Kechris and Marks, \cite{Pikh_survey} for an introductory article by Pikhurko, and \cite{Notices} for an overview of some recent developments. We also direct the reader to \cite{KechrisDST,AnushDST} for general descriptive set theory background. In descriptive combinatorics, 
    instead of simply
    solving a given LCL problem $\Pi$ on a group $\G$, we work with an {action} $\G \acts X$ of $\G$, and our goal is to solve $\Pi$ on every orbit of this action in a ``uniform'' 
    fashion.

    \begin{definition}[Solutions to LCL problems on group actions]
        Let $\G \acts X$ be an action of a group $\G$. 
        For an LCL problem $\Pi = (W, \Lambda, \mathcal{A})$ on $\G$, a \emphd{solution} to $\Pi$ on $X$, or a \emphd{$\Pi$-labeling} of $X$, is a function $f \colon X \to \Lambda$ such that for all $x \in X$, the mapping $W \to \Lambda \colon w \mapsto f(w \cdot x)$ is in $\mathcal{A}$.  
    \end{definition}

    Given an action $\G \acts X$, any map $f \colon X \to \Lambda$ gives rise to a family of labelings of $\G$ parameterized by the points of $X$. Namely, for each point $x \in X$, we have the corresponding mapping $f_x \colon \G \to \Lambda$ defined by $f_x(\gamma) \defeq f(\gamma \cdot x)$ for all $\gamma \in \G$. 
    It is then clear that $f \colon X \to \Lambda$ is a solution to an LCL problem $\Pi$ on $X$ if and only if for all $x \in X$, $f_x \colon \G \to \Lambda$ is a solution to $\Pi$ on $\G$.

    Recall that a group action $\G \acts X$ is \emphd{free} if 
    $\gamma \cdot x \neq x$ for all $x \in X$ and $\gamma \in \G \setminus \set{\bm{1}}$. 
    It is not hard to see that if $\G \acts X$ is a {free} 
    action, 
    then any LCL problem $\Pi = (W, \Lambda, \mathcal{A})$ with a solution on $\G$ also has a solution on $X$. 
    Indeed, given a $\Pi$-labeling $f \colon \G \to \Lambda$ 
    of $\G$, we can obtain a $\Pi$-labeling $f^* \colon X \to \Lambda$ of $X$ as follows: Form a set $X^* \subseteq X$ consisting of one point from each $\G$-orbit and then define $f^*(\gamma \cdot x) \defeq f(\gamma)$ for all $\gamma \in \G$ and $x \in X^*$. 
    Note, however, that this construction invokes the Axiom of Choice to establish the existence of the set $X^*$. In particular, there is no reason to expect the function 
    $f^*$ 
    to have any desirable topological or measure-theoretic
    regularity properties.
    
    \begin{exmp}[Proper $2$-coloring on $\Z$]\label{exmp:2coloring}
        Consider the group $\Z$. The Cayley graph $G \defeq \Cay(\Z, \set{1})$ 
        is a path which is infinite in both directions. Let $\Pi$ be the proper $2$-coloring problem for $G$ 
        (as described in Example~\ref{exmp:coloring}). It is then easy to see that $\Z$ has a $\Pi$-labeling, namely the map $\Z \to 2$ sending each $k \in \Z$ to $0$ if $k$ is even and $1$ if $k$ is odd. (The only other $\Pi$-labeling of $\Z$ sends $k$ to $0$ if $k$ is odd and $1$ if $k$ is even.) 
        It follows that every free action $\Z \acts X$ admits a $\Pi$-labeling $f \colon X \to 2$. However, such a $\Pi$-labeling $f$ often must exhibit various pathological properties.
        
        For example, suppose that $T \colon X \to X$ is an aperiodic invertible measure-preserving transformation on a standard probability space $(X,\mu)$ such that $T^2$ is ergodic. 
        (There are many examples of such transformations, for instance the rotation of a circle by the angle $\alpha \pi$ for fixed $\alpha \in \R \setminus \Q$.) Iterating $T$ gives rise to a free action $\Z \acts (X,\mu)$ defined by $k \cdot x \defeq T^k(x)$ for all $k \in \Z$ and $x \in X$. 
        We claim that, although $X$ does have a $\Pi$-labeling, no $\Pi$-labeling of $X$ can be {$\mu$-measurable}. 
%
            Indeed, if $f \colon X \to 2$ were a $\mu$-measurable $\Pi$-labeling of $X$, then $A \defeq f^{-1}(0)$ and $B \defeq f^{-1}(1)$ would be $\mu$-measurable sets that partition $X$ and satisfy $T(A) = B$, hence $\mu(A) = \mu(B) = 1/2$. At the same time, $T^2(A) = A$, 
        so by the ergodicity of $T^2$ we must have $\mu(A) \in \set{0,1}$, a contradiction. 
    \end{exmp}
    
    By adding extra regularity constraints
    , we arrive at the following definition, which introduces the ``big four'' 
    complexity classes studied in descriptive combinatorics:

    \begin{tcolorbox}[colframe=white]
    \begin{definition}[Descriptive complexity classes]\label{defn:complexity}
        Let $\G$ be a finitely generated group and let 
        $\Pi = (W, \Lambda, \mathcal{A})$ be an LCL problem on $\G$. We equip the set $\Lambda$ with the discrete topology. The classes $\BOREL(\G)$, $\MEAS(\G)$, $\BAIRE(\G)$, and $\CONT(\G)$ are defined as follows:


        \begin{tcolorbox}

        $\Pi \in \BOREL(\G)$ if and only if 
            every free Borel action $\G \acts X$ on a standard Borel space $X$ has a Borel $\Pi$-labeling $f \colon X \to \Lambda$.

        \medskip

        $\Pi \in \MEAS(\G)$ if and only if 
            every 
            free Borel action $\G \acts (X,\mu)$ on a standard probability space $(X,\mu)$ has a $\mu$\=/measurable $\Pi$-labeling $f \colon X \to \Lambda$.

            \medskip

        $\Pi \in \BAIRE(\G)$ if and only if 
            every free Borel action $\G \acts X$ on a Polish space $X$ has a Baire-measurable $\Pi$-labeling $f \colon X \to \Lambda$.

            \medskip

        $\Pi \in \CONT(\G)$ if and only if 
            every free continuous action $\G \acts X$ on a zero\-/dimensional 
            Polish space $X$ has a continuous $\Pi$-labeling $f \colon X \to \Lambda$.
            \end{tcolorbox}
    \end{definition}
    \end{tcolorbox}

    \begin{remks}\label{remk:classes}
        \begin{enumerate}[label=\ep{\normalfont\roman*},wide]
            \item Since the topology on $\Lambda$ is discrete, a labeling $f \colon X \to \Lambda$ is Borel, $\mu$\-/measurable, Baire-measurable, or continuous if and only if for every label $\lambda \in \Lambda$, its preimage $f^{-1}(\lambda)$ is a Borel, $\mu$-measurable, Baire-measurable, or clopen subset of $X$ respectively.
            
            \item In the definition of $\BOREL(\G)$, $X$ can be assumed to be any particular uncountable standard Borel space, for example, the unit interval $[0,1]$, the real line $\R$, or the Cantor space $2^\N$. This results in no loss of generality due to the Borel isomorphism theorem \cite[Thm.~15.6]{KechrisDST}.

            \item\label{item:measure} We stress that in the definition of $\MEAS(\G)$, the action $\G \acts (X,\mu)$ need not be measure\-/preserving. (However, see Question~\ref{ques:pmp}.)
            

            \item Recall that a topological space is \emphd{zero-dimensional} if it has a base consisting of clopen sets. The assumption that $X$ is zero-dimensional in the definition of $\CONT(\G)$ guarantees that the family of continuous functions $X \to \Lambda$ is sufficiently rich to make the definition nontrivial.
        \end{enumerate}
    \end{remks}

    The following relations between these classes hold for every finitely generated group $\G$:
    \[
        \CONT(\G) \,\subseteq\, \BOREL(\G) \,\subseteq\,\MEAS(\G) \cap \BAIRE(\G).
    \]
    The inclusion $\BOREL(\G) \subseteq \MEAS(\G) \cap \BAIRE(\G)$ is immediate since every Borel function on a standard probability space (resp.~on a Polish space) is measurable (resp.~Baire\-/measurable). The inclusion $\CONT(\G) \subseteq \BOREL(\G)$ follows from standard results in descriptive set theory \cite[\S13]{KechrisDST}, which imply that if $\G \acts X$ is a Borel action on a standard Borel space $X$, then there exists a compatible zero-dimensional Polish topology $\tau$ on $X$ making the action continuous.

    More is known for specific groups $\G$. For example, when $\G = \Z$ the four classes coincide: 
    \[
        \CONT(\Z) \,=\, \BOREL(\Z) \,=\, \MEAS(\Z) \,=\, \BAIRE(\Z).
    \]
    See \cite{paths} by Greb\'ik and Rozho\v{n} for a thorough study of this case. On the other hand, for a free group $\mathbb{F}_n$ of rank $n \geq 2$, the four classes are distinct and linearly ordered by inclusion:
    \[
        \CONT(\mathbb{F}_n) \,\subsetneq\, \BOREL(\mathbb{F}_n) \,\subsetneq\, \MEAS(\mathbb{F}_n) \,\subsetneq\, \BAIRE(\mathbb{F}_n).
    \]
    This chain of strict inclusions is established in \cite{trees} by Brandt, Chang, Greb\'ik, Grunau, Rozho\v{n}, and Vidny\'anszky,\footnote{Technically, the results in \cite{trees} are stated for LCL problems on unlabeled regular trees, but virtually the same arguments apply in the setting of $\mathbb{F}_n$-actions.} building on and synthesizing numerous earlier results in descriptive combinatorics, distributed computing, ergodic theory, and probability theory 
    \cite{BerDist,Ber_cont,Marks,CS3,entropy,ConleyKechris,MeasurableBrooks,CM16}. 
    For sufficiently large $n$, the three classes $\BOREL(\mathbb{F}_n)$, $\MEAS(\mathbb{F}_n)$, and $\BAIRE(\mathbb{F}_n)$ are distinguished by the proper $k$-coloring problem (see Example \ref{exmp:coloring}), since the smallest $k$ such that the proper $k$-coloring problem for the standard Cayley graph of $\mathbb{F}_n$ belongs to a given complexity class is $2n+1$ for the class $\BOREL(\mathbb{F}_n)$ \cite{KST,Marks}, 
    $(2+o(1))n/\log n$ for the class $\MEAS(\mathbb{F}_n)$ \cite{Ber_Lebesgue,BerDist,RV}, 
    and $3$ for the class $\BAIRE(\mathbb{F}_n)$ \cite{CM16}. 

    In this paper we are interested in another natural class of groups, namely $\Z^n$ for $n \geq 2$. In contrast to the $\Z$ case, for $n \geq 2$ we have the following strict inclusion:
    \[
        \CONT(\Z^n) \,\subsetneq\, \BOREL(\Z^n),
    \]
    because the proper $3$-coloring problem for the standard Cayley graph of $\Z^n$ belongs to the class $\BOREL(\Z^n)$ but not $\CONT(\Z^n)$, as shown by Gao, Jackson, Krohne, and Seward \cite{BorelAbelian,ContinuousAbelian}. In general, the class $\CONT(\Z^n)$ is relatively well understood, thanks to the work of 
    Gao and Jackson \cite{GaoJack}, Gao, Jackson, Krohne, and Seward \cite{ContinuousAbelian}, Greb\'ik and Rozho\v{n} \cite{GRgrids}, and the second named author \cite{Ber_cont}. For example, it is known that $\CONT(\Z^n)$ is precisely the class of LCL problems on $\Z^n$ solvable by a sublogarithmic-time algorithm in the \LOCAL model of distributed computation. This is a general result that holds for arbitrary countable groups \cite[Thm.~1.15]{Ber_cont}; for a detailed discussion of the $\Z^n$ case, see \cite[\S7.1]{GRgrids}.

    Much less is known about the classes $\BOREL(\Z^n)$, $\MEAS(\Z^n)$, and $\BAIRE(\Z^n)$. In particular, it has been a major open problem whether any of them are distinct (see, e.g., \cites[\S6]{paths}[\S8]{GRgrids}):

    \begin{question}\label{ques:separation}
        For $n \geq 2$, are any of $\BOREL(\Z^n)$, $\MEAS(\Z^n)$, and $\BAIRE(\Z^n)$ distinct?
    \end{question}

    Here we partially answer Question~\ref{ques:separation} by exhibiting an LCL problem on $\Z^n$ for each $n \geq 2$ that belongs to $\MEAS(\Z^n)$ but not $\BAIRE(\Z^n)$:

    \begin{tcolorbox}
	    \begin{theorem}\label{theo:separation}
	        For every integer $n \geq 2$, we have \[\MEAS(\Z^n) \setminus \BAIRE(\Z^n) \, \neq \, \0.\] As a consequence, $\BAIRE(\Z^n) \neq \MEAS(\Z^n)$  and $\BOREL(\Z^n) \neq \MEAS(\Z^n)$.
	    \end{theorem}
	\end{tcolorbox}

    It remains an open problem to determine whether $\BOREL(\Z^n) = \BAIRE(\Z^n)$.

    In addition to making progress on Question~\ref{ques:separation}, Theorem~\ref{theo:separation} refutes a piece of conventional wisdom in descriptive combinatorics. Although measure and Baire category may behave quite differently \cite{Oxtoby}, in practice it is usually easier to solve an LCL problem Baire-measurably rather than $\mu$-measurably, and a natural expectation is that $\MEAS(\G) \subseteq \BAIRE(\G)$ for every finitely generated group $\G$. While this is true when $\G$ is a free group, and no counterexamples have been known before, Theorem~\ref{theo:separation} shows that this inclusion fails for $\G = \Z^n$ with $n \geq 2$.

    \subsection{Shift actions: the universal domain for descriptive combinatorics?}\label{subsec:uni}

    While the complexity classes described in Definition~\ref{defn:complexity} contain LCL problems that can be solved on \emph{all} free actions of $\G$ with specified regularity properties, it also makes sense to restrict one's attention to particular actions of special interest. One such action is the {shift}:

    \begin{definition}[Shift action]
        Let $\G$ be a group and let $X$ be a set. We use $\Shift(X,\G) \defeq X^\G$ to denote the set of all functions $x \colon \G \to X$. The \emphd{shift action} $\G \acts \Shift(X,\G)$ is given by the formula
    \[
        (\gamma \cdot x)(\delta) \,\defeq\, x(\delta\gamma) \quad \text{for all $\gamma$, $\delta \in \G$ and $x \in \Shift(X,\G)$}.
    \]
    \end{definition}

    We will only use the above definition for countable (in fact, finitely generated) groups $\G$. If $\G$ is countable and $X$ is a standard Borel, (zero-dimensional) Polish, or standard probability space, then so is $\Shift(X,\G)$ equipped with a suitable product structure, and the shift action $\G \acts \Shift(X,\G)$ is Borel (as well as continuous and measure-preserving when these make sense). 
    
    
    In general, the shift action $\G \acts \Shift(X,\G)$ is not free, so it does not directly fit into the context provided by Definition~\ref{defn:complexity}. To remedy this, we will only work on the \emphd{free part} of the shift action, i.e., the set $\Free(X,\G)$ of all $x \in \Shift(X,\G)$ such that $\gamma \cdot x \neq x$ for every $\gamma \in \G \setminus \set{\bm{1}}$. For a countably infinite group $\G$ and a Polish space $X$ with $|X| \geq 2$, the set $\Free(X,\G)$ is comeager in $\Shift(X,\G)$. Similarly, if $\G$ is countably infinite and $(X,\mu)$ is a standard probability space, then $\Free(X,\G)$ is conull in $\Shift(X,\G)$ (with respect to the product measure $\mu^\G$) unless $\mu$ is concentrated on a single point. Therefore, in these cases there is no difference from the point of view of Baire category and measure theory whether one considers the whole space $\Shift(X,\G)$ or just $\Free(X,\G)$.

    The following result shows that shift actions are in some sense ``universal'' for Borel and continuous solutions to LCL problems:

    \begin{theorem}[{Seward--Tucker-Drob \cite{STD}, AB/Seward \cite[Thm.~1.12]{Ber_cont}}]\label{theo:universal}
        Fix a finitely generated group $\G$. In the following statements, $\Pi$ is an LCL problem on $\G$.

        \begin{enumerate}[label=\ep{\normalfont\roman*}]
            \item\label{item:Borel_universal} For every standard Borel space $X$ of cardinality at least $2$, we have
            \[
                \BOREL(\G) \,=\, \set{\Pi \,:\, \text{$\Pi$ has a Borel solution on $\Free(X,\G)$}}.
            \]

            \item For every zero-dimensional Polish space $X$ of cardinality at least $2$, we have
            \[
                \CONT(\G) \,=\, \set{\Pi \,:\, \text{$\Pi$ has a continuous solution on $\Free(X,\G)$}}.
            \]
         \end{enumerate}
    \end{theorem}

    \begin{remk}
        In the papers \cite{STD,Ber_cont} the above results are established in the case when $X$ has exactly $2$ points; the generalization to arbitrary $X$ of cardinality at least $2$ follows since any injection $2 \to X$ gives rise to a continuous $\G$-equivariant map $\Free(2,\Gamma) \to \Free(X,\Gamma)$. (Recall that $2$ denotes the set $\set{0,1}$, which we endow with the discrete topology.) 
    \end{remk}

    A natural question is whether there is a version of Theorem~\ref{theo:universal} for the classes $\MEAS(\G)$ and $\BAIRE(\G)$. For $\BAIRE(\G)$, some partial results are known. By a theorem of Keane and Weiss \cite[Thm.~2]{Weiss_generic}, for arbitrary Polish spaces $X$, $Y$ of cardinality at least $2$, the shift actions $\G \acts \Free(X,\G)$ and $\G \acts \Free(Y,\G)$ are isomorphic modulo meager sets, which implies that they admit Baire-measurable solutions to exactly the same LCL problems.\footnote{The statement in \cite{Weiss_generic} is for zero-dimensional spaces only. The result for general Polish spaces follows since every Polish space contains a comeager zero-dimensional subspace.} Thus, we may without loss of generality restrict our attention to a specific choice of a Polish space $X$, say $X = 2$, and give the following definition:

    \begin{tcolorbox}[colframe=white]
    \begin{definition}[$\BS$]\label{defn:BS}
        Fix a finitely generated group $\G$. 
        An LCL problem $\Pi$ on $\G$ is in the class $\BS(\G)$ if and only if it has a Baire-measurable solution on $\Free(2,\G)$.
    \end{definition}
    \end{tcolorbox}
    
    When $\G$ is a free group, a version of Theorem~\ref{theo:universal} for the class $\BAIRE(\G)$ does hold:

    \begin{theorem}[{\cite[Prop.~6.1, Prop.~6.4, and proof of Thm.~6.6]{trees}\footnote{Technically, the proof of \cite[Thm.~6.6]{trees} is presented for the group $(\Z/2\Z)^{* n}$ rather than $\mathbb{F}_n$, but exactly the same arguments apply in the $\mathbb{F}_n$ case. (The only relevant feature of the group is that its Cayley graph is a tree.) }}]
        For all $n \in \N$, we have 
        \[
            \BAIRE(\mathbb{F}_n) \,=\, \BS(\mathbb{F}_n). 
        \]
    \end{theorem}

    On the other hand, we show that an analogous statement for the group $\Z^n$ with $n \geq 2$ fails; this is the first known instance where the shift action is not ``universal'' for Baire-measurable solutions to LCL problems (on any finitely generated group, not just $\Z^n$).

    \begin{tcolorbox}
	    \begin{theorem}\label{theo:not_universal}
	        For every integer $n \geq 2$, we have $\BAIRE(\Z^n) \subsetneq \BS(\Z^n)$. 


         %
	    \end{theorem}
	\end{tcolorbox}


    The role of shift actions in the context of $\MEAS(\Gamma)$ 
    is the topic of the next subsection. 

    \subsection{\ep{Finitary} factors of i.i.d.}

    We view the unit interval $[0,1]$ as a probability space with the usual Lebesgue measure. (This is just a concrete example of an atomless standard probability space, which are all isomorphic \cite[Thm.~17.41]{KechrisDST}.) The study of measurable functions on $\Shift([0,1],\G)$ (with respect to the product measure) is an important topic in probability theory and ergodic theory, known as the \emphd{theory of factors of i.i.d.~processes}\footnote{Here ``i.i.d.'' stands for ``independent and identically distributed,'' referring to the random variables $x(\gamma)$ for $\gamma \in \G$ and $x$ drawn at random from $\Shift([0,1], \G)$.}. A significant place in this theory is occupied by so-called finitary functions. 
    Since we prefer to work with free actions, we define finitary functions on $\Free([0,1],\G)$.

    \begin{definition}[Finitary functions]\label{defn:finitary}
        Let $\G$ be a finitely generated group 
        and let $\Lambda$ be a finite set. We say that two elements $x$, $x' \in \Free([0,1],\G)$ \emphd{agree} on a subset $D \subseteq \G$ if $x(\gamma) = x'(\gamma)$ for all $\gamma \in D$. A measurable function $f \colon \Free([0,1],\G) \to \Lambda$ is \emphd{finitary} if for almost every $x \in \Free([0,1],\G)$, there exists a finite set $D_x \subseteq \G$ such that $f(x) = f(x')$ whenever $x' \in \Free([0,1],\G)$ agrees with $x$ on $D_x$.
    \end{definition}

    Holroyd, Schramm, and Wilson \cite{Holroyd} initiated the study of LCL problems that can be solved by finitary functions. This line of research was continued in the works of Holroyd \cite{Holroyd1}, Brandt, Hirvonen, Korhonen, Lempi\"ainen, \"Osterg\aa{}rd, Purcell, Rybicki, Suomela, and Uzna\'nsk \cite{CS8}, Spinka \cite{Spinka}, Brandt, Chang, Greb\'ik, Grunau, Rozho\v{n}, and Vidny\'anszky \cite{trees}, Greb\'ik and Rozho\v{n} \cite{GRgrids}, and Bencs, Hru{\v{s}}kov{\'{a}}, and T{\'{o}}th \cite{FIIDSchreier}. Following Greb\'ik and Rozho\v{n} \cite{GRgrids}, we define the following two complexity classes:

    \begin{tcolorbox}[colframe=white]
    \begin{definition}[\texttt{FIID} and \texttt{FFIID}]
        Let $\G$ be a finitely generated group and let 
        $\Pi = (W, \Lambda, \mathcal{A})$ be an LCL problem on $\G$. The classes $\FIID(\G)$ and $\FFIID(\G)$ are defined as follows:


        \begin{tcolorbox}
        $\Pi \in \FIID(\G)$  if and only if $\Pi$ has a measurable solution on $\Free([0,1],\G)$.

        \medskip

        $\Pi \in \FFIID(\G)$  if and only if $\Pi$ has a finitary solution on $\Free([0,1],\G)$.
        \end{tcolorbox}

        
        (Here ``\texttt{FIID}'' stands for ``factor of i.i.d.,'' and ``\texttt{FFIID}'' for ``finitary factor of i.i.d.'')
        
    \end{definition}
    \end{tcolorbox}

    Clearly, 
    $\MEAS(\G) \subseteq \FIID(\G)$ and $\FFIID(\G) \subseteq \FIID(\G)$ for every finitely generated group $\G$. 
    A natural question is whether any of these three classes can be different (see, e.g., \cites[\S7]{trees}[Rmk.~2.11]{GRgrids}). For instance, for $\G = \Z$, they coincide \cite{paths}:
    \[
        \MEAS(\Z) \,=\, \FIID(\Z) \,=\, \FFIID(\Z).
    \]
    On the other hand, there exist finitely generated groups $\G$ with $\FIID(\G) \neq \MEAS(\G)$; in other words, for such $\G$, there is an LCL problem $\Pi$ such that $\Pi \in \FIID(\G)$ but $\Pi \notin \MEAS(\G)$ \ep{notice the analogy with Theorem~\ref{theo:not_universal}}. This fact was recently observed by Kun for the group $\G = (\Z/2\Z)^{\ast d}$ with $d \geq 3$ in his breakthrough paper \cite{KunHall}. Specifically, Kun proved that the so-called balanced orientation problem is not in $\MEAS((\Z/2\Z)^{\ast d})$ \cite[Corl.~1.6]{KunHall}; on the other hand, it is in $\FIID((\Z/2\Z)^{\ast d})$ by a result of Thornton \cite{RileyOrient}. An alternative construction is to take $\G = (\Z/d\Z) \ast (\Z/d\Z)$ for $d \geq 3$ and let $\Pi$ be the problem asking to assign labels $0$ and $1$ so that exactly one point from each $(\Z/d\Z)$-orbit gets a~$1$. It is not hard to see that this is equivalent to finding a perfect matching in an auxiliary $d$-regular acyclic graph. As a consequence, $\Pi \in \FIID((\Z/d\Z) \ast (\Z/d\Z))$ by a theorem of Lyons and Nazarov \cite{LyonsNazarov} and $\Pi \notin \MEAS((\Z/d\Z) \ast (\Z/d\Z))$ due to \cite[Thm.~1.2]{KunHall}. In view of these constructions, it seems very plausible that $\FIID(\mathbb{F}_n) \neq \MEAS(\mathbb{F}_n)$ 
    for any $n \geq 2$, but this is still an open problem. 
    Similarly, it is not known whether $\FIID(\Z^n) = \MEAS(\Z^n)$ for $n \geq 2$. 

    As for the relationship between $\FIID(\G)$ and $\FFIID(\G)$, it has been poorly understood prior to this work; indeed, the following has been an open question \cites[430]{Oberwolfach}[Remk.~2.11]{GRgrids}:
     
    \begin{question}\label{ques:ffiid}
        Do we have $\FIID(\G) = \FFIID(\G)$  
        for every finitely generated group $\G$?
    \end{question}

    The class $\FFIID(\G)$ can be further stratified by measuring the sizes of the sets $D_x$ appearing in Definition~\ref{defn:finitary} \cite{Holroyd}, and this stratification largely parallels the complexity hierarchy arising in distributed computing \cite{GRgrids,trees}. Thus, a positive answer to Question~\ref{ques:ffiid} would establish a close link between distributed algorithms on the one hand and ergodic theory and probability theory on the other. We remark that there exist measurable functions on $\Free([0,1],\G)$ that are not finitary, but this on its own does not solve Question~\ref{ques:ffiid} since an LCL problem with a non-finitary measurable solution may also have a different solution that is finitary.
    
    Here we give a negative answer to Question~\ref{ques:ffiid} for $\G = \Z^n$ with $n \geq 2$. In fact, we show that the class $\FFIID(\Z^n)$ is distinct from both $\FIID(\Z^n)$ and $\MEAS(\Z^n)$:

    \begin{tcolorbox}
	    \begin{theorem}\label{theo:ffiid}
	        For every integer $n \geq 2$, we have \[\MEAS(\Z^n) \setminus \FFIID(\Z^n) \, \neq \, \0.\] As a consequence, $\FFIID(\Z^n) \neq \MEAS(\Z^n)$  and $\FFIID(\Z^n) \neq \FIID(\Z^n)$.
	    \end{theorem}
	\end{tcolorbox}

    As mentioned above, we still do not know whether $\FIID(\Z^n) = \MEAS(\Z^n)$.
    
    Another intriguing open question concerns the relationship between the classes $\FFIID(\Z^n)$ and $\BAIRE(\Z^n)$. Currently, neither of the inclusions $\FFIID(\Z^n) \subseteq \BAIRE(\Z^n)$, $\BAIRE(\Z^n) \subseteq \FFIID(\Z^n)$ is ruled out. That being said, we do have the following:

    \begin{tcolorbox}
	    \begin{theorem}\label{theo:FFIID_Baire}
	        For every finitely generated group $\G$, $\FFIID(\G) \subseteq \BS(\G)$.
	    \end{theorem}
	\end{tcolorbox}

    The proof of Theorem~\ref{theo:FFIID_Baire} relies on 
    the general analysis of the Baire-measurable combinatorics of $\Free(2,\G)$ performed in \cite{BerGeneric} by the second named author. 
    In the case $\G = \Z^n$ with $n \geq 2$, the inclusion $\FFIID(\Z^n) \subsetneq \BS(\Z^n)$ is strict (see Theorem~\ref{theo:main} below).


    

    \begin{remk}
        The reader may be wondering why we only define \emph{measurable} finitary functions and do not introduce versions of Definition~\ref{defn:finitary} in the Borel and Baire-measurable contexts. The reason is that such an extension of Definition~\ref{defn:finitary} would yield no new complexity classes. Indeed, let us for simplicity consider the space $\Free(2, \G)$ in place of $\Free([0,1], \G)$, and say that a Borel (or Baire-measurable) function $f \colon \Free(2, \G) \to \Lambda$ is \emphd{Borel-finitary} (resp.~\emphd{Baire-finitary}) if for every point (resp.~for a comeager set of points) $x \in \Free(2,\G)$, there exists a finite set $D_x \subseteq \G$ such that $f(x) = f(x')$ whenever $x' \in \F(2,\G)$ agrees with $x$ on $D_x$. It is not hard to see that a function $f \colon \Free(2, \G) \to \Lambda$ is Borel-finitary if and only if it is continuous \ep{this is essentially the definition of the product topology on $\Shift(2, \G)$}, and every Baire-measurable function is equal to some continuous, hence
        Baire-finitary function on a comeager set \cite[Prop.~9.8]{AnushDST}. Note that, in particular, this observation yields the inclusion $\CONT(\G) \subseteq \FFIID(\G)$.
    \end{remk}

    \subsection{Computable combinatorics}

    Another way to gauge the complexity of LCL problems is by using computability theory instead of descriptive set theory. This avenue of research has a long history dating back to the 1970s and 1980s; see the survey \cite{Gasarch} by Gasarch for an overview. 
    Qian and the fourth named author recently discovered parallels between the computable setting and descriptive combinatorics \cite{ASIalgorithms}, making computable combinatorics part of the general theory of local complexity.
    In the following definition, we say that an action of a 
    group $\G \acts X$ on a computable set $X \subseteq \N$ is \emphd{computable} if 
    for each $\gamma \in \G$, the function $X \to X : x \mapsto \gamma \cdot x$ is computable.  

    \begin{tcolorbox}[colframe=white]
    \begin{definition}[Computable complexity]\label{defn:comp_comp}
        Fix a finitely generated group $\G$. 
        An LCL problem $\Pi = (W, \Lambda, \mathcal{A})$ on $\G$ belongs to the class $\COMP(\G)$ if and only if every free computable action $\G \acts X$ on a computable set $X \subseteq \N$ has a computable $\Pi$-labeling $f \colon X \to \Lambda$.
%
    \end{definition}
    \end{tcolorbox}

    Note that if a finitely generated group $\G$ admits a free computable action $\G \acts X$ and $X \neq \0$, then the word problem on $\G$ is decidable: To decide whether $\sigma_1 \cdots \sigma_k = \mathbf{1}$ for a tuple of generators $(\sigma_1, \ldots, \sigma_k)$, we take any $x \in X$, compute $(\sigma_1 \cdots \sigma_k) \cdot x$, and compare the result with $x$. 
    Therefore, if the word problem on $\G$ is undecidable, the class $\COMP(\G)$ vacuously includes all LCL problems on $\G$. In other words, Definition~\ref{defn:comp_comp} is only interesting when $\G$ has a decidable word problem, and we shall make this assumption throughout this section. 
    
    Although Definitions~\ref{defn:complexity} and \ref{defn:comp_comp} consider rather different types of group actions (Borel vs.~computable), it turns out that every finitely generated group $\G$ 
    satisfies the following inclusion:
    \[
        \CONT(\G) \,\subseteq\, \COMP(\G).
    \]
    This is a consequence of \cite[Thm.~1.4.1]{WeilacherThesis} and \cite[Thm.~1.15]{Ber_cont}. An interesting question is to determine how the class $\COMP(\G)$ relates to $\BOREL(\G)$, $\MEAS(\G)$, and $\BAIRE(\G)$. The fourth named author proved \cite[Thm.~4]{FelixTrees}\footnote{Technically, \cite[Thm.~4]{FelixTrees} is a result for LCL problems on so-called highly computable regular trees, but the arguments easily transfer to the setting of computable $\mathbb{F}_n$-actions.} that for the free group $\mathbb{F}_n$,
    \[
       \COMP(\mathbb{F}_n) \,=\, \BAIRE(\mathbb{F}_n). 
    \]
    Furthermore, various problems have been shown to belong to both $\BAIRE(\G)$ and $\COMP(\G)$ for quite general classes of groups $\Gamma$, often with similar techniques; compare, e.g., \cite[Thm.~4.1]{Kierstead} and \cite[Corl.~1.5]{ASIalgorithms}, \cite{Schmerl2} and \cite{MeasurableBrooks}, or \cite{Schmerl1} and \cites[Thm.~B]{CM16}[Thm.~1.5]{FelixTwoEnded}. 
    In view of these facts, the following question appears natural: 
    
    \begin{question}\label{q:comp_vs_baire}
        Do the classes $\COMP(\G)$ and $\BAIRE(\G)$ coincide 
        for every finitely generated group $\Gamma$ with a decidable word problem?
    \end{question}
    
    We give a negative answer to this question for $\G = \Z^n$ with $n \geq 2$:
    
    \begin{tcolorbox}
	    \begin{theorem}\label{theo:comp}
	        For every integer $n \geq 2$, we have $\COMP(\Z^n) \setminus \BAIRE(\Z^n) \neq \0$.
	    \end{theorem}
	\end{tcolorbox}



    We do not know whether the opposite inclusion $\BAIRE(\Z^n) \subseteq \COMP(\Z^n)$ holds; indeed, it could even be that $\BAIRE(\G) \subseteq \COMP(\G)$ for \emph{every} finitely generated group $\G$. We remark that there do exist natural combinatorial problems that can be solved Baire-measurably but not computably. Namely, Manaster and Rosenstein \cite{Manaster2} constructed, for each $d \in \N$, a computable $d$-regular bipartite graph with no computable perfect matching. On the other hand, Bowen, Conley, and the fourth named author \cite{BaireMatchings} showed that when $d$ is odd, every Borel $d$-regular bipartite graph on a Polish space has a Baire-measurable perfect matching. Unfortunately, the graphs constructed by Manaster and Rosenstein are not generated by free group actions, so they do not directly shed any light on the relationship between the classes $\COMP(\G)$ and $\BAIRE(\G)$.

    \subsection{Summary, further results, and an overview of the remainder of the paper}

    All the results mentioned above, with the exception of Theorem \ref{theo:FFIID_Baire},
    follow from the existence of a single LCL problem $\Pi$:

    \begin{tcolorbox}
	    \begin{theorem}\label{theo:main}
	        For every $n \geq 2$, there exists an LCL problem $\Pi$ on $\Z^n$ such that:
            \begin{enumerate}[label=\ep{\normalfont{\roman*}}]
                \item\label{item:meas_yes} $\Pi \in \MEAS(\Z^n)$,

                \item\label{item:Baire_no} $\Pi \notin \BAIRE(\Z^n)$,

                \item\label{item:shift_yes} $\Pi \in \BS(\Z^n)$, 

                \item\label{item:ffiid_no} $\Pi \notin \FFIID(\Z^n)$,

                \item\label{item:comp_yes} $\Pi \in \COMP(\Z^n)$.
            \end{enumerate}
	    \end{theorem}
	\end{tcolorbox}


    It is clear that Theorem~\ref{theo:main} implies Theorems~\ref{theo:separation}, \ref{theo:not_universal}, \ref{theo:ffiid}, and \ref{theo:comp}. The main challenge in proving Theorem~\ref{theo:main} is to achieve the desired behavior with a \emph{locally} checkable problem. For example, a version of Theorem~\ref{theo:not_universal} \ep{non-universality of the shift in the Baire-measurable context} for \emph{non-local} problems was proved before in \cite[Corl.~2.12]{BerGeneric}; however, that result holds for arbitrary countably infinite groups, which shows that the subtle distinctions between, say, $\Z^n$ and $\mathbb{F}_n$ are lost when non-local problems are considered.

    The starting point for our construction 
    is a type of combinatorial structure called a \emph{toast}. This concept originates in the work of Conley and Miller \cite{CM16} (the term ``toast'' was coined by Miller) and has been widely used in descriptive combinatorics ever since \cite{bowen2021perfect,BPZ,trees,ContinuousAbelian,BorelAbelian,gao2022forcing,marks2017borel}. Related notions also appear in the theory of random processes \cite{FIIDSchreier,Holroyd,Spinka}. Recently, Greb\'ik and Rozho\v{n} initiated a systematic study of toast-based constructions, under the name of \textsf{TOAST} algorithms \cite{GRgrids}, and closely related \textsf{ASI} algorithms were studied by Qian and the fourth named author \cite{ASIalgorithms}. A central role in our arguments is played by \emph{toasts with rectangular pieces}, which we simply call \emph{rectangular toasts}. Briefly, a {rectangle} in a free action $\Z^n \acts X$ is a set of the form $R \cdot x$, where $x \in X$ and $R \subset \Z^n$ is an $n$-dimensional axis-parallel rectangle in $\Z^n$, and a rectangular toast $\mathcal{R}$ for the action $\Z^n \acts X$ is a nested family of rectangles (see Fig.~\ref{fig:waffle}). We formally define rectangular toasts and 
    establish their relevant properties in \S\ref{sec:waffles}. It turns out that one can always construct Borel rectangular toasts $\mathcal{R}$ such that the union of the rectangles in $\mathcal{R}$ is conull (Proposition~\ref{prop:MeasurableWaffle}), but it may be impossible to make it comeager (Proposition~\ref{prop:noBaireWaffle}). In other words, building a rectangular toast with a large union is a problem that can be solved measurably but not Baire-measurably, indicating that rectangular toasts may be employed toward proving Theorem~\ref{theo:separation}. Unfortunately, the definition of a rectangular toast is plainly non-local. Nevertheless, in \S\ref{sec:baking}, we find a way to incorporate rectangular toasts in an LCL problem and show that the resulting problem  
    witnesses Theorem~\ref{theo:main}. 
    %


        \begin{figure}[t]
			\centering
			\begin{tikzpicture}[yscale=1.8,xscale=1.8]
                \node (CONT) at (0,-1.6) {\large $\CONT(\Z^n)$};
                \node (BOREL) at (0,-0.5) {\large $\BOREL(\Z^n)$};
                \node (MEAS) at (2.1,1) {\large $\MEAS(\Z^n)$};
                \node (BAIRE) at (-2.1,1) {\large $\BAIRE(\Z^n)$};
                \node (FIID) at (2.1,5.6) {\large $\FIID(\Z^n)$};
                \node (BS) at (-2.1,5.6) {\large $\BS(\Z^n)$};
                \node (FFIID) at (3.8,3.3) {\large $\FFIID(\Z^n)$};
                \node (COMP) at (-3.8,3.3) {\large $\COMP(\Z^n)$};

                 \draw[-{Stealth},red,inner sep=1.5,very thick,dotted] (COMP) -- (FFIID) node [pos=0.33,sloped,fill=white] {\scriptsize\textbf{Thm.~\ref{theo:main}}};

                \draw[-{Stealth},blue,inner sep=1,very thick] (CONT) -- (BOREL);
                \draw[-{Stealth},blue,inner sep=1,very thick] (CONT) to[bend right=25] (FFIID);
                \draw[-{Stealth},blue,inner sep=1.5,very thick] (BOREL) -- (MEAS) node [midway,sloped,fill=white] {\scriptsize\textbf{Thm.~\ref{theo:separation}}};
                \draw[-{Stealth},red,inner sep=1.5,very thick,dotted] (MEAS) -- (FFIID) node [midway,sloped,fill=white] {\scriptsize\textbf{Thm.~\ref{theo:ffiid}}};
                \draw[-{Stealth},red,inner sep=1.5,very thick,dotted] (MEAS) -- (BAIRE) node [midway,sloped,fill=white] {\scriptsize\textbf{Thm.~\ref{theo:separation}}};
                \draw[-{Stealth},very thick,inner sep=0.5,circle] (MEAS) -- (FIID) node [pos=0.4,fill=white] {\large $=$?};
                \draw[-{Stealth},very thick,inner sep=0.5,circle] (BOREL) -- (BAIRE) node [midway,fill=white] {\large $=$?};
                \draw[-{Stealth},blue,very thick,inner sep=1.5] (FFIID) -- (FIID) node [midway,sloped,fill=white] {\scriptsize\textbf{Thm.~\ref{theo:ffiid}}};

                \draw[-{Stealth},red,very thick,dotted,inner sep=1.5] (FIID) -- (BAIRE) node [pos=0.35,sloped,fill=white] {\scriptsize\textbf{Thm.~\ref{theo:separation}}};

                \draw[-{Stealth},blue,inner sep=1.5,very thick] (BAIRE) -- (BS) node [midway,sloped,fill=white] {\scriptsize\textbf{Thm.~\ref{theo:not_universal}}};
                
                \draw[{Stealth}-{Stealth},red,very thick,dotted,inner sep=1.5] (BS) to node [midway,sloped,fill=white] {\scriptsize\textbf{\cite{Sequel}}} (COMP);
                
                \draw[{Stealth}-{Stealth},red,very thick,dotted,inner sep=1.5] (BS) to node [midway,sloped,fill=white] {\scriptsize\textbf{Thm.~\ref{theo:main_extra} and \cite{Sequel}}}(FIID);

                \draw[{Stealth}-{Stealth},red,very thick,dotted,inner sep=1.5] (BS) to node [pos=0.6,sloped,fill=white] {\scriptsize\textbf{Thm.~\ref{theo:main_extra} and \cite{Sequel}}} (MEAS);

                \draw[-{Stealth},blue,inner sep=1,very thick] (CONT) to[bend left=25] (COMP);

                \draw[-{Stealth},red,inner sep=1.5,very thick,dotted] (COMP) -- (BAIRE) node [midway,sloped,fill=white] {\scriptsize\textbf{Thm.~\ref{theo:comp}}};
                \draw[{Stealth}-{Stealth},red,very thick,dotted,inner sep=1.5] (MEAS) -- (COMP) node [midway,sloped,fill=white] {\scriptsize\textbf{Thm.~\ref{theo:main_extra} and \cite{Sequel}}};

                \draw[-{Stealth},blue,inner sep=1.5,very thick] (FFIID) -- (BS) node [pos=0.45,sloped,fill=white] {\scriptsize\textbf{Thms.~\ref{theo:FFIID_Baire} and \ref{theo:main}}};

                \draw[{Stealth}-{Stealth},red,inner sep=1.5,very thick,dotted] (FIID) to  node [midway,sloped,fill=white] {\scriptsize\textbf{Thm.~\ref{theo:main_extra} and \cite{Sequel}}} (COMP);

                

			\end{tikzpicture}
   \begin{minipage}[t]{15cm}
\captionof{figure}{\label{fig:relations}
   Complexity classes of LCL problems on $\Z^n$ with $n \geq 2$.  
   Here:}\small
   \begin{itemize}
       \item {\protect\tikz[baseline=-0.7ex]{\protect\node (A1) at (0,0) {\footnotesize$A$};
                \protect\node (B1) at (1.8,0) {\footnotesize$B$};
                \protect\draw[-{Stealth},blue,thick] (A1) -- (B1);}} indicates that $A$ is a strict subset of $B$,
        \item {\protect\tikz[baseline=-0.7ex]{\protect\node (A2) at (0,0) {\footnotesize$A$};
                \protect\node (B2) at (1.8,0) {\footnotesize$B$};
                \protect\draw[-{Stealth},black,thick,inner sep=0.5,circle] (A2) -- (B2) node [midway,fill=white] {$=$?};}} means that $A \subseteq B$ and we do not know whether $A = B$,
        \item {\protect\tikz[baseline=-0.7ex]{\protect\node (A3) at (0,0) {\footnotesize$A$};
                \protect\node (B3) at (1.8,0) {\footnotesize$B$};
                \protect\draw[-{Stealth},red,thick,dotted] (A3) -- (B3);}}
                indicates that $A$ is not a subset of $B$. 
   \end{itemize}
   All known inclusions/non-inclusions between these classes 
   follow from the ones shown. The relations established in this paper and in the forthcoming paper \cite{Sequel} are labeled accordingly.
   \end{minipage}
	\end{figure}
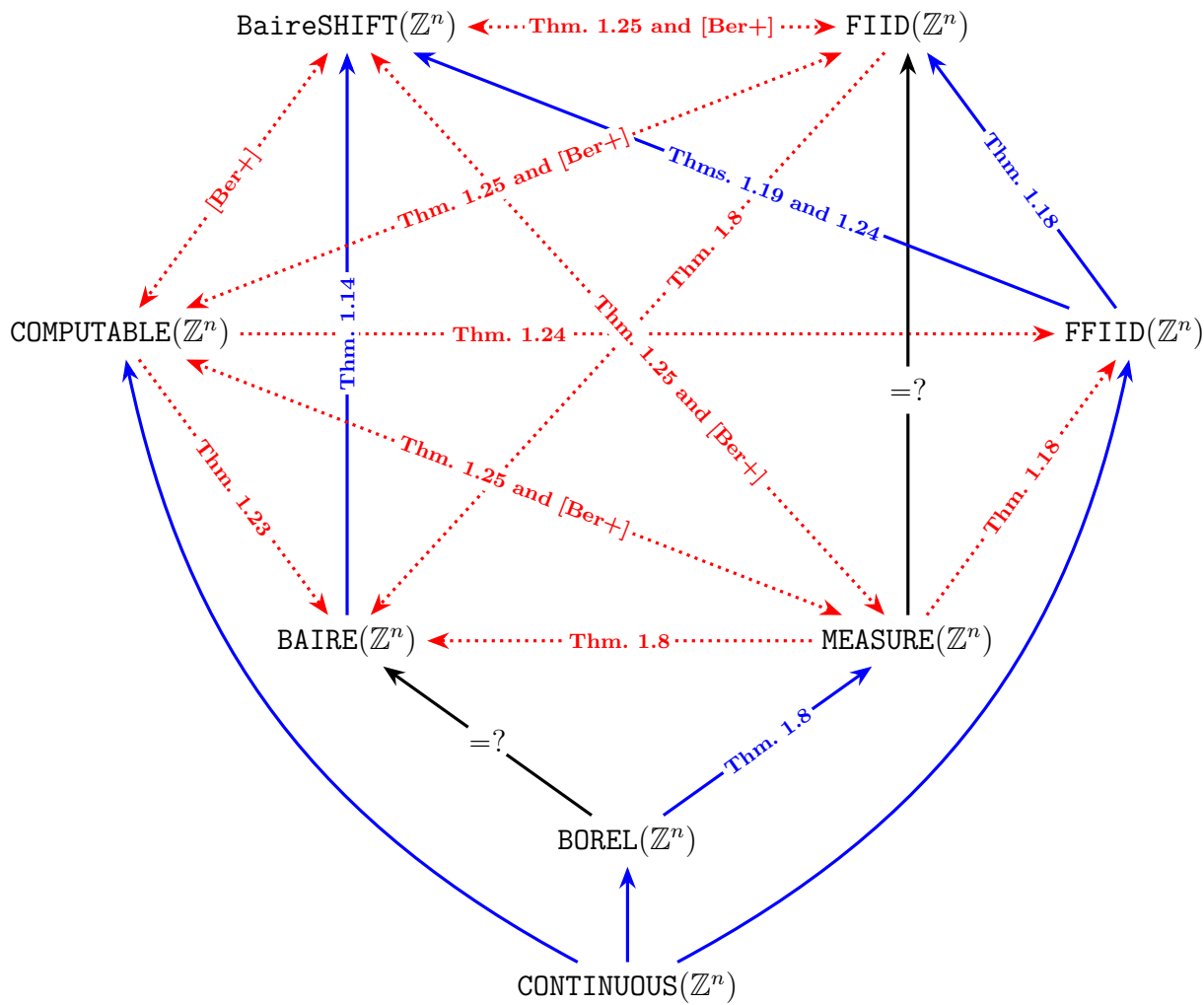

    The LCL problem we use to prove Theorem~\ref{theo:main} can be enhanced with extra layers of structure to derive further separation results between different complexity classes. In particular, in \S\ref{sec:extra} we modify the problem witnessing Theorem~\ref{theo:main} to obtain the following:

    \begin{tcolorbox}
	    \begin{theorem}\label{theo:main_extra}
	        For every $n \geq 2$, there exists an LCL problem $\Pi$ on $\Z^n$ such that:
            \begin{enumerate}[label=\ep{\normalfont{\roman*}}]
                \item $\Pi \notin \FIID(\Z^n)$,

                \item $\Pi \in \BS(\Z^n)$,  

                \item $\Pi \in \COMP(\Z^n)$.
            \end{enumerate}
	    \end{theorem}
	\end{tcolorbox}

    It follows that we have $\COMP(\Z^n) \not\subseteq \FIID(\Z^n)$, and hence
    $
        \COMP(\Z^n) \not \subseteq \MEAS(\Z^n)
    $. 
    Theorem~\ref{theo:main_extra} also makes a first step toward refuting the possibility that $\BAIRE(\Z^n)$ is a subset of $\MEAS(\Z^n)$, since it yields an LCL problem that is in $\BS(\Z^n)$ but not in $\MEAS(\Z^n)$. 

    We conclude the paper with the proof of Theorem~\ref{theo:FFIID_Baire} in \S\ref{sec:FFIID_Baire} and then some final remarks, including a list of open problems, in \S\ref{sec:final}.

    The method of adding extra structure to 
    the LCL problem from Theorem~\ref{theo:main} 
    that we use to establish Theorem~\ref{theo:main_extra} is very flexible and can be employed to obtain a variety of other results. Specifically, in the forthcoming second part of this study \cite{Sequel}, we will prove the following:

    \begin{theo}[\cite{Sequel}]\label{theo:announce}
        For every integer $n \geq 2$, the classes $\MEAS(\Z^n)$, $\COMP(\Z^n)$, and $\BS(\Z^n)$ are pairwise incomparable with respect to inclusion.
    \end{theo}

    Although the proof of Theorem~\ref{theo:announce} builds on the work in this paper, it requires a considerable amount of additional technical machinery to be developed; that is why we decided to present it separately. Theorem~\ref{theo:announce} clearly strengthens Theorem~\ref{theo:separation} (since $\BAIRE(\Z^n)$ is a strict subset of $\BS(\Z^n)$). In view of the inclusion $\FFIID(\Z^n) \subset \BS(\Z^n)$ (Theorem~\ref{theo:FFIID_Baire}), it also strengthens Theorem~\ref{theo:ffiid}. The incomparability between the three classes in Theorem~\ref{theo:announce} makes the contrast between the groups $\Z^n$ and $\mathbb{F}_n$ especially stark, since for the free group $\mathbb{F}_n$, the complexity classes defined in this paper---with the possible exception of $\FFIID(\mathbb{F}_n)$---are linearly ordered by inclusion. Theorem~\ref{theo:announce} also shows that, while $\CONT(\Z^n)$ is the \emph{minimum} element in the poset of the complexity classes studied in this paper, it does not have a \emph{maximum} element. (By contrast, for $\mathbb{F}_n$, the class $\BAIRE(\mathbb{F}_n) = \BS(\mathbb{F}_n) =  \COMP(\mathbb{F}_n)$ is maximum.)

    The current state of our knowledge regarding the complexity classes of LCL problems on $\Z^n$ for $n \geq 2$, including the results of this paper and those of \cite{Sequel}, is summarized in Fig.~\ref{fig:relations}.






    \section{Rectangular toasts}\label{sec:waffles}

    \subsection{Basic definitions}

    For the remainder of the paper, we fix an integer $n \geq 2$. Before stating the definition of a toast, we need some additional terminology. 
    We shall only define toasts for $\Z^n$-actions, although the definition can be easily generalized to other structures. For $\gamma$, $\delta \in \Z^n$, we write $\dist(\gamma, \delta) \defeq \|\gamma - \delta\|_1$. Similarly, if $\Z^n \acts X$ is a free action of $\Z^n$ on a set $X$ and $x$, $y \in X$ are in the same $\Z^n$-orbit, then there is a unique element $\gamma \in \Z^n$ such that $\gamma \cdot x = y$, and we let $\dist(x,y) \defeq \|\gamma\|_1$. If $x$ and $y$ are in distinct $\Z^n$-orbits, we set $\dist(x,y) \defeq \infty$. Note that $\dist$ is an \ep{extended} metric on $X$.\footnote{As a word of caution, when $X$ is a Polish space, the topology on $X$ is induced by a metric, but that metric is different from $\dist$ \ep{unless $X$ is countable and discrete}.} We say $x$, $y \in X$ are \emphd{adjacent} if $\dist(x,y) = 1$. As usual, for $x \in X$ and $A$, $B \subseteq X$, we let \[\dist(x, A) \,\defeq\, \inf \set{\dist(x,a) \,:\, a \in A} \quad \text{and} \quad \dist(A,B) \,\defeq\, \inf \set{\dist(a,b) \,:\, a \in A, \, b \in B}.\] 
    The \emphd{boundary} of $A \subseteq X$ is the set $\partial A$ of all $a \in A$ such that $\dist(a, X \setminus A) = 1$. 
    We write $\finset{X}$ to denote the set of all finite subsets of $X$.

    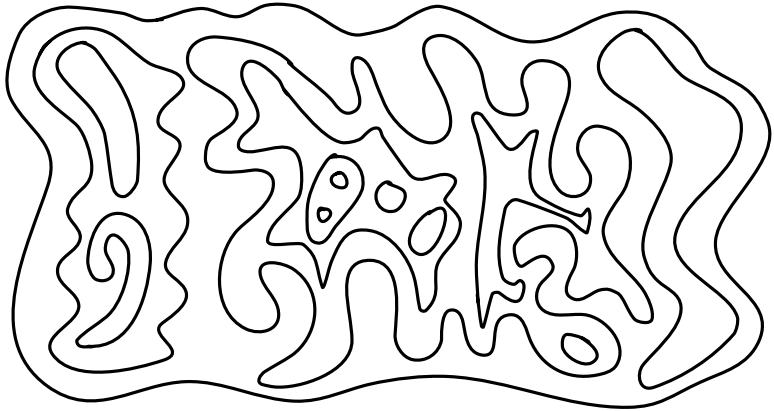
\begin{figure}[t]
			\centering
			\begin{tikzpicture}[x=0.75pt,y=0.75pt,yscale=-0.5,xscale=0.5]
\draw  [very thick] [line join = round][line cap = round] (187.53,32.54) .. controls (181.73,32.89) and (175.91,34.17) .. (170.12,33.58) .. controls (161.82,32.73) and (154.5,27.3) .. (146.25,26.03) .. controls (129.36,23.44) and (112.41,20.87) .. (95.34,20.06) .. controls (76.25,19.16) and (50.84,29.91) .. (42.62,48.43) .. controls (34.58,66.51) and (24.75,94.49) .. (35.26,113.88) .. controls (52.79,146.21) and (85.62,157.76) .. (71.96,200) .. controls (52.95,258.81) and (3.93,353.07) .. (71.13,402.33) .. controls (119.82,438.03) and (161.12,399.63) .. (211.65,398.04) .. controls (249.66,396.85) and (283.75,416.61) .. (321.7,417.77) .. controls (349.16,418.61) and (370.81,406.56) .. (396.73,400.66) .. controls (421.15,395.09) and (450.2,391.64) .. (475.7,392.86) .. controls (537.38,395.82) and (595.02,424.84) .. (657.02,423.99) .. controls (689.73,423.54) and (700.51,412.12) .. (730.71,398.39) .. controls (753.53,388.02) and (778.51,382.29) .. (789.46,356.41) .. controls (808,312.61) and (741.73,299.14) .. (745.46,255.89) .. controls (748.64,219.11) and (814.07,175.5) .. (797.5,134.2) .. controls (783.55,99.44) and (769.39,100.03) .. (741.61,82.51) .. controls (731.98,76.43) and (729.97,65.93) .. (722.97,57.35) .. controls (702.53,32.31) and (681.42,22.41) .. (649.42,25.39) .. controls (616.72,28.44) and (592.85,56.34) .. (558.91,55.09) .. controls (525.57,53.85) and (503.3,25.3) .. (468.1,19.42) .. controls (452.33,16.79) and (427.58,36.87) .. (420.65,39.65) .. controls (410.38,43.76) and (399.37,45.86) .. (388.49,47.83) .. controls (367.04,51.73) and (345.1,22.7) .. (319.88,18.49) .. controls (310.83,16.98) and (301.36,16.33) .. (292.35,18.09) .. controls (281.93,20.13) and (273.27,29.15) .. (262.66,29.68) .. controls (250.13,30.31) and (238.49,20.88) .. (225.96,21.4) .. controls (210.72,22.03) and (196.46,29.16) .. (181.71,33.04) ;
\draw  [very thick] [line join = round][line cap = round] (670.49,44.18) .. controls (667.89,40.78) and (656.93,43.23) .. (653.12,45.27) .. controls (639.42,52.61) and (622.2,69.72) .. (628.66,87.54) .. controls (635.55,106.55) and (653.97,111.98) .. (669.11,124.39) .. controls (682.77,135.59) and (704.45,157.38) .. (700.23,178.25) .. controls (694.84,204.95) and (664.91,244.94) .. (673.69,272.95) .. controls (676.8,282.88) and (679.19,288.81) .. (687.11,296.29) .. controls (695.99,304.68) and (709.47,308.68) .. (710.59,322.53) .. controls (712.37,344.71) and (680.87,369.42) .. (671.87,381.77) .. controls (665.39,390.65) and (668.3,408.81) .. (682.92,403.86) .. controls (713.23,393.62) and (742.72,374.06) .. (763.91,350.05) .. controls (766.18,347.48) and (768.84,337.64) .. (767.46,332.93) .. controls (759.19,304.72) and (694.98,285.4) .. (705.36,248.1) .. controls (716.92,206.56) and (772.95,188.28) .. (771.46,140.27) .. controls (770.54,110.82) and (733.12,104.11) .. (714.58,91.44) .. controls (697.42,79.71) and (684.38,47.8) .. (664.67,43.54) ;
\draw  [very thick] [line join = round][line cap = round] (64.07,67.07) .. controls (77.85,40.9) and (115.09,33.51) .. (139.29,49.36) .. controls (151.38,57.28) and (161.13,68.45) .. (173.33,76.2) .. controls (184.7,83.42) and (205.47,82.69) .. (208.7,95.78) .. controls (211.82,108.46) and (183.28,116.46) .. (183,134.94) .. controls (182.79,147.97) and (207.17,147.97) .. (205.19,161.28) .. controls (203.29,174.09) and (183.9,186.71) .. (190.59,200.35) .. controls (198.1,215.66) and (216.1,221.26) .. (212.64,240.2) .. controls (210.03,254.49) and (194.02,263) .. (190.84,277.49) .. controls (187.84,291.15) and (208.63,295.12) .. (211.56,308.27) .. controls (214.69,322.38) and (189.86,336.22) .. (185.71,340.53) .. controls (174.46,352.22) and (214.74,367.62) .. (188.87,376.29) .. controls (161.97,385.3) and (97.88,401.39) .. (76.35,370.08) .. controls (63.55,351.45) and (104.04,339.03) .. (103.19,321.34) .. controls (102.56,308.27) and (83.61,302.14) .. (81.83,288.44) .. controls (79.59,271.17) and (103.56,267.79) .. (104.67,252.49) .. controls (105.42,241.99) and (88.63,233.99) .. (92.14,222.94) .. controls (98.52,202.84) and (117.89,204.19) .. (116.16,178.94) .. controls (115.44,168.46) and (109.27,148.26) .. (101.81,140.02) .. controls (90.75,127.82) and (53.56,107.98) .. (58.4,82.61) .. controls (59.87,74.89) and (63.43,67.61) .. (67.28,60.76) .. controls (68.83,58) and (71.95,56.48) .. (74.28,54.35) ;
\draw  [very thick] [line join = round][line cap = round] (110.04,55.78) .. controls (96.07,54.29) and (73.11,70.8) .. (84.05,87.25) .. controls (92.52,99.99) and (104.76,109.78) .. (114.63,121.48) .. controls (118.96,126.61) and (127.81,140.67) .. (129.38,149.1) .. controls (132.45,165.58) and (133.97,182.32) .. (136.68,198.87) .. controls (137.69,205.03) and (141.34,211.97) .. (148.72,211.25) .. controls (154.18,210.72) and (161.93,199.75) .. (162.33,194.58) .. controls (165.26,156.31) and (162.81,128.52) .. (143.39,95.88) .. controls (139.95,90.11) and (122.35,65.08) .. (120.01,62.78) .. controls (112.65,55.56) and (102.33,55.63) .. (107.73,55.63) ;
\draw  [very thick] [line join = round][line cap = round] (112.56,362.38) .. controls (156.03,365.2) and (172,314.44) .. (175.15,279.07) .. controls (176.94,259.02) and (170.6,233.73) .. (147.78,228.76) .. controls (122.12,223.17) and (112.25,258.93) .. (112.21,276.8) .. controls (112.2,281.9) and (115.43,292.28) .. (120.3,294.76) .. controls (128.03,298.69) and (140.14,294.69) .. (137.67,283.17) .. controls (136.41,277.3) and (131.58,272.63) .. (130.22,266.79) .. controls (129.42,263.36) and (130.18,259.59) .. (131.4,256.28) .. controls (137.42,240.01) and (151.76,259.72) .. (152.91,266.1) .. controls (156.46,285.84) and (146.58,319.92) .. (130.51,333.08) .. controls (121.05,340.84) and (109.96,346.78) .. (101.66,355.77) .. controls (100.18,357.38) and (101.92,360.98) .. (103.88,361.94) .. controls (107.63,363.77) and (112.2,362.63) .. (116.36,362.97) ;
\draw  [color={rgb, 255:red, 0; green, 0; blue, 0 }  ,draw opacity=1 ][very thick] [line join = round][line cap = round] (318.15,77.68) .. controls (300.74,58.19) and (266.95,51.83) .. (242.09,50.15) .. controls (230.07,49.34) and (215.23,55.65) .. (213.13,68.85) .. controls (208.27,99.55) and (258.08,100.03) .. (262.71,125.87) .. controls (264.87,137.93) and (243.65,147.47) .. (237.21,152.6) .. controls (231.37,157.24) and (230.73,167.28) .. (230.89,173.96) .. controls (230.95,176.23) and (232.36,180.19) .. (234.05,181.9) .. controls (242.99,190.92) and (278.07,184.38) .. (280.46,184.37) .. controls (288.04,184.31) and (297.99,186.8) .. (299.26,195.71) .. controls (302.22,216.53) and (271.83,235.72) .. (262.95,245.38) .. controls (239.72,270.68) and (235.82,327.05) .. (271.19,344.57) .. controls (286.81,352.32) and (309.53,346) .. (303.94,323.96) .. controls (300.85,311.75) and (289.71,303.7) .. (286.58,291.7) .. controls (281.27,271.29) and (309.99,277.95) .. (319.68,284.59) .. controls (346.96,303.29) and (347.4,340.66) .. (326.63,363.76) .. controls (324,366.7) and (321.08,369.44) .. (317.8,371.65) .. controls (309.03,377.58) and (299.56,382.4) .. (290.63,388.08) .. controls (286.87,390.46) and (280.6,398.97) .. (288.75,401.1) .. controls (292.56,402.09) and (296.43,402.96) .. (300.34,403.37) .. controls (325.16,405.96) and (378.96,388.63) .. (378.72,356.51) .. controls (378.66,348.3) and (377.39,340.13) .. (376.26,332) .. controls (373.71,313.71) and (365.2,274.55) .. (397.61,275.62) .. controls (432,276.75) and (423.01,332.68) .. (422.18,350.94) .. controls (421.73,360.88) and (429.51,373.71) .. (439.74,375.85) .. controls (453.37,378.69) and (466.46,369.84) .. (469.04,356.02) .. controls (469.18,355.28) and (467.73,328.54) .. (477.03,326.08) .. controls (484.89,323.99) and (489.08,330.51) .. (490.64,337.67) .. controls (493.25,349.58) and (493.98,369.16) .. (510.18,371.26) .. controls (532.73,374.18) and (511.34,331.93) .. (538.69,330.52) .. controls (551.93,329.83) and (554.73,347.15) .. (558.27,355.67) .. controls (569.77,383.39) and (596.1,397.04) .. (625.35,391.48) .. controls (632.15,390.19) and (641.22,388.05) .. (645.23,381.76) .. controls (658.71,360.62) and (634.16,332.58) .. (612.68,330.42) .. controls (599.11,329.05) and (583.04,334.25) .. (571.98,326.27) .. controls (563.25,319.98) and (563.19,308.15) .. (568.58,300.03) .. controls (572.93,293.49) and (583.97,287.44) .. (581.16,280.11) .. controls (574.69,263.26) and (560.61,275.22) .. (551.61,280.99) .. controls (544.25,285.72) and (541.11,264.38) .. (542.83,260.33) .. controls (550.2,242.92) and (587.3,232.54) .. (599.65,248.83) .. controls (605.96,257.15) and (606.8,272.58) .. (603.16,282.77) .. controls (601.54,287.3) and (589.46,304.96) .. (597.78,311.92) .. controls (609.44,321.67) and (618.91,307.76) .. (629.05,305.51) .. controls (636.7,303.81) and (639.09,304.94) .. (646.42,308.62) .. controls (652.78,311.81) and (656.89,328.25) .. (664.02,335.3) .. controls (666.77,338.01) and (679.76,340.93) .. (681.34,335) .. controls (683.64,326.38) and (676.52,312.02) .. (672.95,306.49) .. controls (660.23,286.8) and (634.81,272.54) .. (632.16,247.8) .. controls (630.62,233.48) and (648.7,213.27) .. (653.22,204.69) .. controls (663.67,184.85) and (662.99,151.52) .. (637.98,142.74) .. controls (634.72,141.59) and (631.34,140.71) .. (627.92,140.27) .. controls (612.48,138.3) and (601.64,157.6) .. (608.19,170.8) .. controls (611.22,176.91) and (615.85,182.19) .. (618.55,188.46) .. controls (619.6,190.92) and (619.51,193.79) .. (619.24,196.45) .. controls (617.61,212.32) and (598.5,215.34) .. (588.85,206.02) .. controls (563.88,181.89) and (590.78,132.05) .. (596.4,110.08) .. controls (599.18,99.18) and (599.91,85.75) .. (588.9,79.25) .. controls (585.93,77.5) and (582.78,75.94) .. (579.43,75.11) .. controls (575.76,74.2) and (571.9,73.82) .. (568.13,74.12) .. controls (555.16,75.16) and (547.56,87.69) .. (550.77,100.07) .. controls (553.78,111.68) and (564.06,129.87) .. (540.22,128.83) .. controls (506.97,127.37) and (516.73,95.1) .. (505.05,75.6) .. controls (496.93,62.06) and (478.19,43.53) .. (460.21,49.91) .. controls (445.42,55.15) and (451.59,72.08) .. (457,81.92) .. controls (464.82,96.14) and (491.97,138.37) .. (469.63,153.64) .. controls (450.64,166.62) and (423.2,138.82) .. (414.78,123.94) .. controls (410.48,116.35) and (407.51,108.06) .. (404.13,100.02) .. controls (400.68,91.84) and (397.11,74.18) .. (386.07,71.26) .. controls (378.08,69.15) and (378.61,79.02) .. (378.97,84.48) .. controls (379.74,96.16) and (396.5,121.03) .. (377.74,126.66) .. controls (367.85,129.62) and (362.15,117.57) .. (356.92,112.25) .. controls (343.54,98.66) and (328.53,88.93) .. (314.75,76.15) ;
\draw  [color={rgb, 255:red, 0; green, 0; blue, 0 }  ,draw opacity=1 ][very thick] [line join = round][line cap = round] (617.66,357.84) .. controls (615.16,356.43) and (612.87,354.54) .. (610.16,353.6) .. controls (605.67,352.04) and (601,350.59) .. (596.25,350.39) .. controls (589.95,350.13) and (590.7,361.45) .. (591.12,363.37) .. controls (593.29,373.13) and (611.45,383.77) .. (620.32,379.1) .. controls (635.19,371.29) and (617.3,357.9) .. (609.77,353.75) ;
\draw  [color={rgb, 255:red, 0; green, 0; blue, 0 }  ,draw opacity=1 ][very thick] [line join = round][line cap = round] (505.93,318.28) .. controls (505.39,303.75) and (503.53,289.2) .. (504.31,274.68) .. controls (505.96,243.94) and (514.91,216.45) .. (512.64,185.7) .. controls (511.79,174.18) and (508.29,163.02) .. (505.93,151.71) .. controls (504.25,143.61) and (495.61,129.1) .. (503.62,127) .. controls (505.2,126.58) and (508.63,127.73) .. (509.78,129.22) .. controls (517.54,139.23) and (524.44,149.91) .. (532.47,159.7) .. controls (542.41,171.83) and (551.95,152.67) .. (558.07,147.08) .. controls (559.99,145.32) and (566.18,143.67) .. (565.82,146.24) .. controls (563.91,159.89) and (547.37,200.5) .. (569.86,213.37) .. controls (582.22,220.44) and (595.1,227.37) .. (608.98,230.58) .. controls (612.38,231.37) and (616.07,221.14) .. (617.41,224.37) .. controls (620.37,231.47) and (617.93,239.96) .. (615.98,247.4) .. controls (615.6,248.84) and (613.05,246.67) .. (611.79,245.87) .. controls (605.35,241.82) and (600.58,235.03) .. (593.54,232.16) .. controls (578,225.83) and (562.37,219.63) .. (546.33,214.7) .. controls (544.18,214.04) and (541.45,214.26) .. (539.62,215.59) .. controls (536.66,217.74) and (533.73,220.78) .. (533.01,224.37) .. controls (528.06,249.2) and (523.55,269.45) .. (530.5,292.34) .. controls (531.82,296.71) and (538.25,299.81) .. (542.04,299.79) .. controls (544.6,299.77) and (546.37,295.34) .. (548.8,296.14) .. controls (555.25,298.25) and (550.93,310.29) .. (549.93,312.56) .. controls (544.89,324.07) and (530.67,310.62) .. (527.14,308.32) .. controls (525.95,307.54) and (520.91,307.73) .. (520.24,309.65) .. controls (516.43,320.5) and (514.88,332.05) .. (511.21,342.95) .. controls (510.99,343.61) and (509.47,342.7) .. (509.34,342.01) .. controls (507.36,331.54) and (506.48,320.9) .. (505.05,310.34) ;
\draw  [color={rgb, 255:red, 0; green, 0; blue, 0 }  ,draw opacity=1 ][very thick] [line join = round][line cap = round] (314.4,108.62) .. controls (309.37,105.64) and (308.73,98.32) .. (305.13,93.72) .. controls (301.51,89.1) and (275.93,62.88) .. (268.73,80.7) .. controls (261.37,98.87) and (289.15,118.3) .. (283.18,133.92) .. controls (279.41,143.76) and (262.87,148.7) .. (264.24,160.7) .. controls (264.9,166.49) and (276.07,166.35) .. (279.73,165.83) .. controls (304.1,162.41) and (324.66,151.35) .. (327.67,184.82) .. controls (328.25,191.27) and (328.76,206.69) .. (323.92,213.33) .. controls (313.73,227.31) and (294.61,238.14) .. (293.54,255.41) .. controls (292.96,264.71) and (324.63,255.44) .. (329.99,259.95) .. controls (345.73,273.18) and (341.84,282.2) .. (348.93,302.86) .. controls (349.32,303.98) and (350.66,301.12) .. (351.05,300) .. controls (353.93,291.82) and (355.33,283.13) .. (358.7,275.14) .. controls (367.69,253.78) and (374.75,241.55) .. (397.66,245.74) .. controls (418.93,249.63) and (435.64,272.84) .. (441.27,292.5) .. controls (444.24,302.86) and (444.74,313.79) .. (447.09,324.31) .. controls (447.66,326.9) and (451.66,327.83) .. (453.45,326.09) .. controls (471.58,308.43) and (457.69,293.2) .. (466.13,277.01) .. controls (473.95,262) and (486.47,257.51) .. (488.08,239.82) .. controls (489.16,227.84) and (470.15,218.17) .. (472.54,213.58) .. controls (476.1,206.72) and (487.41,197.28) .. (481.57,192.22) .. controls (468.31,180.76) and (448.48,203.61) .. (436.48,188.47) .. controls (426.95,176.44) and (418.19,163.79) .. (409.7,150.99) .. controls (408.54,149.23) and (408.26,143.1) .. (404.67,142.21) .. controls (397.55,140.44) and (391.24,150.42) .. (387.45,152.86) .. controls (383.04,155.7) and (373.31,157.52) .. (368.71,156.76) .. controls (348.5,153.42) and (319.94,114.4) .. (310.11,102.94) ;
\draw  [color={rgb, 255:red, 0; green, 0; blue, 0 }  ,draw opacity=1 ][very thick] [line join = round][line cap = round] (363.04,171.06) .. controls (362.58,164.88) and (354.22,173.68) .. (353.32,174.96) .. controls (349.67,180.1) and (346.09,185.36) .. (343.45,191.09) .. controls (333.47,212.74) and (332.71,219.76) .. (333.64,244.11) .. controls (333.77,247.44) and (335.71,254.57) .. (339.56,256.64) .. controls (351.61,263.12) and (356.74,248.49) .. (362.1,241.4) .. controls (365.42,237) and (387.55,214.13) .. (389.52,201.74) .. controls (392.49,183.1) and (378.58,167.01) .. (358.3,171.6) ;
\draw  [color={rgb, 255:red, 0; green, 0; blue, 0 }  ,draw opacity=1 ][very thick] [line join = round][line cap = round] (418.04,198.24) .. controls (414.58,193.64) and (408.15,198.77) .. (405.65,202.33) .. controls (397.77,213.57) and (407.69,229.61) .. (421.78,225.81) .. controls (428.47,224.01) and (437.37,209.59) .. (429.33,204.21) .. controls (425.11,201.38) and (419.99,200.2) .. (415.32,198.19) ;
\draw  [color={rgb, 255:red, 0; green, 0; blue, 0 }  ,draw opacity=1 ][very thick] [line join = round][line cap = round] (457.89,226.27) .. controls (447.23,230.94) and (425.43,257.18) .. (443.59,268.25) .. controls (465.32,281.5) and (490.53,201.21) .. (456.31,228.59) ;
\draw  [color={rgb, 255:red, 0; green, 0; blue, 0 }  ,draw opacity=1 ][very thick] [line join = round][line cap = round] (366.29,187.7) .. controls (353.77,188.24) and (355.58,201.87) .. (368.31,202.4) .. controls (379.86,202.88) and (370.48,180.45) .. (362.25,188.69) ;
\draw  [color={rgb, 255:red, 0; green, 0; blue, 0 }  ,draw opacity=1 ][very thick] [line join = round][line cap = round] (350.56,224.5) .. controls (341.89,216.04) and (342.65,243.53) .. (353.37,234.61) .. controls (362.57,226.95) and (353.45,220.48) .. (345.87,225.38) ;
\end{tikzpicture}
\caption{A toast for $\Z^2$.}\label{fig:toast}
    \end{figure}

    \begin{definition}[Toasts]\label{defn:toast}
        Let $\Z^n \acts X$ be a free action of $\Z^n$ on a set $X$. A family $\mathcal{T} \subseteq \finset{X}$ of finite sets is a \emphd{$q$-toast}, where $q \in \N$, if the following two conditions hold:
    \begin{enumerate}[label=\ep{\normalfont{}\small{T}\arabic*}]
        \item\label{item:nested} for all $K$, $L \in \mathcal{T}$, either $K \cap L = \0$, or $K \subseteq L$, or $L \subseteq K$, and
        \item\label{item:far_boundaries} for distinct $K$, $L \in \mathcal{T}$, we have $\dist(\partial K, \partial L) > q$.
    \end{enumerate}
    A set $K \in \mathcal{T}$ is called a \emphd{piece} of the $q$-toast $\mathcal{T}$. The \emphd{boundary} of a $q$-toast $\mathcal{T}$ is the set $\partial \mathcal{T} \defeq \bigcup_{K \in \mathcal{T}} \partial K$. A $q$-toast $\mathcal{T}$ is \emphd{complete}\footnote{Usually completeness is treated as part of the requirements for being a $q$-toast. However, it will be more convenient for us to define it separately.} if $\bigcup\mathcal{T} = X$, i.e., if for all $x \in X$, there is $K \in \mathcal{T}$ with $x \in K$. See Fig.~\ref{fig:toast} for an illustration in the $n = 2$ case.
    \end{definition}

    If $X$ is a standard Borel space, then the set $\finset{X}$ of all finite subsets of $X$ also carries a natural standard Borel structure; namely, a subset of $\finset{X}$ is Borel if and only if its preimage in $\bigsqcup_{k = 0}^\infty X^k$ under the map $(x_0, \ldots, x_{k-1}) \mapsto \set{x_0, \ldots, x_{k-1}}$ is Borel \ep{see, e.g., \cite[\S2.1]{BerFelixASI}}. We say that a family $\mathcal{T}$ of finite subsets of $X$ (for example, a $q$-toast) is \emphd{Borel} if it is a Borel subset of $\finset{X}$. Note that if $\Z^n \acts X$ is a free Borel action of $\Z^n$ and $\mathcal{T} \subseteq \finset{X}$ is a Borel $q$-toast, then every element $x \in X$ may belong to only countably many pieces $K \in \mathcal{T}$. By the Luzin--Novikov uniformization theorem \cite[Thm.~18.10]{KechrisDST}, it follows that the sets $\bigcup\mathcal{T}$ and $\partial \mathcal{T}$ are Borel.
    
    The following fundamental fact was established by Gao, Jackson, Krohne, and Seward:

    \begin{theorem}[{Gao--Jackson--Krohne--Seward \cite[Corl.~4.13]{BorelAbelian}}]\label{theo:BorelToast}
        If $\Z^n \acts X$ is a free Borel action on a standard Borel space, then for each $q \in \N$, there exists a Borel complete $q$-toast $\mathcal{T} \subseteq \finset{X}$.
    \end{theorem}

    We approach Theorem~\ref{theo:main} by studying $q$-toasts with additional constraints on the shapes of their pieces. Specifically, we consider $q$-toasts all of whose pieces are $n$-dimensional rectangles. We use $[a,b]$ to denote the closed interval from $a$ to $b$ in $\Z$, i.e., $[a,b] \defeq \set{c \in \Z \,:\, a \leq c \leq b}$.

\begin{figure}[t]
			\centering
			\begin{tikzpicture}[scale=0.8]
                \draw[very thick] (0,0.5) rectangle (1,2.5);
                \draw[very thick] (2,1) rectangle (3,2);
                \draw[very thick] (2.5,-0.5) rectangle (3.5,0.5);
                \draw[very thick] (-0.5,-0.5) rectangle (2,0);
                \draw[very thick] (-1,-1) rectangle (4,3);
                \draw[very thick] (4.5,-0.5) rectangle (6,2.5);
                \draw[very thick] (-1.5,-2) rectangle (6.5,4);
                \draw[very thick] (-3.5,-2) rectangle (-2,4);
                \draw[very thick] (-3,-1.5) rectangle (-2.5,0.5);
                \draw[very thick] (-3,1.5) rectangle (-2.5,3.5);
                \draw[very thick] (7,-2) rectangle (8,4);
                \draw[very thick] (-4,-2.5) rectangle (8.5,4.5);
			\end{tikzpicture}
   \caption{A rectangular toast for $\Z^2$.}\label{fig:waffle}
	\end{figure}
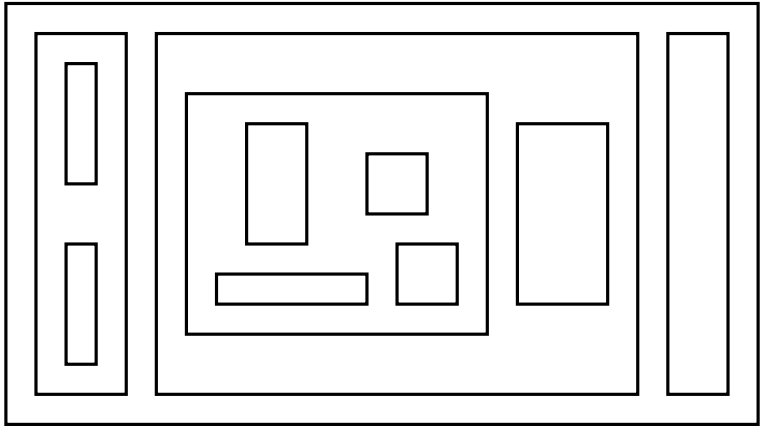

    \begin{definition}[Rectangles and rectangular toasts]
        A \emphd{rectangle} in $\Z^n$ is a set of the form \[[a_1, b_1] \times \cdots \times [a_n, b_n],\] where $a_1 \leq b_1$, \ldots, $a_n \leq b_n$ are integers. The values $b_i - a_i$ for $1 \leq i \leq n$ are called the \emphd{side lengths} of the rectangle. A rectangle all of whose side lengths are equal \ep{to $N$} is called an \ep{$N$-}\emphd{square}.
        
        Given a free action $\Z^n \acts X$, a subset of $X$ is a \emphd{rectangle} \ep{resp.~a \emphd{square}}
        if it has the form $R \cdot x \defeq \set{\gamma \cdot x \,:\, \gamma \in R}$ for a point $x \in X$ and a rectangle \ep{resp.~a square}
        $R \subset \Z^n$. The \emphd{side lengths} of $R \cdot x$ are those of $R$. 
        A \emphd{rectangular} \ep{resp.~\emphd{square}} \emphd{$q$-toast} is a $q$-toast all of whose pieces are rectangles \ep{resp.~squares} with side lengths at least $q$. See Fig.~\ref{fig:waffle} for an illustration in the $n = 2$ case.
    \end{definition}

    \begin{remk}\label{remk:only_bdry_matters}
        Since we are assuming that $n \geq 2$, if $K$, $L$ are two rectangles with $\partial K \cap \partial L \neq \0$, then either $K \cap L = \0$, or $K \subset L$, or $L \subset K$. Hence, when $\mathcal{T}$ is a set of rectangles, condition \ref{item:far_boundaries} in Definition~\ref{defn:toast} implies \ref{item:nested}.
    \end{remk}

    Let $\Z^n \acts X$ be a free action of $\Z^n$ and let $\mathcal{R} \subseteq \finset{X}$ be a rectangular $q$-toast. An important role in our arguments will be played by the points of $X$ whose entire $\Z^n$-orbit is covered by $\bigcup \mathcal{R}$. That is, we say that a point $x \in X$ is \emphd{$\mathcal{R}$-complete} if $\Z^n \cdot x \subseteq \bigcup \mathcal{R}$. Equivalently, $x\in X$ is $\mathcal{R}$-complete if $x \notin \Z^n \cdot (X \setminus \bigcup \mathcal{R})$. The following simple observation will be helpful:
    
    \begin{lemma}\label{lemma:Wcomplete}
        Let $\Z^n \acts X$ be a free action of $\Z^n$ and let $\mathcal{R} \subseteq \finset{X}$ be a rectangular $q$-toast for some $q \geq 1$. If $x \in X$ is an $\mathcal{R}$-complete point, then each $y \in \Z^n \cdot x$ is contained in infinitely many pieces of $\mathcal{R}$, and the set $\Z^n \cdot x$ includes pieces of $\mathcal{R}$ all of whose side lengths are arbitrarily large.
    \end{lemma}
    \begin{scproof}
        It suffices to show that for each rectangle $K \in \mathcal{R}$ with $K \subseteq \Z^n \cdot x$, there exists another rectangle $K' \in \mathcal{R}$ such that $K' \supsetneq K$ (which implies that every side length of $K'$ strictly exceeds the corresponding side length of $K$). To this end, take a point $y \notin K$ adjacent to the boundary of $K$ and note that any rectangle $K' \in \mathcal{R}$ containing $y$ must satisfy $K' \supsetneq K$. 
    \end{scproof}

    
    In \S\S\ref{subsec:waffle_measure}--\ref{subsec:waffle_comput} 
    we present a series of facts about rectangular toasts used to prove Theorem~\ref{theo:main}.

    \subsection{Rectangular toasts and measure theory}\label{subsec:waffle_measure}

    We start by showing that every free Borel action of $\Z^n$ admits a square $q$-toast that is ``almost'' complete from the point of view of measure theory. 

    \begin{prop}\label{prop:MeasurableWaffle}
        Let $\Z^n \acts (X,\mu)$ be a free Borel action of $\Z^n$ on a standard probability space $(X,\mu)$. Then for each $q \in \N$, there exists a Borel square $q$-toast $\mathcal{R} \subseteq \finset{X}$ with $\mu(\bigcup \mathcal{R}) = 1$.  
        %
    \end{prop}

    This statement seems well known to experts in the area, although we were unable to find a proof in the literature. For completeness, we prove it here. First we need to review some concepts from the theory of Borel graphs. 
    A graph $G$ is \emphd{Borel} if its vertex set $V(G)$ is a standard Borel space and its edge set $E(G)$ is a Borel subset of $\fs{V(G)}{2}$. \ep{Here we write $\fs{V(G)}{2} \subset \finset{V(G)}$ for the set of all $2$-element subsets of $V(G)$.} Borel graphs are the main objects of study in descriptive combinatorics; see \cite{KechrisMarks,Pikh_survey} for an overview.
    
    A \emphd{subgraph} of a graph $G$ is a graph $H$ such that $V(H) \subseteq V(G)$ and $E(H) \subseteq E(G)$. We write $H \subseteq G$ to mean that $H$ is a subgraph of $G$. If $G$ and $H$ are Borel graphs and $V(H) \subseteq V(G)$ is a Borel set equipped with the relative Borel $\sigma$-algebra, then we call $H$ a \emphd{Borel subgraph} of $G$. 
    
    \begin{definition}[Hyperfinite graphs]
        A graph $G$ is \emphd{component-finite} if every connected component of $G$ is finite. A Borel graph $G$ is \emphd{hyperfinite} if there is an increasing sequence $G_0 \subseteq G_1 \subseteq \cdots \subseteq G$ of component-finite Borel subgraphs of $G$ whose union is $G$.
    \end{definition}

    The systematic study of hyperfinite graphs was initiated by Weiss \cite{WeissHyperfinite} and Slaman and Steel \cite{SS}, with important foundational work done by Dougherty, Jackson, and Kechris \cite{DJK} and Jackson, Kechris, and Louveau \cite{jackson2002countable}, among others. For an overview of this topic, see \cite{KechrisCBER}. We shall use the following fact:

    \begin{theorem}[{Weiss (unpublished); see \cite{jackson2002countable}}]\label{theo:hyperfinite}
        Let $\Z^n \acts X$ be a free Borel action of $\Z^n$ on a standard Borel space $X$ and let $\bm{0} \notin S \subseteq \Z^n$ be a finite subset. If $G$ is the graph with vertex set $X$ and edge set $E(G) \defeq \set{\set{x, \sigma \cdot x} \,:\, x \in X, \, \sigma \in S}$, then $G$ is hyperfinite.
    \end{theorem}

    We remark that Theorem~\ref{theo:hyperfinite} is a corollary of Theorem~\ref{theo:BorelToast}, since the existence of a Borel $1$-toast implies hyperfiniteness \cite{BorelAbelian}. A long-standing open question raised by Weiss is whether a version of Theorem~\ref{theo:hyperfinite} holds for every amenable countable group in place of $\Z^n$ \cite{WeissHyperfinite}; for partial results see \cite{jackson2002countable} by Jackson, Kechris, and Louveau and \cite{conley2020borel} by Conley, Jackson, Marks, Seward, and Tucker-Drob.
    However, if null sets are ignored (as indeed they are in this subsection), Weiss's question has a positive answer due to Ornstein and Weiss \cite{ornstein1980ergodic}. 
    

    A set $I \subseteq V(G)$ is \emphd{independent} in a graph $G$ if no two vertices in $I$ are adjacent. The \emphd{chromatic number} $\chi(G)$ of $G$ is the minimum $k \in \N$ \ep{if it exists} such that $V(G)$ can be covered by $k$ independent sets; if there is no such $k$, we set $\chi(G) \defeq \infty$. Equivalently, $\chi(G)$ is the least $k$ such that $G$ has a proper $k$-coloring, i.e., a labeling $f \colon V(G) \to k$ such that adjacent vertices receive different labels \ep{we have already discussed this concept in Example~\ref{exmp:coloring}}. The chromatic number of a graph is of fundamental importance in graph theory; see \cite{Coloring1,Coloring2} for an overview. Naturally, it is also central to the study of Borel graphs \cite{KST}. For a Borel graph $G$ and a probability Borel measure $\mu$ on $V(G)$, the \emphd{$\mu$-approximate chromatic number} $\chi_\mu^{\mathsf{ap}}(G)$ of $G$ is the minimum $k \in \N$ \ep{if it exists} such that for any $\epsilon > 0$, we can find Borel independent sets $I_1$, \ldots, $I_k \subseteq V(G)$ with $\mu(I_1 \cup \ldots \cup I_k) \geq 1 - \epsilon$. A graph $G$ is \emphd{locally finite} if every vertex of $G$ has finitely many neighbors. We shall use the following result of Conley and Kechris:

    \begin{theorem}[{Conley--Kechris \cite[Thm.~3.8]{ConleyKechris}}]\label{theo:appchi}
        If $G$ is a locally finite hyperfinite Borel graph and $\mu$ is a probability Borel measure on $V(G)$, then $\chi_\mu^{\mathsf{ap}}(G) \leq \chi(G)$.
    \end{theorem}

    The following lemma is a consequence of Theorems~\ref{theo:hyperfinite} and \ref{theo:appchi}. 
    For measure-preserving actions, it essentially follows from the Ornstein--Weiss construction of quasi-tilings, cf.~\cite[Thm.~5]{ornstein.weiss}.

    \begin{lemma}\label{lemma:almost_cover_by_squares}
            Let $\Z^n \acts (X,\mu)$ be a free Borel action on a standard probability space. For all $q \in \N$ and $\epsilon > 0$, there exist $N = N(n, q,\epsilon) \in \N$ \ep{depending only on $n$, $q$, and $\epsilon$ but not on the action $\Z^n \acts X$ or the measure $\mu$} and a Borel set $\mathcal{D} \subseteq \finset{X}$ of $N$-squares such that $N \geq q$ and:
            \begin{itemize}
                \item $\mu( \bigcup \mathcal{D} ) \geq 1 - \epsilon$, and
                \item for distinct $C$, $D \in \mathcal{D}$, we have $\dist(C,D) > q$. 
         \end{itemize}
    \end{lemma}
    \begin{scproof}
        Fix $r \in \N$ and $\epsilon > 0$ and let $N \geq q$ be a large natural number to be determined later. Define a new probability measure $\nu$ on $X$ by the formula
        \[
            \nu(A) \,\defeq\, \frac{1}{(N+1)^n} \sum_{\delta \in [0, N]^n} \mu(\delta \cdot A). 
        \]
        \ep{Note that if the action $\Z^n \acts (X,\mu)$ is measure-preserving, then $\nu = \mu$.} Let $G$ be the graph with vertex set $X$ and edge set
        \[
            E(G) \,\defeq\, \big\{\set{x,y} \in \fs{X}{2} \,:\,x \neq y \text{ and } \dist([0,N]^n \cdot x , \, [0,N]^n \cdot y) \,\leq\, q\big\}.
        \]
        By Theorem~\ref{theo:hyperfinite} applied with
        \[
            S \,\defeq\,\set{\gamma \in \Z^n \,:\, \gamma \neq 0 \text{ and } \dist([0,N]^n, \, [0,N]^n + \gamma) \,\leq\,q},
        \]
        the graph $G$ is hyperfinite. Therefore, by Theorem~\ref{theo:appchi}, $\chi^\mathsf{ap}_\nu(G) \leq \chi(G)$. Note that each connected component of $G$ is isomorphic to the Cayley graph $\Cay(\Z^n, S)$, and hence $\chi(G) = \chi(\Cay(\Z^n, S))$. The graph $\Cay(\Z^n, S)$ can be partitioned into $k \defeq (N+q+1)^n$ independent sets, namely
        \[
            I_\eta \,\defeq\, (N+q+1)\Z^n + \eta, \quad \text{for each } \eta \in [0,N+q]^n.
        \]
        Therefore, $\chi(\Cay(\Z^n, S)) \leq k$.\footnote{In fact, $\chi(\Cay(\Z^n, S)) = k$, but we only need the upper bound.} It follows that $\chi^\mathsf{ap}_\nu(G)  \leq k$. Therefore, there exist Borel independent sets $I_1$, \ldots, $I_k \subseteq X$ in $G$ such that $\nu (I_1 \cup \ldots \cup I_k) \geq 1 - \epsilon/2$. Let $I$ be the set among $I_1$, \ldots, $I_k$ whose $\nu$-measure is maximum. Then $\nu(I) \geq (1-\epsilon/2)/k$. We claim that the following set $\mathcal{D}$ of $N$-squares satisfies the requirements of the lemma:
        \[
            \mathcal{D} \,\defeq\, \set{[0,N]^n \cdot x \,:\, x \in I}.
        \]
        Since $I$ is independent in $G$, distinct $C$, $D \in \mathcal{D}$ satisfy $\dist(C,D) > q$, as desired. It remains to compute a lower bound on $\mu(\bigcup \mathcal{D})$. To this end, we observe that
        \[
            \bigcup \mathcal{D} \,=\, [0,N]^n \cdot I \,=\, \bigcup_{\delta \in [0,N]^n} (\delta \cdot I).
        \]
        Thus, as the sets $\delta \cdot I$ for $\delta \in [0,N]^n$ are pairwise disjoint, we have
        \[
            \mu\left(\bigcup \mathcal{D}\right) \,=\, \sum_{\delta \in [0,N]^n} \mu(\delta \cdot I) \,=\, (N+1)^n \,\nu(I) \,\geq\, \left(\frac{N+1}{N+ q + 1}\right)^n \, (1 - \epsilon/2) \,\geq\, 1 - \epsilon,
        \]
        where the last inequality holds for all $N$ large enough as a function of $n$, $q$, and $\epsilon$, as desired.
    \end{scproof}

    
    \begin{scproof}[ of Proposition~\ref{prop:MeasurableWaffle}]
        Fix a free Borel action $\Z^n \acts (X,\mu)$ on a standard probability space and $q \in \N$. Without loss of generality, we may assume that $q \geq 1$. 
        For each $i \in \N$, we recursively define parameters $\epsilon _i > 0$ and $L_i$, $N_i \in \N$ as follows:
        \[
            L_i \,\defeq\, 1 + q + \sum_{j < i} N_j, \quad \epsilon_i \,\defeq\, 2^{-i} (1 + 2L_i)^{-n}, \quad N_i \,\defeq\, N(n, q, \epsilon_i).
        \]
        Here $N(\cdot, \cdot, \cdot)$ is the function from Lemma~\ref{lemma:almost_cover_by_squares}. We also define a probability measure $\nu_i$ on $X$ by
        \[
            \nu_i(A) \,\defeq\, \frac{1}{(1 + 2L_i)^n} \sum_{\delta \in [-L_i,L_i]^n} \mu(\delta \cdot A).
        \]
        %
        %
        %
        \ep{We again note that $\nu_i = \mu$ if the action $\Z^n \acts (X,\mu)$ is measure-preserving.} Lemma~\ref{lemma:almost_cover_by_squares} yields a Borel set $\mathcal{D}_i \subseteq \finset{X}$ of $N_i$-squares such that:
        \begin{enumerate}[label=\ep{\normalfont\small{}{D}\arabic*}]
                \item\label{item:almost_cover} $\nu_i( \bigcup \mathcal{D}_i ) \geq 1 - \epsilon_i$, and
                \item\label{item:far} for distinct $C$, $D \in \mathcal{D}_i$, we have $\dist(C,D) > q$. 
         \end{enumerate}
        Let $\mathcal{D} \defeq \bigcup_{i \in \N} \mathcal{D}_i$ and define
        \[
            \mathcal{Q}_i \,\defeq\, \big\{C \in \bigcup_{j < i} \mathcal{D}_j \,:\, \text{$\dist(\partial C, \partial D) \leq q$ for some $D \in \mathcal{D}_i$}\big\} \quad \text{and} \quad \mathcal{Q}\,\defeq\, \bigcup_{i \in \N} \mathcal{Q}_i.
        \]
        Note that for each $C \in \mathcal{Q}_i$, there are only countably many sets $D \in \mathcal{D}_i$ as in the definition of $\mathcal{Q}_i$, and thus $\mathcal{Q}_i$ is a Borel set by the Luzin--Novikov theorem \cite[Thm.~18.10]{KechrisDST}, which implies that $\mathcal{Q}$ is also Borel. We claim that the set $\mathcal{R} \defeq \mathcal{D} \setminus \mathcal{Q}$ is a desired square $q$-toast.

        Clearly, $\mathcal{R}$ is a Borel set of squares of side length at least $q$. Consider any two distinct sets $C$, $D \in \mathcal{R}$ and let $i$, $j$ be such that $C \in \mathcal{D}_j$, $D \in \mathcal{D}_i$. Without loss of generality, assume that $j \leq i$. If $j = i$, then $\dist(C,D) > q$ by \ref{item:far}, while if $j < i$, then $\dist(\partial C, \partial D) > q$ because $C \notin \mathcal{Q}_i$. By Remark~\ref{remk:only_bdry_matters}, it follows that $\mathcal{R}$ is indeed a square $q$-toast. It remains to verify that $\mu(\bigcup \mathcal{R}) = 1$.
        
        First we note that by \ref{item:almost_cover},
        \begin{equation}\label{eq:bound1}
            \mu\left(X \setminus \bigcup \mathcal{D}_i\right) \,\leq\, (1 + 2L_i)^n \,\nu_i\left(X \setminus \bigcup \mathcal{D}_i\right) \, \leq\, 2^{-i}.
        \end{equation}
        Next we bound $\mu(\bigcup \mathcal{Q}_i)$. Suppose $j < i$ and $C \in \mathcal{D}_j \cap \mathcal{Q}_i$. For each point $x \in C$, the set $[-N_j,N_j]^n \cdot x$ includes $C$, so, letting $D \in \mathcal{D}_i$ be such that $\dist(\partial C, \partial D) \leq q$, we have
        \[
            \left([-q - N_j,\, q + N_j]^n \cdot x]\right) \,\cap\, \partial D \,\neq\, \0.
        \]
        Since the sets in $\mathcal{D}_i$ are a distance at least $2$ apart, it follows that the set $[-1- q-N_j, 1 + q+N_j]^n \cdot x$ contains a point not in $\bigcup \mathcal{D}_i$. As $1 + q + N_j \leq L_i$, we conclude that
        \[
            \bigcup \mathcal{Q}_i \,\subseteq\, [-L_i, L_i]^n \cdot \left(X \setminus \bigcup \mathcal{D}_i\right).
        \]
        Therefore, using \ref{item:almost_cover} again, we obtain
        \begin{equation}\label{eq:bound2}
            \mu\left(\bigcup \mathcal{Q}_i\right) \,\leq\, \sum_{\delta \in [-L_i, L_i]^n} \mu\left(\delta \cdot \left(X \setminus \bigcup \mathcal{D}_i\right)\right) \,=\, (1 + 2L_i)^n \, \nu_i\left(X \setminus \bigcup \mathcal{D}_i\right) \,\leq\, 2^{-i}.
        \end{equation}
        Since $\sum_{i \in \N} 2^{-i} < \infty$, bounds \eqref{eq:bound1} and \eqref{eq:bound2} and the Borel--Cantelli lemma imply that for $\mu$-almost every point $x \in X$, there exists $j \in \N$ such that for all $i \geq j$, we have $x \in \bigcup \mathcal{D}_i$ and $x \notin \bigcup \mathcal{Q}_i$. We claim that every such point $x$ is in $\bigcup \mathcal{R}$. Indeed, since $x \in \bigcup \mathcal{D}_j$, there exists a set $C \in \mathcal{D}_j$ with $x \in C$. As $x \notin \mathcal{Q}_i$ for all $i > j$, we have $C \in \mathcal{R}$, so $x \in \bigcup \mathcal{R}$, as desired.
    \end{scproof}


    \subsection{Rectangular toasts and Baire category}




    Gao, Jackson, Krohne, and Seward showed \cite[Thm.~1.6]{gao2022forcing} that the loss of a null set in Proposition~\ref{prop:MeasurableWaffle} is necessary, that is, a free Borel action of $\Z^n$ may not have a Borel complete rectangular $1$-toast.\footnote{The arguments in \cite{gao2022forcing} are for $\Z^2$, but the generalization to $\Z^n$ is straightforward.} Note that this crucially uses the assumption that $n \geq 2$ (and fails for $n = 1$). Their obstruction is Baire-categorical in nature; they construct a certain profinite action of $\Z^n$ and show that $\bigcup\mathcal{R}$ cannot be comeager for any Borel rectangular $1$-toast $\mathcal{R}$ on this action. 
    Here we observe that this conclusion actually holds for a rather general class of actions.

    \begin{definition}[Minimal actions]\label{def:minimal}
        Let $\Gamma \acts X$ be a continuous action of a group $\Gamma$ on a compact space $X$. The action of $\gamma \in \G$ on $X$ is \emphd{minimal} if for every $x \in X$, the set $\set{\gamma^k \cdot x \,:\, k \in \N}$ is dense. 
    \end{definition}

    In the following, we let $\mathsf{e}_1$, \ldots, $\mathsf{e}_n$ be the standard generators of $\Z^n$.

    \begin{prop}\label{prop:noBaireWaffle}
        Let $\Z^n \acts X$ be a free continuous action on a compact Polish space $X \neq \0$. If the action of $\mathsf{e}_1$ is minimal, then for any Borel rectangular $1$-toast $\mathcal{R} \subseteq \finset{X}$, $\bigcup \mathcal{R}$ is not comeager. 
    \end{prop}

        For a concrete example of such an action, 
        take $X = \mathbb{S}^1$, the circle. For $\alpha \in \R$, let $T_\alpha \colon \mathbb{S}^1 \to \mathbb{S}^1$ be the rotation by the angle $2\pi\alpha$. Fix $\Q$-linearly independent numbers $\alpha_1$, \ldots, $\alpha_n \in \R \setminus \Q$ and define 
        \[
            (k_1, \ldots, k_n) \cdot x \,\defeq\, T_{k_1\alpha_1 \,+\, \cdots \,+\, k_n \alpha_n} (x) \quad \text{for all } (k_1, \ldots, k_n) \in \Z^n \text{ and } x \in \mathbb{S}^1.
        \]
        The $\Q$-linear independence of $\alpha_1$, \ldots, $\alpha_n$ ensures that this action of $\Z^n$ is free, and since $\alpha_1 \notin \Q$, the set $\set{T_{k\alpha_1}(x) \,:\, k \in \N}$ 
        is dense for every  $x \in \mathbb{S}^1$ \cite[Prop.~1.3.3]{Katok}.

        \begin{scproof}[ of Proposition~\ref{prop:noBaireWaffle}]
        
        Note that since the action $\Z^n \acts X$ is continuous, a Borel set $A \subseteq X$ is meager if and only if $\Z^n \cdot A$ is meager. Let $\mathcal{R} \subseteq \finset{X}$ be a Borel rectangular $1$-toast. Since the boundaries of the pieces of $\mathcal{R}$ are at distance at least $2$ from each other, the set $D \defeq X \setminus \partial \mathcal{R}$ satisfies
        \[
            D \,\cup\, \bigcup_{i = 1}^n (\mathsf{e}_i \cdot D) \,\cup\, \bigcup_{i = 1}^n ((-\mathsf{e}_i) \cdot D) \,=\, X.
        \]
        It follows that $D$ is nonmeager. Let $U \subseteq X$ be a nonempty open set such that $U \setminus D$ is meager (which exists by \cite[Prop.~9.8]{AnushDST}). Since the action of $\mathsf{e}_1$ is minimal, we have
        $
            \bigcup_{k \in \N} ((-k\mathsf{e}_1) \cdot U) = X
        $, and the compactness of $X$ yields an integer $N \in \N$ such that
        $
            \bigcup_{0 \leq k \leq N} ((-k\mathsf{e}_1) \cdot U) = X
        $.
        It follows that the set $\bigcup_{0 \leq k \leq N} ((-k\mathsf{e}_1) \cdot D)$ is comeager, and thus its complement, i.e.,
        \[
            M \,\defeq\, X \,\setminus\, \bigcup_{0 \leq k \leq N} ((-k\mathsf{e}_1) \cdot D) \,=\, \big\{y \in X \,:\, (k\mathsf{e}_1) \cdot y \in \partial \mathcal{R} \text{ for all } 0 \leq k \leq N\big\}
        \]
        is meager.  
        Now suppose that $x \in X$ is an $\mathcal{R}$-complete point. By Lemma~\ref{lemma:Wcomplete}, the set $\Z^n \cdot x$ includes a rectangle $K \in \mathcal{R}$ all of whose side lengths exceed $N$. Since $n \geq 2$, there is $y \in \Z^n \cdot x$ such that the boundary of $K$ contains all of the points $y$, $\mathsf{e}_1 \cdot y$, $(2\mathsf{e}_1) \cdot y$, \ldots, $(N\mathsf{e}_1) \cdot y$, which implies that $y \in M$. In other words, the set of all $\mathcal{R}$-complete points is a subset of $\Z^n \cdot M$, and hence it is meager. Therefore, the set $\Z^n \cdot (X \setminus \bigcup\mathcal{R})$ is comeager, so $X \setminus \bigcup\mathcal{R}$ is nonmeager, as desired.
    \end{scproof}

    \subsection{Rectangular toasts and shifts}

    Our next proposition asserts that the shift action $\Z^n \acts \Free(2, \Z^n)$ does \emph{not} satisfy the conclusion of Proposition~\ref{prop:noBaireWaffle}:

    \begin{prop}\label{prop:shiftWaffle}
        For each $q \in \N$, there exists a Borel square $q$-toast $\mathcal{R} \subseteq \finset{\Free(2,\Z^n)}$ 
        such that the set $\bigcup \mathcal{R}$ is comeager.
    \end{prop}
    \begin{scproof}
        Define $\rho \colon \Free(2, \Z^n) \to 2$ by $\rho(x) \defeq x(\bm{0})$. Call a square $K$ in $\Free(2, \Z^n)$ \emphd{safe} if
        \begin{itemize}
            \item the side length of $K$ is at least $q$, and
        
            \item for all $x \in X$, if $x \in \partial K$, then $\rho(x) = 1$, and if $0 < \dist(x, \partial K) \leq q$, then $\rho(x) = 0$.
        \end{itemize}
        Let $\mathcal{R}$ be the set of all safe squares. We claim that $\mathcal{R}$ is a square $q$-toast such that $\bigcup \mathcal{R}$ is comeager, as desired. If $K$, $L \in \mathcal{R}$ are two distinct squares, then $\dist(\partial K, \partial L) > q$, because $\rho$ is equal to $1$ on $\partial K$ and $0$ on all the elements of $X \setminus \partial L$ at distance at most $q$ from $\partial L$. By Remark~\ref{remk:only_bdry_matters}, it follows that $\mathcal{R}$ is a square $q$-toast. It remains to check that $\bigcup \mathcal{R}$ is comeager.

        Note that a point $x \in \Free(2, \Z^n)$ belongs to $\bigcup \mathcal{R}$ if and only if there exists a square $R \subset \Z^n$ in $\Z^n$ with side length at least $q$ such that $R \ni \bm{0}$ and for all $\gamma \in \Z^n$,
        \begin{itemize}
            \item if $\gamma \in \partial R$, then $x(\gamma) = 1$, and
            \item if $0 < \dist(\gamma, \partial R) \leq q$, then $x(\gamma) = 0$.
        \end{itemize}
        In particular, the set $\bigcup \mathcal{R}$ is open, since the existence of such a square $R$ is witnessed by the restriction of $x$ to a finite subset of $\Z^n$. We claim $\partial \mathcal{R}$ is dense. Indeed, the topology on $\Free(2, \Z^n)$ is generated by clopen sets of the form
        \[
            U_{S, \phi} \,\defeq\, \set{x \in \Free(2, \Z^n) \,:\, \rest{x}{S} = \phi},
        \]
        where $S \subset \Z^n$ is finite and $\phi \colon S \to 2$. Given such a set $U_{S, \phi}$, we need to show $U_{S, \phi} \cap \bigcup \mathcal{R} \neq \0$. 
        Take an arbitrary square $R$ 
        in $\Z^n$ with side length at least $q$, $\set{\bm{0}} \cup S \subset R$, and $\dist(S, \partial R) > q$ \ep{which exists because $S$ is a finite set}, and define $x \colon \Z^n \to 2$ by
        \[
            x(\gamma) \,\defeq\, \begin{cases}
                \phi(\gamma) &\text{if } \gamma \in S, \\
                1 &\text{if } \gamma \in \partial R, \\
                0 &\text{otherwise}.
            \end{cases}
        \]
        Then $x \in U_{S, \phi}$ and $R \cdot x$ is a safe square containing $x$, so $x \in U_{S, \phi} \cap \bigcup \mathcal{R}$, as desired. This shows that $\bigcup \mathcal{R}$ is a dense open---hence comeager---set, and the proof is complete.
    \end{scproof}

    \subsection{Rectangular toasts and finitary factors of i.i.d.}

    Before analyzing rectangular toasts from the point of view of finitary factors of i.i.d., it will be helpful to address a certain technical subtlety in Definition~\ref{defn:finitary}. Notice that if $f \colon \Free([0,1],\G) \to \Lambda$ is equal to a finitary function almost everywhere, then $f$ itself need not be finitary; however, such $f$ must have a somewhat weaker property that we call being \emph{weakly finitary}. To define this notion, we first need to introduce a useful piece of terminology.
    
    Suppose $x \in \Free([0,1],\Gamma)$ and $D \subseteq \Gamma$. The product measure $\Leb^{\G \setminus D}$ on $[0,1]^{\Gamma \setminus D}$ then naturally gives rise to a probability measure on the set of all points $x' \in \Free([0,1],\G)$ that agree with $x$ on $D$. 
    When in the sequel we use the phrase ``\emphd{for almost all $x'$ that agree with $x$ on $D$},'' the words ``almost all'' refer to this measure.


    \begin{definition}[Weakly finitary functions]\label{defn:weakly_finitary}
        Let $\G$ be a finitely generated group and let $\Lambda$ be a finite set. 
        A measurable function $f \colon \Free([0,1],\G) \to \Lambda$ is \emphd{weakly finitary} if for almost every $x \in \Free([0,1],\G)$, there exists a finite set $D_x \subseteq \G$ such that $f(x) = f(x')$ for almost all $x' \in \Free([0,1],\G)$ that agree with $x$ on $D_x$. 
        %
    \end{definition}

    Definition~\ref{defn:weakly_finitary} differs from Definition~\ref{defn:finitary} in that $f(x')$ must equal $f(x)$ for \emph{almost all} $x'$ agreeing with $x$ on $D_x$, and not necessarily for \emph{all} such $x$. A straightforward application of Fubini's theorem shows that if $f \colon \Free([0,1],\G) \to \Lambda$ is weakly finitary (for instance, if $f$ is finitary) and $f' \colon \Free([0,1],\G) \to \Lambda$ agrees with $f$ almost everywhere, then $f'$ is also weakly finitary.

    In practice, when working with a measurable function $f$, it is often convenient to replace $f$ by a \emph{Borel} function $f'$ that agrees with $f$ almost everywhere. It is therefore useful to observe 
    that for weakly finitary $f$, such $f'$ can be chosen to be {finitary} \ep{in the strong sense of Definition~\ref{defn:finitary}}:
    
    \begin{prop}\label{prop:Borel_finitary}
        Let $\G$ be a finitely generated group and let $\Lambda$ be a finite set. If $f \colon \Free([0,1],\G) \to \Lambda$ is weakly finitary, then there exists a Borel finitary function $f' \colon \Free([0,1],\G) \to \Lambda$ that is equal to $f$ almost everywhere. 
    \end{prop}
    \begin{scproof}
        Upon modifying $f$ on a null set, we may assume that $f$ is Borel (and still weakly finitary). Fix an increasing sequence $(D_k)_{k \in \N}$ of finite subsets of $\G$ such that $\bigcup_{k \in \N} D_k = \G$ and define $f' \colon \Free([0,1],\G) \to \Lambda$ as follows: If there exist $k \in \N$ and \ep{necessarily unique} $\lambda_k \in \Lambda$ such that $f(x') = \lambda_k$ for almost all $x'$ that agree with $x$ on $D_k$, then let $f'(x) \defeq \lambda_k$ for the smallest such~$k$; otherwise, let $f'(x)$ be any fixed element $\lambda \in \Lambda$. The function $f'$ is Borel by \cite[Thm.~17.25]{KechrisDST}. Moreover, since $f$ is weakly finitary, for almost all $x$, there exist $D_k$ and $\lambda_k$ as above, and $f'(x') = f'(x) = \lambda_k$ for all $x'$ that agree with $x$ on $D_k$. Hence, $f'$ is finitary. It is also clear from the construction that $f$ and $f'$ are equal almost everywhere.
     \end{scproof}
    
    The above proposition will not be explicitly invoked in the sequel, but it explains why we can safely replace ``finitary'' with the more robust notion of ``weakly finitary.'' 
    
    After these preliminaries, we are ready to show that the square toasts given by Proposition~\ref{prop:MeasurableWaffle} applied to $\Free([0,1],\Z^n)$ must in a certain sense be non-finitary. 
    In the following, we use $\Leb$ to denote the 
    Lebesgue measure on $[0,1]$. 

    \begin{prop}\label{prop:noFFIID}
        Suppose that $\mathcal{R} \subseteq \finset{\Free([0,1],\Z^n)}$ is a Borel rectangular $1$-toast such that $\Leb^{\Z^n}(\bigcup \mathcal{R}) = 1$. Then 
        the characteristic function of the set $\partial \mathcal{R}$ 
        is not weakly finitary.
    \end{prop}

    We again stress that Proposition~\ref{prop:noFFIID} relies on the assumption that $n \geq 2$.
    
    \begin{scproof}
        Let $f \colon \Free([0,1],\Z^n) \to \set{0,1}$ be the characteristic function of $\partial \mathcal{R}$ and suppose for contradiction that it is weakly finitary.  
        For a finite set $D \subset \Z^n$ and a function $\phi \colon D \to [0,1]$, we say that $\phi$ \emphd{weakly determines} 
        $f$ if the value $f(x)$ is the same for
        almost
        all 
        $x \in \Free([0,1],\Z^n)$ with $\rest{x}{D} = \phi$. Let $B(D) \subseteq [0,1]^D$ be the set of all $\phi \colon D \to [0,1]$ that weakly determine $f$, and for $i \in \set{0,1}$, let $B_i(D)$ be the set of all $\phi \in B(D)$ such that the common value of $f$ on almost all points extending $\phi$ is $i$. The sets $B(D)$, $B_0(D)$, $B_1(D)$ are $\Leb^D$-measurable by Fubini's theorem.

        Since $f$ is weakly finitary, 
        for almost all $x \in \Free([0,1],\Z^n)$, there is a finite set $D_x \subset \Z^n$ such that $\rest{x}{D_x}$ weakly determines $f$ and, furthermore,
        \[\rest{x}{D_x} \,\in\, B_{f(x)}(D_x).\] Since the boundaries of the pieces of $\mathcal{R}$ are a distance at least $2$ apart, 
        $f(x) = 0$ for a set of points $x \in \Free([0,1],\Z^n)$ of positive measure. As there are only countably many finite subsets of $\Z^n$, 
        Fubini's theorem shows that there exists a finite set $D \subset \Z^n$ such that $\epsilon \defeq \Leb^D(B_0(D)) > 0$.

        We shall employ the following simple observation:

        \begin{claim}\label{claim:disjointify}
            If $S \subset \Z^n$ is a finite set, then there exists a subset $S' \subseteq S$ of size $|S'| \geq |S|/|D|^2$ such that the sets $D + \gamma$ for $\gamma \in S'$ are pairwise disjoint.
        \end{claim}
        \begin{claimproof}
            Let $S' \subseteq S$ be a maximal subset such that the sets $D + \gamma$ for $\gamma \in S'$ are pairwise disjoint. Then for each $\gamma \in S \setminus S'$, there is $\gamma' \in S'$ with $(D + \gamma) \cap (D + \gamma') \neq \0$, i.e., $\gamma \in D - D + \gamma'$. On the other hand, for each $\gamma' \in S'$, there are fewer than $|D - D| \leq |D|^2$ elements $\gamma \in S \setminus S'$ in $D - D + \gamma'$ (note that $\gamma'$ itself should not be counted). Thus, $(|D|^2 - 1) |S'| \geq |S| - |S'|$, i.e., $|S'| \geq |S|/|D|^2$. 
        \end{claimproof}

        Let $x \in \Free([0,1],\Z^n)$ be a random point sampled from the distribution $\Leb^{\Z^n}$ 
        and let $\Xi$ be the random variable equal to the number of rectangles in $\mathcal{R}$ that contain $x$. We claim that the expectation of $\Xi$ is finite. 
        Consider a fixed rectangle $R$ in $\Z^n$ with side lengths $N_1$, \ldots, $N_n$ such that $\bm{0} \in R$. If $R \cdot x \in \mathcal{R}$, then for all $\delta \in \partial R$, we have $f(\delta \cdot x) = 1$, and hence, {with probability 1}, $\rest{(\delta \cdot x)}{D} \notin B_0(D)$. Let $S \subseteq \partial R$ be a set of size $|S| \geq |\partial R|/|D|^2$ such that the sets $D + \delta$ for $\delta \in S$ are pairwise disjoint, which exists by Claim~\ref{claim:disjointify}. Then the random events $\rest{(\delta \cdot x)}{D} \notin B_0(D)$ for $\delta \in S$ are independent, so
        \begin{align*}
            \mathbb{P}[R \cdot x \in \mathcal{R}] \,\leq\, \mathbb{P}\big[\rest{(\delta \cdot x)}{D} \notin B_0(D) \text{ for all } \delta \in S\big] \,&=\, (1 - \epsilon)^{|S|} \\
            &\leq\, (1 - \epsilon)^{|\partial R|/|D|^2} \,\leq\, \prod_{i=1}^n (1-\epsilon)^{N_i/|D|^2},
        \end{align*}
        where in the last step we use the crude bound $|\partial R| \geq \sum_{i=1}^N N_i$. (This is where the assumption $n \geq 2$ is invoked.) Note that there are $\prod_{i=1}^n N_i$ many rectangles in $\Z^n$ with side lengths $N_1$, \ldots, $N_n$ that include $\bm{0}$. Therefore, by the linearity of expectation,
        \begin{align*}
            \mathbb{E}[\Xi] \,=\, \sum_{R} \mathbb{P}[R \cdot x \in \mathcal{R}] \, \leq\, \sum_{N_1, \ldots, N_n} \prod_{i=1}^n \left(N_i (1-\epsilon)^{N_i /|D|^2}\right) \,=\, \left(\sum_{N = 1}^\infty N (1-\epsilon)^{N/|D|^2}\right)^n \,<\, \infty, 
        \end{align*}
        as desired. (Here the first summation is over all rectangles $R$ in $\Z^n$ containing $\bm{0}$, and the second one is over all integers $N_1$, \ldots, $N_n \geq 1$.)
        
        It follows that $\Xi < \infty$ almost always, i.e., almost every point $x \in \Free([0,1],\Z^n)$ is contained in finitely many rectangles in $\mathcal{R}$. By Lemma~\ref{lemma:Wcomplete}, this means that the set of all $\mathcal{R}$-complete points is null, i.e., $\Leb^{\Z^n}(\Z^n \cdot (X \setminus \bigcup \mathcal{R})) = 1$. Since the action $\Z^n \acts \Free([0,1], \Z^n)$ is measure-preserving, we conclude that the set $X \setminus \bigcup \mathcal{R}$ has positive measure, which is a contradiction.
    \end{scproof}

    \subsection{Computable rectangular toasts}\label{subsec:waffle_comput}

    Now we turn our attention to the computability-theoretic aspects of rectangular toasts. We say that a set $\mathcal{R} \subseteq \finset{\N}$ is \emphd{computable} if there is an algorithm that, given a tuple $(x_1, \ldots, x_k)$ of natural numbers, decides if $\set{x_1, \ldots, x_k} \in \mathcal{R}$. 
    %
    %
    The following proposition was pointed out to us by Riley Thornton in a personal communication as a potential tool for resolving Question \ref{q:comp_vs_baire}.

    \begin{prop}\label{prop:compWaffle}
        For every free computable action $\Z^n \acts X$ on a computable set $X \subseteq \N$ and for each $q \in \N$, there exists a computable complete square $q$-toast $\mathcal{R} \subseteq \finset{X}$. 
    \end{prop}
    \begin{scproof}
        Without loss of generality, assume that $X = \N$. We shall describe an algorithm that enumerates a sequence $R_0$, $R_1$, \ldots{} with the following properties:
        \begin{enumerate}[label=\ep{\normalfont\small{R}\arabic*}]
            \item\label{item:induction_large}\label{item:induction_first} each $R_i$ is an $N_i$-square in $X$ for some $N_i \geq i + q$,
            \item\label{item:induction_far} for all $i \neq j$, we have $\dist(\partial R_i, \partial R_j) > q$,
            \item\label{item:induction_cover}\label{item:induction_last} for all $x \in \N$, $x \in R_x$. 
        \end{enumerate}
        The set $\mathcal{R} \defeq \set{R_i \,:\, i \in \N}$ is a desired computable compete square $q$-toast. Indeed, $\mathcal{R}$ is a square $q$-toast by \ref{item:induction_large} and \ref{item:induction_far}, and it is complete by \ref{item:induction_cover}. Furthermore, $\mathcal{R}$ is computable because, given a tuple $(x_1, \ldots, x_k)$ such that $R \defeq \set{x_1, \ldots, x_k}$ is an $N$-square, 
        we can decide if $R$ is in $\mathcal{R}$ by computing $R_0$, \ldots, $R_N$. If $R = R_i$ for some $0 \leq i \leq N$, then $R \in \mathcal{R}$; otherwise, $R \notin \mathcal{R}$ by \ref{item:induction_large}.

        To enumerate $R_0$, $R_1$, \ldots, we proceed as follows. 
        Note that, since the action $\Z^n \acts \N$ is computable, given two rectangles $R$ and $R'$, 
        we can computably determine whether $\dist(\partial R, \partial R') > q$.
        Assume that $R_0$, \ldots, $R_{i-1}$ have already been computed. We then cycle through the integers $j \geq i + q$ until we find one such that $\dist(\partial ([-j, j]^n \cdot i), \partial R_t) > q$ for all $t \in \set{0,\ldots, i-1}$. Note that such $j$ must exist because the set $\set{R_0, \ldots, R_{i-1}}$ is finite. We then let $N_i \defeq 2j$ and $R_i \defeq [-j, j]^n \cdot i$. This construction clearly satisfies properties \ref{item:induction_first}--\ref{item:induction_last}, as desired.
    \end{scproof}

    \section{LCL problems based on rectangular toasts}\label{sec:baking}

    \subsection{The rectangular $q$-toast problem}\label{subsec:qwaffleproblem}

    In view of Propositions~\ref{prop:MeasurableWaffle}, \ref{prop:noBaireWaffle}, \ref{prop:shiftWaffle}, \ref{prop:noFFIID}, and \ref{prop:compWaffle}, our goal now is to reduce the task of finding a complete rectangular $q$-toast for some $q \geq 1$ to solving an LCL problem. The challenge is to ``localize'' the seemingly ``global'' notion of a complete rectangular $q$-toast. We do this in two steps. First, for each \ep{large enough} $q$, we describe an LCL problem $\Waffle(q)$ on $\Z^n$ whose solutions encode rectangular $q$-toasts \ep{there are some technical caveats to this statement that are discussed below}. After that, in \S\S\ref{subsec:p2c}, \ref{subsec:CW}, we explain how to modify the problem to force the rectangular $q$-toast to be {complete}. In the latter step we shall make use of the proper $2$-coloring problem (see Examples~\ref{exmp:coloring} and \ref{exmp:2coloring}), which is typically hard to solve ``globally,'' but admits well-behaved ``local'' solutions.

    To begin with, we convert a toast into a labeling using the following definition:

    \begin{definition}[From $\mathcal{R}$ to $f_\mathcal{R}$]\label{defn:fromWtof}
    Let $\Z^n \acts X$ be a free action. We define the \emphd{outer boundary} of a set $A \subseteq X$ to be $\partial^+ A \defeq \partial(X \setminus A)$, i.e., $x \in \partial^+ A$ if and only if $x \notin A$ and $\dist(x,A) = 1$. Given a $q$-toast $\mathcal{R} \subseteq \finset{X}$ for $q \geq 1$, we let $\partial^+\mathcal{R} \defeq \bigcup_{K \in \mathcal{R}} \partial^+ K$ and define a map $f_\mathcal{R} \colon X \to \RBG$ by 
%
%
        \[
            f_\mathcal{R}(x)\,\defeq\, \begin{cases}
                \mathtt{R} &\text{if } x \in \partial \mathcal{R},\\
                \mathtt{B} &\text{if } x \in \partial^+\mathcal{R},\\
                \mathtt{G} &\text{otherwise}.
            \end{cases}
        \]
    \end{definition}

    \begin{figure}[t]
		\centering
		\includegraphics[scale=0.9]{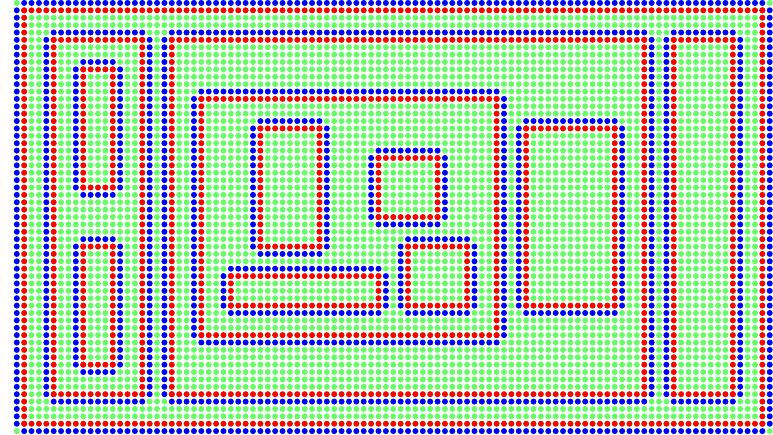}
   \caption{The labeling corresponding to the rectangular toast from Fig.~\ref{fig:waffle}.}\label{fig:waffle_coloring}
	\end{figure}
    
    To phrase Definition~\ref{defn:fromWtof} another way, we color the points on the boundary of each piece \texttt{R}\emph{ed}, the points adjacent to some piece from the outside \texttt{B}\emph{lue}, and the other points \texttt{G}\emph{reen}. 
    This construction is illustrated in Fig.~\ref{fig:waffle_coloring}. Note that, since $q \geq 1$, $\partial \mathcal{R} \cap \partial^+ \mathcal{R} = \0$, so $f_\mathcal{R}$ is well-defined.
    
    Roughly, the solutions to our LCL problem $\Waffle(q)$ will be the $\RBG$-colorings that ``locally'' look as if they came from this construction:

    \begin{definition}\label{def:waffle_lcl}
        Given $q \geq 4$, we define the LCL problem $\Waffle(q) \defeq ([0,2q]^n, \RBG, \mathcal{A})$, where 
        \[ \mathcal{A} \,\defeq\, \big\{\phi \colon [0,2q]^n \to \set{\mathtt{R}, \mathtt{B}, \mathtt{G}} \;:\; \phi = \rest{f_\mathcal{R}}{[0,2q]^n} \text{ for some rectangular $q$-toast $\mathcal{R} \subseteq \finset{\Z^n}$}\big\}. \]
    \end{definition}

    Here and in the sequel, we view $\Z^n$ as acting on itself by translation. Note that in Definition~\ref{def:waffle_lcl}, as well as in the rest of this section, we for convenience assume $q \geq 4$. 
    
    If $\Z^n \acts X$ is a free action, then a function $f \colon X \to \RBG$ is an $\Waffle(q)$-labeling of $X$ if and only if for every $2q$-square $Q$ in $X$, there exists a rectangular $q$-toast $\mathcal{R} \subseteq \finset{X}$ such that $f$ agrees with $f_\mathcal{R}$ on $Q$. In particular, for any rectangular $q$-toast $\mathcal{R} \subseteq \finset{X}$, the map
    $f_\mathcal{R}$ is an $\Waffle(q)$-labeling by definition. We wish to establish a converse to this fact: 
    given an $\Waffle(q)$-labeling $f$, we would like to find a rectangular $q$-toast $\mathcal{R}$ with $f = f_\mathcal{R}$. Unfortunately, this is not always possible. For instance, consider the map $f \colon \Z^2 \to \RBG$ given by 
    \[ f(x,y) \,\defeq\, \begin{cases}
        \mathtt{R} & \text{if } x = 0, \\ 
        \mathtt{B} & \text{if } x = 1, \\ 
        \mathtt{G} & \text{otherwise}.
    \end{cases} \]
    This is an $\Waffle(q)$-labeling of $\Z^2$. Indeed, the labeling around each point of the form $(0,y)$ locally looks as if that point is on the right boundary of a very large rectangle. However, this labeling clearly does not come from a rectangular toast. 
    Notice, though, that $f$ would come from a rectangular toast if we allowed the ``infinite rectangle'' $(-\infty,0] \times (-\infty, \infty)$ to be included. It turns out that, modulo allowing such ``infinite rectangles,'' $\Waffle(q)$-labelings really \emph{do} correspond to rectangular $q$-toasts. 
    
    Before stating this result precisely, we need some extra notation. 
    Let $\mathsf{e}_1$, \ldots, $\mathsf{e}_n$ be the standard generators of $\Z^n$. 
    \label{page:G(X)} Given a free action $\Z^n \acts X$, we let $G(X)$ be the graph with $V(G(X)) \defeq X$ and
    \[
        E(G(X)) \,\defeq\, \{\set{x,y} \in [X]^2 \,:\, \dist(x,y) = 1\} \,=\, \set{\set{x, \mathsf{e}_i \cdot x} \,:\, x\in X, \, 1 \leq i \leq n}.
    \]
    Note that $\dist$ is precisely the path metric on $G(X)$. 
    In the special case when $X = \Z^n$ and $\Z^n \acts \Z^n$ is the translation action of $\Z^n$ on itself, we let $G_n \defeq G(\Z^n)$ and observe that $G_n$ is the Cayley graph $\Cay(\Z^n, \set{\mathsf{e}_1, \ldots, \mathsf{e}_n})$. 
    Note that for every free action $\Z^n \acts X$, every connected component of $G(X)$ is isomorphic to $G_n$. 
    For any graph $G$ and a set $U \subseteq V(G)$, we let $G[U]$ be the \emphd{subgraph} of $G$ \emphd{induced} by $U$, i.e., the graph with vertex set $U$ and edge set $E(G) \cap [U]^2$. 
    An \emphd{interval} in $\Z$ is a set of the form $[a,b]$, $[a, +\infty) \defeq \set{c \in \Z \,:\, a \leq c}$, $(-\infty,b) \defeq \set{c \in \Z \,:\, c \leq b}$, or $(-\infty, +\infty) \defeq \Z$, where $a \leq b$ are integers.

    \begin{lemma}\label{lem:infinite_rectangles}
        Let $f \colon \Z^n \to \RBG$ be an $\Waffle(q)$-labeling of $\Z^n$. Consider an arbitrary connected component $C$ of the {induced sub}graph $G_n[f\inv(\mathtt{R})]$. Let $R \defeq V(C)$ and let $B$ be the set of all points in $f\inv(\mathtt{B})$ having a neighbor in $R$. Then $R = \partial S$ and $B = \partial^+ S$ for some set $S$ of the form $S = I_1 \times \cdots \times I_n$, where $I_1$, \ldots, $I_n$ are {\ep{possibly infinite}} intervals in $\Z$.
        %
    \end{lemma}
    \begin{scproof}
    Let $E \subseteq E(G_n)$ be the set of all edges with one endpoint in $B$ and the other one in $R$; we call such edges \emphd{bad}. Let $H$ be the subgraph of $G_n$ obtained by deleting the bad edges; that is, $V(H) \defeq \Z^n$ and $E(H) \defeq E(G_n) \setminus E$.   

    The following claim is the key to our argument:
    
    \begin{claim}[$R$ and $B$ are separated in $H$]\label{claim:2comps}
        The set $R$ is contained in a single component of the graph $H$, and that component is disjoint from $B$ \ep{hence $H$ has at least two components}.
    \end{claim}

   With a little more work, one can show that $H$ has \emph{exactly two} components, one containing $R$ and the other one $B$; we will not need this, however.

    
    Intuitively, Claim~\ref{claim:2comps} can be thought of as a discretized version of the Jordan--Brouwer separation theorem from topology, which says that a connected closed topological $(n-1)$-manifold embedded in $\R^n$ separates $\R^n$ into exactly two connected sets. 
    This intuition informs our proof of the claim, which uses a discretized version of the mod 2 intersection theory that appears in some proofs of the Jordan--Brouwer theorem.

    \begin{claimproof}[Proof of Claim~\ref{claim:2comps}]
        It suffices to prove the following subclaim: \textsl{Every cycle in $G_n$ contains an even number of bad edges}. This implies the claim as follows.
        First, it is clear that $R$ is contained in a single component of $H$ since $R = V(C)$, $C$ is connected, and no edge in $C$ is bad. Hence, if the claim failed, there would be a path $P$ from $R$ to some $y \in B$ using no bad edges. By the definition of $B$, there must be some $x \in R$ adjacent to $y$ in $G_n$. Since $R$ is in a single component of $H$ we may assume $x$ is the starting point of $P$. Adding the edge $\set{y,x}$ to $P$ yields a cycle in $G_n$ with a single bad edge \ep{namely $\set{y,x}$}, contradicting the subclaim. 

        It is well known and easy to see that the binary cycle space of $G_n$ is spanned by the $4$-cycles. \ep{For example, this follows from the fact that the relations $\mathsf{e}_i + \mathsf{e}_j = \mathsf{e}_j + \mathsf{e}_i$ for $1 \leq i < j \leq n$ give a presentation of $\Z^n$.} 
        Thus it suffices to verify the subclaim for 4-cycles.

        So, let $\Sigma$ be a 4-cycle in $G_n$. If none of its vertices are in $R$, then $\Sigma$ contains no bad edges. Now suppose that $V(\Sigma) \cap R \neq \0$ and take an arbitrary vertex $w \in V(\Sigma) \cap R$. Since $f$ is an $\Waffle(q)$-labeling, 
        there exists a rectangular $q$-toast $\mathcal{R}$ such that $f$ and $f_\mathcal{R}$ agree on $V(\Sigma)$. Since $f(w) = \mathtt{R}$, there is a rectangle $K \in \mathcal{R}$ such that $w \in \partial K$. Moreover, as $q \geq 4$, 
        $\Sigma$ cannot meet any other rectangle in $\mathcal{R}$. Thus, an edge of $\Sigma$ is bad if and only if it has one endpoint in $K$ and the other one not in $K$. Since the cycle starts and ends at $w \in K$, the number of such edges must be even.
    \end{claimproof}

    Let $S$ be the vertex set of the component of $H$ that contains $R$. 
    We claim that $\partial S = R$. Indeed, any $x \in \partial S$ must be an endpoint of a bad edge, and since $x \notin B$ by Claim~\ref{claim:2comps}, it follows that $x \in R$. For the other direction, if $x \in R$, then, by the definition of $\Waffle(q)$, there is a point $y \in B$ adjacent to $x$. By Claim~\ref{claim:2comps}, $y \notin S$, so $x \in \partial S$, as desired. A similar argument shows that $\partial^+ S = B$.
    

    It remains to show that $S$ is a product of intervals. We will establish this by showing that $S$ is ``locally convex'' in the following sense:

    \begin{claim}\label{claim:convex}
        Let $(w,x,y,z)$ be a 4-cycle in $G_n$ with $w$, $x$, $y \in S$. Then $z \in S$. 
    \end{claim}
    \begin{claimproof}
        Since $\partial S = R$ and $\partial^+ S = B$, if this fails
        we must have $z \in B$ while $w$, $y \in R$. However, we claim that if $(w,x,y,z)$ is a 4-cycle in $G_n$ with $f(w) = f(y) = \mathtt{R}$, then $f(z) \neq \mathtt{B}$. 
        Indeed, let $\mathcal{R}$ be a rectangular $q$-toast such that $f$ and $f_\mathcal{R}$ agree on the vertices $w$, $x$, $y$, $z$ and their neighbors. 
        Then there are $K$, $K' \in \mathcal{R}$ such that $w \in \partial K$ and $y \in \partial K'$.  Since $\dist(w,y) = 2$ and $q \geq 4$, 
        we must have $K = K'$. As $w$, $y \in K$ and $K$ is a rectangle, it follows that $z \in K$ as well. 
        For $f(z)$ to be $\mathtt{B}$, then, there would need to be some other rectangle $K'' \neq K$ in $\mathcal{R}$ with $\dist(z,K'') = 1$. But then $\dist(y, \partial K'') = 2$, which contradicts the definition of a $q$-toast as $y \in \partial K$. 
    \end{claimproof}

    It remains to show that any $G_n$-connected set $S$ satisfying Claim~\ref{claim:convex} is a product of intervals. For any $x$, $y \in \Z$, let $I(x,y) \defeq [\min(x,y), \max(x,y)]$ be the interval between $x$ and $y$, and for any $x = (x_1,\ldots,x_n)$, $y = (y_1,\ldots,y_n) \in \Z$, 
    let $I(x,y) \defeq I(x_1,y_1) \times \cdots \times I(x_n,y_n)$. It suffices to show that for all $x$, $y \in S$, $I(x,y) \subseteq S$.
    

    We proceed by induction on the length of an $xy$-path $P$ in $G_n[S]$. In the base case $x = y$, so $I(x,y) = I(x,x) = \set{x}$ and there is nothing to show. For the inductive step, let $z$ be the penultimate vertex of $P$. Then $I(x,z) \subseteq S$ by the inductive hypothesis.
    Without loss of generality, we may assume that 
    $y = \mathsf{e}_1 + z$. 
    Then
    \[
        I(x,y) \,\subseteq\, I(x,z) \cup \big(\{y_1\} \times I(x_2,y_2) \times \cdots \times I(x_n,y_n)\big).
    \]
    Calling the second set in this union $T$, it remains to show that $T \subseteq S$. To this end, take any $w \in T$. Noting that $T$ is $G_n$-connected, we argue that $w \in S$ by induction on the length of a $yw$-path $Q$ in $G_n[T]$. In the base case $w = y \in S$ as desired. For the inductive step, let $v$ be the penultimate vertex of $Q$. 
    Then $(v,w,w - \mathsf{e}_1, v - \mathsf{e}_1)$ is a 4-cycle in $G_n$. Since $v \in S$ by the inductive hypothesis and $w - \mathsf{e}_1$, $v - \mathsf{e}_1$ are in $I(x,z) \subseteq S$, we have $w \in S$ by Claim~\ref{claim:convex}, as desired.
    \end{scproof}


    We can now state the main result of this section: if an $\Waffle(q)$-labeling of $\Z^n$ does not come from a rectangular $q$-toast, then it must contain infinite rays of red, blue, and green points. 

    \begin{lemma}\label{lem:waffle_dichotomy}
        Let $f \colon \Z^n \to \RBG$ be an $\Waffle(q)$-labeling of $\Z^n$. Exactly one of the following holds:
        \begin{enumerate}[label=\ep{\normalfont\small{RT}\arabic*}]
            \item\label{item:waffle} $f = f_\mathcal{R}$ for a rectangular $q$-toast $\mathcal{R}$, or 
            \item\label{item:infinite_rays} there exist $x \in \Z^n$ and orthogonal
            $\mathsf{d}$, $\mathsf{d}' \in \{\pm \mathsf{e}_1,\ldots, \pm \mathsf{e}_n\}$
            such that for all $i \in \N$,
            \[
                f(i\mathsf{d}+x) \,=\, \mathtt{R}, \quad f(i\mathsf{d} + \mathsf{d}' + x) \,=\,\mathtt{B}, \quad \text{and} \quad f(i\mathsf{d} + 2\mathsf{d}' + x) \,=\, \mathtt{G}.
            \] 
        \end{enumerate}
    \end{lemma}
    \begin{scproof}
        Options \ref{item:waffle} and \ref{item:infinite_rays} are mutually exclusive since the rectangles in a rectangular $q$-toast are finite.
        Let $\mathcal{C}$ be the set of all components of the graph $G_n[f\inv(\mathtt{R})]$.
        For each $C \in \mathcal{C}$, fix a set $S_C$ such that $\partial S_C = V(C)$ and $\partial^+S_C$ is the set of all blue points with a neighbor in $V(C)$,
        given by Lemma~\ref{lem:infinite_rectangles}. Suppose for some $C \in \mathcal{C}$, $S \defeq S_C$ is not a \ep{finite} rectangle.
        Since $S$ is a product of intervals at least one of which is infinite and $\partial S \neq \0$, we can choose $x \in \partial S$ and $\mathsf{d} \in \{\pm \mathsf{e}_1,\ldots, \pm \mathsf{e}_n\}$ so that the set $L \defeq \{i\mathsf{d} + x \,:\, i \in \N\}$ is contained in $\partial S \subseteq f\inv(\mathtt{R})$. Furthermore, there is some $\mathsf{d}' \in \{\pm \mathsf{e}_1,\ldots, \pm \mathsf{e}_n\}$ orthogonal to $\mathsf{d}$ such that $L + \mathsf{d}' \subseteq \partial^+ S \subseteq f\inv(\mathtt{B})$. 
        The definition of $\Waffle(q)$ implies that for any $y \in \Z^n$, if $f(y) = \mathtt{R}$ and $f( \mathsf{d}' + y) = \mathtt{B}$, then $f(2\mathsf{d}' + y) = \mathtt{G}$. 
        Hence, in this case \ref{item:infinite_rays} holds.
        
        
        
        Now assume $S_C$ is a finite rectangle for all $C \in \mathcal{C}$ and 
        let $\mathcal{R} \defeq \set{S_C \,:\, C \in \mathcal{C}}$. 
        Then $\partial \mathcal{R} = f\inv(\mathtt{R})$, and, 
        since every blue vertex has a red neighbor, $\partial^+ \mathcal{R} = f\inv(\mathtt{B})$. Therefore, $f = f_\mathcal{R}$. To confirm \ref{item:waffle}, it remains to verify that $\mathcal{R}$ is a rectangular $q$-toast.
        
        Note first that all side lengths of each rectangle in $\mathcal{R}$ are at least $q$. Indeed, if, say, the $i$-th side length of $K \in \mathcal{R}$ is $\ell \leq q - 1$, then there exists a point $x \in \partial^+ K$ such that
        \[
            f(x) \,=\, \mathtt{B}, \quad f(\mathsf{e}_i + x) \,=\, f(2\mathsf{e}_i + x) \,=\, \cdots \,=\, f((\ell+1) \mathsf{e}_i + x) \,=\, \mathtt{R}, \quad f((\ell+2)\mathsf{e}_i + x) \,=\, \mathtt{B},
        \]
        which cannot occur in an $\Waffle(q)$-labeling since $\ell+2 \leq q+1 < 2q$.
        
        Now suppose that there are $K \neq K' \in \mathcal{R}$ with $\dist(\partial K, \partial K') \leq q$. Let $x \in \partial K$, $x' \in \partial K'$ witness this. Then there is a square of side length $q$, call it $Q$, containing both $x$ and $x'$. Using the definition of $\Waffle(q)$, let $\mathcal{R}'$ be a rectangular $q$-toast 
        such that $f$ agrees with $f_{\mathcal{R}'}$ on $Q$, and let $L$, $L' \in \mathcal{R}'$ be such that $x \in \partial L$ and $x' \in \partial L'$. Since $\dist(x,x') \leq q$ and $\mathcal{R}'$ is a $q$-toast, we must have $L = L'$. Since every side length of $L$ is at least $q$ by the definition of a rectangular $q$-toast, the graph $G_n[\partial L \cap Q]$ is connected. But $x' \in \partial L \cap Q \subseteq f\inv(\mathtt{R})$, so $x'$ is in the same component of $G_n[f\inv(\mathtt{R})]$ as $x$, a contradiction.
    \end{scproof}

    \subsection{Partial 2-colorings} \label{subsec:p2c}

    At this point, we need to modify the problem $\Waffle(q)$ in order to (\emph{a}) avoid case \ref{item:infinite_rays} in Lemma \ref{lem:waffle_dichotomy}, and (\emph{b}) force the rectangular $q$-toast $\mathcal{R}$ in case \ref{item:waffle} to be complete. 
    We will achieve these goals by requiring the green points 
    to additionally carry a proper $2$-coloring. 
    \ep{A similar idea appeared in the proof of \cite[Thm.~6.1]{CS8}.} In this section we shall develop tools that will be employed to work with such $2$-colorings later. 

    \begin{remk}
        It turns out that the $\Z^n$-actions we shall work with automatically avoid case \ref{item:infinite_rays} away from a null/meager set, so the $2$-coloring is only truly needed to guarantee completeness.
    \end{remk}

    \begin{definition}[Partial $2$-colorings]
        Let $\Z^n \acts X$ be a free action. A \emphd{partial 2-coloring} of $X$ is a function $c \colon X \to \set{0,1,\blank}$ such that if $\dist(x, y) = 1$ and $c(x)$, $c(y) \in \set{0,1}$, then $c(x) \neq c(y)$. We say that a point $x \in X$ is \emphd{uncolored} if $c(x) = \blank$ and \emphd{colored} otherwise. The set of all colored points is denoted by $\dom(c)$. 
    \end{definition}

        Note that a partial $2$-coloring is a solution to a certain LCL problem with labels $\set{0,1,\blank}$. 
    
    As we saw in Example \ref{exmp:2coloring}, although the standard Cayley graph of $\Z$ admits a proper 2-coloring, in the definable (e.g., Borel, measurable, etc.) context proper 2-colorings may fail to exist. The same is true for $\Z^n$ with $n \geq 2$, and, as we shall see, 
    similar arguments impose restrictions on 
    the domain of any {Borel} {partial} 2-coloring. 
    Intuitively, Propositions~\ref{prop:fiid_partial2color} and \ref{prop:baire_partial2color} imply that under suitable assumptions on a $\Z^n$-action $\Z^n \acts X$, the domain of a Borel partial $2$-coloring $c \colon X \to \set{0,1,\blank}$ must ``shatter'' the graph $G(X)$ into several smaller components; for instance, it is not possible for $\dom(c)$ to have a connected intersection with every $\Z^n$-orbit.  We will also show in Proposition~\ref{prop:component_finite_2color} that if a Borel set $U \subseteq X$ ``shatters'' $X$ into sufficiently small pieces---namely if $G(X)[U]$ is component-finite---then we can find a Borel partial $2$-coloring with domain $U$. 


    \begin{definition}[$c$-Separated points]
        Let $\Z^n \acts X$ be a free action and let $c \colon X \to \set{0,1,\blank}$ be a partial $2$-coloring. We say that points $x$, $y \in X$ are \emphd{$c$-separated} if $x$ and $y$ belong to $\dom(c)$ and lie in different components of the graph $G(X)[\dom(c)]$.
    \end{definition}

    Recall that $\Leb$ denotes the Lebesgue measure on $[0,1]$, {and $\Leb^{\Z^n}$ the product measure. }

    \begin{prop}\label{prop:fiid_partial2color}
        Let $c \colon \Free([0,1],\Z^n) \to \set{0,1,\blank}$ be a measurable partial $2$-coloring such that 
        $\Leb^{\Z^n}(\dom(c)) > 0$.
        Then for every $\mathsf{d} \in \Sn$, there exist $x \in \dom(c)$ and $k \in \N$ such that the points $x$ and $k\mathsf{d} \cdot x$ are $c$-separated. 
    \end{prop}
    \begin{scproof}
        Without loss of generality, $c\inv(0)$ has positive measure. By the regularity of measure \cite[Thm.~17.10]{KechrisDST}, 
        there exist a finite set $S \subset \Z^n$ and nonempty open sets $U_s \subseteq [0,1]$ for each $s \in S$ such that, if $x \in \Free([0,1],\Z^n)$ is chosen at random according to $\Leb^{\Z^n}$, then
        \[
            \mathbb{P}\big[c(x) = 0 \,\mid\, x(s) \in U_s \text{ for all } s \in S\big] \,\geq\, \frac{2}{3}.
        \]
        Let $k \in \N$ be odd and large enough so that $(k\mathsf{d} + S) \cap S = \0$. Now consider the nonempty open set $W$ of all $x \in \Free([0,1],\Z^n)$ such that for each $s \in S$, we have $x(s)$, $x(k\mathsf{d} + s) \in U_s$. By the union bound, the probability that $c(x) = c(k\mathsf{d} \cdot x) = 0$ given $x \in W$ is at least $1/3$. It follows that there exists some $x \in W$ with $c(x) = c(k\mathsf{d} \cdot x) = 0$. We claim that this point $x$ is as desired. 

        If not, the graph $G(\Free([0,1],\Z^n))$ has a path $P$ from $x$ to $k\mathsf{d} \cdot x$ contained in $\dom(c)$. Note that every component of $G(\Free([0,1],\Z^n))$ is isomorphic to the bipartite graph $G_n$. Thus, since $k$ is odd, every path in $G(\Free([0,1],\Z^n))$ from $x$ to $kd \cdot x$, including $P$, has an odd length. 
        But this is impossible since $c(x) = c(k\mathsf{d} \cdot x) = 0$ and the colors $0$ and $1$ must alternate between $0$ and $1$ along $P$. 
    \end{scproof}

    \begin{prop}\label{prop:baire_partial2color}
        Let $\Z^n \acts X$ be a free continuous action on a Polish space $X$. Suppose that for some $\mathsf{d} \in \Sn$, there is $y \in X$ such that the set $\{2i\mathsf{d} \cdot y \,:\, i \in \N\}$ is dense. Then, if $c \colon X \to \set{0,1,\blank}$ is a Baire-measurable partial $2$-coloring and $\dom(c)$ is nonmeager, then there are $x \in \dom(c)$ and $k \in \N$ such that the points $x$ and $k\mathsf{d} \cdot x$ are $c$-separated. 
        %
    \end{prop}
    \begin{scproof}
        %
        Without loss of generality, $c\inv(0)$ is nonmeager. By \cite[Prop.~9.8]{AnushDST}, there is a nonempty open set $U \subseteq X$ such that $c(x) = 0$ for a comeager set of points $x \in U$. Let $U' \defeq \bigcup_{i \in \Z} ((2i+1)\mathsf{d} \cdot U)$ \ep{note that here the union is over all $i \in \Z$, not just $i \in \N$}. By construction, if $x \in U'$, then $2i\mathsf{d} \cdot x \in U'$ for all $i \in \Z$. Since $U'$ is open and $\{2i\mathsf{d} \cdot y \,:\, i \in \N\}$ is dense, there is $i \in \N$ such that $2i\mathsf{d} \cdot y \in U'$. 
        This implies that in fact $2i\mathsf{d} \cdot y \in U'$ for all $i \in \Z$, so 
     $U'$ is dense. 
        It follows that $U \cap U' \neq \0$, i.e.,  
        $U \cap ((2i+1)\mathsf{d} \cdot U) \neq \0$ for some $i \in \N$. 

        For a comeager set of points $x \in U \cap ((2i+1)\mathsf{d} \cdot U)$, we have $c(x) = c(-(2i+1)\mathsf{d} \cdot x) = 0$. Fix any such $x$ and let $k \defeq 2i + 1$ (if $k < 0$, swap the roles of $x$ and $k\mathsf{d} \cdot x$). The rest of the argument is the same as in the proof of Proposition~\ref{prop:fiid_partial2color}.
    \end{scproof}


    
    On the other hand, we have the following positive results:

    \begin{prop}\label{prop:component_finite_2color}
        Let $\Z^n \acts X$ be a free Borel action on a standard Borel space and let $U \subseteq X$ be a Borel set such that the graph $G(X)[U]$ is component-finite. Then $X$ admits a Borel partial $2$-coloring $c$ with $\dom(c) = U$. 
    \end{prop}
    \begin{scproof}
        This is a standard fact in descriptive combinatorics and an instance of the general principle that Borel combinatorics trivializes on component-finite graphs. In a sentence, we can find a Borel set meeting each component of $G(X)[U]$ exactly once and define $c(x)$ for $x \in U$ to be the distance mod $2$ from $x$ to the chosen point in its component. 
        See \cites[\S2.2]{BerFelixASI}[\S5.3]{Pikh_survey} for details. 
    \end{scproof}



    \begin{prop}\label{prop:comp_2color}
        Let $\Z^n \acts X$ be a free computable action on a computable set $X \subseteq \N$ and let $U \subseteq X$ be a computable subset such that the graph $G(X)[U]$ is component-finite. Then $X$ admits a computable partial $2$-coloring $c$ with domain $U$. 
    \end{prop}
    \begin{scproof}
        We apply the same idea as in the proof of Proposition~\ref{prop:component_finite_2color}. The function $h \colon U \to \finset{U}$ sending a point to the vertex set of its component in $G(X)[U]$ is clearly computable. Now, given $x \in U$, we let $c(x)$ be the distance mod $2$ from $x$ to, say, $\min h(x)$.
    \end{scproof}

    \subsection{The colored rectangular $q$-toast problem}\label{subsec:CW}

    We are now ready to define the LCL problem $\Pi$ that witnesses Theorem~\ref{theo:main}. 

    \begin{definition}
        Let $\pi \colon \set{\mathtt{R},\mathtt{B}, 0, 1} \to \RBG$ and $\tau \colon \set{\mathtt{R},\mathtt{B}, 0, 1} \to \set{0,1,\blank}$ be given by
        \begin{align*}
            \begin{array}{llll}
                \pi(\mathtt{R}) \,\defeq\, \mathtt{R},  &\pi(\mathtt{B}) \,\defeq \, \mathtt{B}, \quad &\pi(0) \,\defeq\, \mathtt{G},  &\pi(1) \,\defeq \, \mathtt{G},\\
            \tau(\mathtt{R}) \,\defeq\, \blank, \quad &\tau(\mathtt{B}) \,\defeq \, \blank, \quad &\tau(0) \,\defeq\, 0, \quad &\tau(1) \,\defeq \, 1.
            \end{array}
        \end{align*}
        Given $q \geq 4$, we let $\CW(q)$ be the LCL problem with labels $\set{\mathtt{R},\mathtt{B}, 0, 1}$ whose solutions 
        are precisely the functions $f$ 
        such that:
        \begin{enumerate}[label=\ep{{\normalfont\small{C}\arabic*}}]
            \item $\pi \circ f$ is a solution to $\Waffle(q)$, and
            \item $\tau \circ f$ is a partial 2-coloring. 
        \end{enumerate}
    \end{definition}

    In other words, to solve $\CW(q)$, we must first find a solution to $\Waffle(q)$ and then a partial 2-coloring whose domain is the set of green points. This is clearly an LCL problem with window $[0,2q]^n$. Our claim now is that taking $\Pi = \CW(q)$ for any $q \geq 4$ fulfils all the requirements of Theorem~\ref{theo:main}. 

    We start with the positive parts of Theorem~\ref{theo:main}. Parts \ref{item:meas_yes} and \ref{item:shift_yes} follow from the following lemma together with Propositions \ref{prop:MeasurableWaffle} and \ref{prop:shiftWaffle} respectively: 

    \begin{lemma}\label{lemma:BorelCW}
        If $q \geq 4$ and $\mathcal{R}$ is a Borel rectangular $q$-toast for some free Borel action $\Z^n \acts X$ on a standard Borel space $X$, then $X$ admits a $\CW(q)$-labeling $f$ such that $\pi \circ f = f_\mathcal{R}$ and whose restriction to the set $\bigcup {\mathcal{R}}$ is Borel. 
    \end{lemma}
    \begin{scproof}
        Note that $f_\mathcal{R}$ is a Borel $\Waffle(q)$-labeling of $X$. 
        Let $U \defeq f_\mathcal{R}\inv(\mathtt{G}) \cap \bigcup \mathcal{R}$. Observe that the graph $G(X)[U]$ is component-finite. Indeed, if $x \in U$, then $f_\mathcal{R}(x) = \mathtt{G}$ and there is a rectangle $K \in \mathcal{R}$ such that $x \in K$. This implies that the component of $G(X)[U]$ containing $x$ is included in $K \setminus \partial K$, which is a finite set, as desired. Hence, by Proposition~\ref{prop:component_finite_2color}, there exists a Borel partial $2$-coloring $c \colon X \to \set{0,1,\blank}$ with domain $U$. On the other hand, the graph $G(X)$ is bipartite, so there exists a \ep{possibly non-Borel} proper $2$-coloring $c' \colon X \to \set{0,1}$ of $G(X)$ \ep{i.e., a partial $2$-coloring with domain $X$}. It remains to set
        \[
            f(x) \,\defeq\, \begin{cases}
                f_\mathcal{R}(x) &\text{if } f_\mathcal{R}(x) \in \set{\mathtt{R}, \mathtt{B}}, \\
                c(x) &\text{if } x \in U,\\
                c'(x) &\text{otherwise}.
            \end{cases}\qedhere
        \]
    \end{scproof}

    Theorem~\ref{theo:main}\ref{item:meas_yes} now follows by invoking Proposition~\ref{prop:MeasurableWaffle} and Lemma \ref{lemma:BorelCW} to obtain a $\CW(q)$-labeling that is Borel on a conull set \ep{namely on $\bigcup \mathcal{R}$ for the square $q$-toast $\mathcal{R}$ provided by Proposition~\ref{prop:MeasurableWaffle}}, and hence measurable. Theorem~\ref{theo:main}\ref{item:shift_yes} follows in exactly the same way from Proposition~\ref{prop:shiftWaffle} and Lemma~\ref{lemma:BorelCW}. 
    To prove Theorem~\ref{theo:main}\ref{item:comp_yes}, we apply the following version of Lemma~\ref{lemma:BorelCW} in the computable setting:

    \begin{lemma}\label{lemma:compCW}
        If $q \geq 4$ and $\mathcal{R}$ is a computable complete rectangular $q$-toast for a free computable action $\Z^n \acts X$, 
        then $X$ admits a computable $\CW(q)$-labeling $f$ such that $\pi \circ f = f_\mathcal{R}$. 
    \end{lemma}
    \begin{scproof}
        The construction is the same as the proof of Lemma~\ref{lemma:BorelCW} but with Proposition~\ref{prop:comp_2color} used in place of Proposition~\ref{prop:component_finite_2color}. (Note that since $\mathcal{R}$ is complete, the ``otherwise'' case in the definition of $f$ does not occur.)
    \end{scproof}

    Theorem~\ref{theo:main}\ref{item:comp_yes} is an immediate consequence of Proposition~\ref{prop:compWaffle} and Lemma~\ref{lemma:compCW}.
    
    Now we turn our attention to the negative parts of Theorem~\ref{theo:main}, i.e., to items \ref{item:Baire_no} and \ref{item:ffiid_no}. Before we proceed, we require one more combinatorial fact about rectangular toasts, stated in the next lemma. It is related to 
    the Hex Theorem \cite{Hex}; 
    we give a direct proof here.

    \begin{lemma}\label{lem:hex}
        Let $\mathcal{R} \subseteq \finset{\Z^n}$ be a rectangular $q$-toast, where $q \geq 4$. In 
        the labeling $f_\mathcal{R}$, 
        any two green points in $\Z^n \setminus \bigcup \mathcal{R}$ are joined by a path in the graph $G_n$ all of whose vertices are green. 
    \end{lemma}
    \begin{scproof}
        %
        It is clear from the definition of $f_\mathcal{R}$ that $f_\mathcal{R}\inv(\mathtt{G}) \setminus \bigcup \mathcal{R}$ is exactly the set of points whose distance from $\bigcup \mathcal{R}$ is at least 2. 
        Let $x$ and $y$ be two such points and let $P = (x = x_0, x_1, \ldots, x_\ell = y)$ be an $xy$-path in $G_n$ minimizing $|V(P) \cap \partial^+\mathcal{R}|$. If $V(P) \cap \partial^+\mathcal{R} = \0$, then $P$ stays at distance at least $2$ from $\bigcup \mathcal{R}$, so every vertex of $P$ is green, as desired. Otherwise, let $K \in \mathcal{R}$ be a rectangle such that $V(P) \cap \partial^+ K \neq \0$ and let 
        $s$, $t$ be the minimum and maximum indices such that $x_s$, $x_{t} \in \partial^+ K$, respectively. Since $x$ has distance at least 2 from $K$, we must have $s > 0$ and, by the choice of $s$, $\dist(x_{s-1},K) = 2$. Similarly, $t < n$ and $\dist(x_{t+1}, K) = 2$. 
        We can then replace the segment of $P$ between $x_{s-1}$ and $x_{t+1}$ by an $x_{s-1} x_{t+1}$-path $Q$ that stays at distance $2$ or $3$ from $K$, resulting in a new $xy$-path $P'$. 
        Since $q \geq 4$, $V(Q) \cap \partial^+ \mathcal{R} = \0$, so we have $|V(P') \cap \partial^+\mathcal{R}| < |V(P) \cap \partial^+\mathcal{R}|$, contradicting the choice of $P$. 
    \end{scproof}

    With Lemma~\ref{lem:hex} in hand, we can now prove Theorem~\ref{theo:main}\ref{item:ffiid_no}:

    \begin{prop}\label{prop:no_fiid}
        If $q \geq 4$, then here is no finitary $\CW(q)$-labeling $f \colon \Free([0,1],\Z^n) \to \set{\mathtt{R}, \mathtt{B}, 0, 1}$. 
    \end{prop}
    \begin{scproof}
        Before proceeding, we note that Lemmas~\ref{lem:waffle_dichotomy} and \ref{lem:hex}, although stated for $q$-toasts on $\Z^n$, can be applied to any free transitive action of $\Z^n$, since such an action is isomorphic to the translation action $\Z^n \acts \Z^n$. To put this another way, given a free action $\Z^n \acts X$ and a rectangular $q$-toast $\mathcal{R} \subseteq \finset{X}$, we may apply  Lemmas~\ref{lem:waffle_dichotomy} and \ref{lem:hex} separately on every $\Z^n$-orbit $\mathcal{O} \subseteq X$.
        
        Suppose 
        $f$ is a finitary $\CW(q)$-labeling of $\Free([0,1],\Z^n)$. For $\mathsf{d} \in \set{\pm \mathsf{e}_1, \ldots, \pm\mathsf{e}_n}$, let $U_\mathsf{d}$ be the set of all points $x \in \Free([0,1], \Z^n)$ such that $f(i\mathsf{d} \cdot x) = \mathtt{G}$ for all $i \in \N$. Since $f$ is a measurable function, $U_\mathsf{d}$ is a measurable set. 

        \begin{claim}\label{claim:Ud_null}
            The set $U_\mathsf{d}$ is null for every $\mathsf{d} \in \set{\pm \mathsf{e}_1, \ldots, \pm\mathsf{e}_n}$.
        \end{claim}
        \begin{claimproof}
            Note that the map $c \colon \Free([0,1], \Z^n) \to \set{0,1,\blank}$ that agrees with $f$ on $U_\mathsf{d}$ and sends each $x \notin U_\mathsf{d}$ to $\blank$ is a partial $2$-coloring with domain $U_\mathsf{d}$. Thus, if $U_\mathsf{d}$ has positive measure, then by Proposition~\ref{prop:fiid_partial2color}, there exist $x \in U_\mathsf{d}$ and $k \in \N$ such that $x$ and $k\mathsf{d} \cdot x$ are $c$-separated. But that is impossible as $(x, \mathsf{d} \cdot x, \ldots, k\mathsf{d} \cdot x)$ is a path from $x$ to $k\mathsf{d} \cdot x$ contained in $U_\mathsf{d}$.
        \end{claimproof}

        Let $N \subseteq \Free([0,1],\Z^n)$ be a null Borel set such that $f$ is Borel on the complement of $N$ (which exists since $f$ is measurable). 
        Let $M \subseteq \Free([0,1],\Z^n)$ be a null Borel set containing $\bigcup_{\mathsf{d}} U_\mathsf{d}$, where the union is over all $\mathsf{d} \in \set{\pm \mathsf{e}_1, \ldots, \pm\mathsf{e}_n}$ (such a set $M$ which exists by Claim~\ref{claim:Ud_null}).
        Consider the following $\Z^n$-invariant Borel subset $X \subseteq \Free([0,1], \Z^n)$:
        \[ X \,\defeq\, \Free([0,1],\Z^n) \setminus (\Z^n \cdot (N \cup M)) \]
        Since the shift action is measure-preserving, $X$ is conull. Let
         \[
             \mathcal{R} \,\defeq\, \set{K \in \finset{X} \,:\, \text{$K$ is a rectangle such that $f(x) = \mathtt{R}$ for all $x \in \partial K$}}.
         \]
         The set $\mathcal{R}$ is clearly Borel. Now we apply Lemma~\ref{lem:waffle_dichotomy} to the $\Waffle(q)$-labeling $\pi \circ f$ on each orbit $\mathcal{O}$ of the shift action. Since, by construction, $X$ is disjoint from every $U_\mathsf{d}$, case \ref{item:waffle} must take place on all orbits $\mathcal{O} \subseteq X$. This means that $\mathcal{R}$ is a rectangular $q$-toast and $\pi \circ f$ agrees with $f_\mathcal{R}$ on $X$.
         
         Since $X$ is conull and $f$ is finitary, it follows that the characteristic function of $\partial \mathcal{R}$ is 
         weakly finitary. By Proposition~\ref{prop:noFFIID}, this implies that $\Leb^{\Z^n}(\bigcup \mathcal{R}) < 1$. Since $q \geq 4$, any point with distance $2$ from some $K \in \mathcal{R}$ is colored green by $f_\mathcal{R}$. Therefore, the set
         \[
            U \,\defeq\, \left\{x \in X \,:\, f_\mathcal{R}(x) = \mathtt{G} \text{ and } x \notin \bigcup\mathcal{R}\right\}
         \]
         has positive measure, and the map $c$ that sends each $x \in U$ to $f(x)$ and all other points to $\blank$ is a Borel partial $2$-coloring with domain $U$. By Proposition \ref{prop:fiid_partial2color}, there exist $x \in U$ and $k \in \N$ such that $x$ and $k\mathsf{e}_1 \cdot x$ are $c$-separated. However, by Lemma~\ref{lem:hex} applied on the $\Z^n$-orbit of $x$, there is a path from $x$ to $k\mathsf{e}_1 \cdot x$ contained in $U$, which is a contradiction. 
    \end{scproof}

    Theorem~\ref{theo:main}\ref{item:Baire_no} is proved in a similar fashion:

    \begin{prop}\label{prop:no_baire}
        Let $\Z^n \acts X$ be a free continuous action on a nonempty compact Polish space $X$ such that for each $\mathsf{d} \in \Sn$, the action of $2\mathsf{d}$ is minimal. Then $X$ does not admit a Baire-measurable $\CW(q)$-labeling. 
    \end{prop}
    \begin{scproof}
        The proof is exactly the same as the proof of Proposition~\ref{prop:no_fiid}, but with Propositions \ref{prop:baire_partial2color} and \ref{prop:noBaireWaffle} taking the place of Propositions \ref{prop:fiid_partial2color} and \ref{prop:noFFIID} respectively. \ep{To apply Proposition~\ref{prop:noBaireWaffle}, we observe that the minimality of the action of $2\mathsf{e}_1$ implies the same for $\mathsf{e}_1$.} 
    \end{scproof}

    The example following Proposition \ref{prop:noBaireWaffle} \ep{action on the circle via independent irrational rotations} again shows that an action as in Proposition~\ref{prop:no_baire} exists, proving item \ref{item:Baire_no} of Theorem \ref{theo:main}.

    \section{Adding extra structure: Proof of Theorem~\ref{theo:main_extra}}\label{sec:extra}

    In this section we prove Theorem~\ref{theo:main_extra}. For simplicity, we only present the proof in the case $n = 2$; the general argument is the same, \emph{mutatis mutandis}. In our construction we start with the problem $\CW(q)$ used to witness Theorem~\ref{theo:main} and then add some extra structure to it in order to ensure that the resulting problem can no longer be solved measurably on the free part of the shift, but still has a Baire-measurable solution on the free part of the shift as well as a computable solution on any free computable action. Similar (but much more complicated) arguments will be used in \cite{Sequel} to prove Theorem~\ref{theo:announce}.

        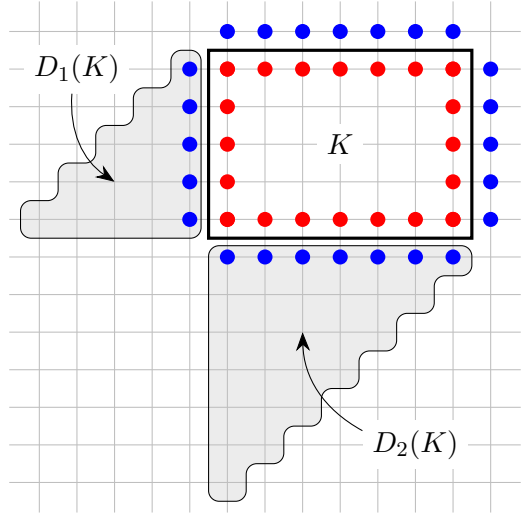
\begin{figure}[t]
			\centering
			\begin{tikzpicture}
                \draw[rounded corners,fill=gray!15,draw=white] (-0.5+0.15,-0.25) -- (-0.5+0.15,2+0.25) -- (-0.5-0.25,2+0.25) -- (-0.5-0.25,1.5+0.25) -- (-1-0.25,1.5+0.25) -- (-1-0.25,1+0.25) -- (-1.5-0.25,1+0.25) -- (-1.5-0.25,0.5+0.25) -- (-2-0.25,0.5+0.25) -- (-2-0.25,0.25) -- (-2.5-0.25,0.25)  -- (-2.5-0.25,-0.25) -- cycle;
                \draw[rounded corners,fill=gray!15,draw=white] (-0.25,-0.5+0.15) -- (3+0.25,-0.5+0.15) -- (3+0.25,-1+0.25) -- (2.5+0.25,-1+0.25) -- (2.5+0.25,-1.5+0.25) -- (2+0.25,-1.5+0.25) -- (2+0.25,-2+0.25) -- (1.5+0.25,-2+0.25) -- (1.5+0.25,-2.5+0.25) -- (1+0.25,-2.5+0.25) -- (1+0.25,-3+0.25) -- (0.5+0.25,-3+0.25) -- (0.5+0.25,-3.5+0.25) -- (0+0.25,-3.5+0.25) -- (0+0.25,-4+0.25) -- (-0.5+0.25,-4+0.25) -- cycle;
   
                \foreach \i in {-5,...,7}
                \draw[lightgray] (\i/2,-3.9) -- (\i/2,2.9);

                \foreach \j in {-7,...,5}
                \draw[lightgray] (-2.9,\j/2) -- (3.9,\j/2);

                \foreach \i in {0,...,6}
                {\fill[fill=red] (\i/2,0) circle (0.1cm);
                \fill[fill=red] (\i/2,2) circle (0.1cm);
                \fill[fill=blue] (\i/2,-0.5) circle (0.1cm);
                \fill[fill=blue] (\i/2,2.5) circle (0.1cm);}
                
                \foreach \j in {0,...,4}
                {\fill[fill=red] (0,\j/2) circle (0.1cm);
                \fill[fill=red] (3,\j/2) circle (0.1cm);
                \fill[fill=blue] (-0.5,\j/2) circle (0.1cm);
                \fill[fill=blue] (3.5,\j/2) circle (0.1cm);}
                
                \draw[very thick] (-0.25,-0.25) rectangle (3+0.25,2+0.25);
                \draw[rounded corners] (-0.5+0.15,-0.25) -- (-0.5+0.15,2+0.25) -- (-0.5-0.25,2+0.25) -- (-0.5-0.25,1.5+0.25) -- (-1-0.25,1.5+0.25) -- (-1-0.25,1+0.25) -- (-1.5-0.25,1+0.25) -- (-1.5-0.25,0.5+0.25) -- (-2-0.25,0.5+0.25) -- (-2-0.25,0.25) -- (-2.5-0.25,0.25)  -- (-2.5-0.25,-0.25) -- cycle;
                \draw[rounded corners] (-0.25,-0.5+0.15) -- (3+0.25,-0.5+0.15) -- (3+0.25,-1+0.25) -- (2.5+0.25,-1+0.25) -- (2.5+0.25,-1.5+0.25) -- (2+0.25,-1.5+0.25) -- (2+0.25,-2+0.25) -- (1.5+0.25,-2+0.25) -- (1.5+0.25,-2.5+0.25) -- (1+0.25,-2.5+0.25) -- (1+0.25,-3+0.25) -- (0.5+0.25,-3+0.25) -- (0.5+0.25,-3.5+0.25) -- (0+0.25,-3.5+0.25) -- (0+0.25,-4+0.25) -- (-0.5+0.25,-4+0.25) -- cycle;

                \node[fill=white] at (1.5,1) {$K$};
                \node[fill=white] (D1) at (-2,2) {$D_1(K)$};
                \node[fill=white] (D2) at (2.5,-3) {$D_2(K)$};

                \draw[-{Stealth[scale=1.5]}] (D1) to[bend right] (-1.5,0.5);
                \draw[-{Stealth[scale=1.5]}] (D2) to[bend left] (1,-1.5);
			\end{tikzpicture}
   \caption{An illustration for the proof of Theorem~\ref{theo:main_extra}. A rectangle $K \in \mathcal{R}$ and the two regions $D_1(K)$ and $D_2(K)$ in which the functions $h_1$ and $h_2$ respectively must equal $\POS$.}\label{fig:rectangle_nbhd}
	\end{figure}
        
        Fix $q \geq 4$ and let $\CWplus$ be the LCL problem 
        whose solutions on any free action $\Z^2 \acts X$ are precisely the triples $(f, h_1, h_2)$ with $f \colon X \to \set{\mathtt{R},\mathtt{B},0,1}$ and $h_1$, $h_2 \colon X \to \set{{\NEG},{\POS}}$ such that:
        \begin{enumerate}[label=\ep{{\normalfont\small{}${}^+$\arabic*}}]
            \item $f$ is a solution to $\CW(q)$,
            \item\label{item:must_be_green} if $f(x) = \mathtt{R}$, then $h_1(x) = h_2(x) = {\NEG}$,
            \item\label{item:start} if $f(x) = \mathtt{B}$ and $f(\mathsf{e}_i \cdot x) = \mathtt{R}$, then $h_i(x) = {\POS}$,
            \item\label{item:continue} if $h_i(\mathtt{e}_i \cdot x) = h_i((\mathtt{e}_1 + \mathtt{e}_2) \cdot x) = \POS$, then $h_i(x) = {\POS}$.
        \end{enumerate}
        Clearly, $\CWplus$ is an LCL problem on $\Z^2$ with label set $\set{\mathtt{R},\mathtt{B},0,1} \times \set{{\NEG},{\POS}}^2$ and window $[0,2q]^2$. We claim that this problem has the properties stated in Theorem~\ref{theo:main_extra}.

        We start by considering what a solution $(f,h_1, h_2)$ to $\CWplus$ on $\Z^2$ \ep{and hence also on an orbit of any free action of $\Z^2$} must look like. Suppose we are in case \ref{item:waffle} of Lemma~\ref{lem:waffle_dichotomy} applied to the $\Waffle(q)$-labeling $\pi \circ f$ and let $\mathcal{R} \subseteq \finset{\Z^2}$ be a rectangular $q$-toast such that $\pi \circ f = f_\mathcal{R}$. Take any rectangle $K \in \mathcal{R}$. For concreteness, say $K = [0,a] \times [0,b]$. Then for all $t \in [0,b]$, we have $f(0,t) = \mathtt{R}$ and $f(-1,t) = \mathtt{B}$, so, by \ref{item:start}, $h_1(-1,t) = \POS$. Applying \ref{item:continue} iteratively, we then see that $h_1(-1 - i, t) = \POS$ for all $0 \leq i \leq b$ and $0 \leq t \leq b - i$. Similarly, $h_2(s,-1-i) = \POS$ for all $0 \leq i \leq a$ and $0 \leq s \leq a - i$.  
        This situation is illustrated in Fig.~\ref{fig:rectangle_nbhd}. Since, by \ref{item:must_be_green}, $f$ cannot be equal to $\mathtt{R}$ when $h_1$ or $h_2$ is $\POS$, we conclude that the following must hold:

        \begin{claim}\label{claim:D}
            Suppose $(f,h_1,h_2)$ is a solution to $\CWplus$ on a free action $\Z^2\acts X$, where $\pi \circ f = f_\mathcal{R}$ for a rectangular $q$-toast $\mathcal{R}$. Let $K = ([0,a]\times[0,b]) \cdot x \in \mathcal{R}$ be a rectangle and let
            \begin{align*}
                D_1(K) \,&\defeq\, \set{(-1-i, t) \cdot x\,:\, 0 \leq i \leq b,\, 0 \leq t \leq b-i},\\
                D_2(K) \,&\defeq\, \set{(s, -1-i) \cdot x\,:\, 0 \leq i \leq a,\, 0 \leq s \leq a-i}.
            \end{align*}
            Then $(D_1(K) \cup D_2(K)) \cap \partial \mathcal{R} = \0$.
        \end{claim}

        It follows that in the setting of Claim~\ref{claim:D}, distinct rectangles $K$, $K' \in \mathcal{R}$ satisfy $D_1(K) \cap D_1(K') = \0$ and $D_2(K) \cap D_2(K') = \0$. Also, note that if $K$ is an $(a \times b)$-rectangle, then
        \begin{align}
            |D_1(K)| + |D_2(K)| \,&=\, \frac{(a+1)(a+2)}{2} + \frac{(b+1)(b+2)}{2} \nonumber\\
            &\geq\, \frac{(a+1)^2 + (b+1)^2}{2} \,\geq\, (a+1)(b+1) \,=\, |K|.\label{eq:large_sum}
        \end{align}
        
        We are now ready to show $\CWplus \notin \FIID(\Z^2)$:

        \begin{claim}
            There is no measurable solution to $\CWplus$ on $\Free([0,1],\Z^2)$ \ep{i.e., $\CWplus \notin \FIID(\Z^2)$}.
        \end{claim}
        \begin{scproof}
            Suppose for contradiction that $(f,h_1,h_2)$ is a measurable $\CWplus$-labeling of $\Free([0,1],\Z^2)$. As in the proof of Proposition~\ref{prop:no_fiid}, we see that there is a Borel rectangular $q$-toast $\mathcal{R}$ such that $\pi \circ f$ agrees with $f_\mathcal{R}$ almost everywhere. For a rectangle $K \in \mathcal{R}$ with side lengths $a$ and $b$, let its \emphd{root} be the point $x$ such that $K = ([0,a]\times[0,b]) \cdot x$, and let $R$ be the set of all roots of the rectangles in $\mathcal{R}$. Also, for each $x \in R$, let $K_x$ be the unique rectangle in $\mathcal{R}$ such that $x$ is its root. As discussed above, Claim~\ref{claim:D} implies that for each $i \in \set{1,2}$, the sets $\set{D_i(K) \,:\, K \in \mathcal{R}}$ are pairwise disjoint. Since the action $\Z^2 \acts \Free([0,1],\Z^2)$ is measure-preserving, it follows that
            \[
                \int_{x \in R} |D_i(K_x)| \,\mathrm{d}\Leb^{\Z^2}(x) \,=\, \Leb^{\Z^2}\left(\bigcup_{K \in \mathcal{R}} D_i(K) \right) \,\leq\, 1.
            \]
            Therefore, by \eqref{eq:large_sum},
            \[
                \int_{x \in R} |K_x| \,\mathrm{d}\Leb^{\Z^2}(x) \,\leq\, \int_{x \in R} \big(|D_1(K_x)| + |D_2(K_x)|\big) \,\mathrm{d}\Leb^{\Z^2}(x) \,\leq\, 2.
            \]
            Using the invariance of the measure again, we see that the expected number of rectangles in $\mathcal{R}$ containing a random point in $\Free([0,1],\Z^2)$ is at most $2$. In particular, 
            almost every point 
            is contained in finitely many rectangles. 
            By Lemma~\ref{lemma:Wcomplete}, this means that the set of all $\mathcal{R}$-complete points is null, i.e., $\Leb^{\Z^2}(\Z^2 \cdot (X \setminus \bigcup \mathcal{R})) = 1$, which implies that $\Leb^{\Z^2}(\bigcup \mathcal{R}) < 1$. Now we use Proposition \ref{prop:fiid_partial2color} and Lemma~\ref{lem:hex} to arrive at a contradiction in exactly the same way as in the last step of the proof of Proposition~\ref{prop:no_fiid}.
        \end{scproof}

        Now we verify the other two properties of $\CWplus$ stated in the theorem.
        
        \begin{claim}\label{claim:Baireshift_big_gaps}
            There is a Baire-measurable solution to $\CWplus$ on $\Free(2,\Z^2)$ \ep{i.e., $\CWplus \in \BS(\Z^2)$}.
        \end{claim}
        \begin{scproof}
            We 
            modify the construction in the proof of Proposition~\ref{prop:shiftWaffle}. Define $\rho \colon \Free(2, \Z^2) \to 2$ by $\rho(x) \defeq x(\bm{0})$. Call a square $K$ in $\Free(2, \Z^2)$ \emphd{safe} if
        \begin{itemize}
            \item the side length of $K$, say $N$, is at least $q$,
        
            \item for all $x \in X$:
            \begin{itemize}
                \item if $x \in \partial K$, then $\rho(x) = 1$,
                \item if $0 < \dist(x, \partial K) \leq q$, then $\rho(x) = 0$,
                \item and, additionally, if $x \notin K$ and $\dist(x, \partial K) \leq 2N$, then $\rho(x) = 0$.
            \end{itemize}
        \end{itemize}
        Let $\mathcal{R}$ be the set of all safe squares. The same argument as in the proof of Proposition~\ref{prop:shiftWaffle} shows that $\mathcal{R}$ is a square $q$-toast and $\bigcup \mathcal{R}$ is a comeager \ep{in fact, dense open} set. \ep{Note that the square $R$ used there to witness the density of $\bigcup \mathcal{R}$ remains safe with the modified definition.} Then Lemma~\ref{lemma:BorelCW} yields a Baire-measurable $\CW(q)$-labeling $f$ of $\Free([0,1],\Z^2)$ such that $\pi \circ f = f_\mathcal{R}$. Note that if $K \in \mathcal{R}$ is an $N$-square, then for all $K' \in \mathcal{R}$ with $K' \not\subseteq K$, we have $\dist(\partial K, \partial K') > 2N$. This implies that for each $i \in \set{0,1}$, the sets $D_1(K)$ and $D_2(K)$ are disjoint from $\partial \mathcal{R}$. Therefore, letting
        \[
            h_i(x) \,\defeq\, \begin{cases}
                \POS &\text{if } x \in D_i(K) \text{ for some } K \in \mathcal{R},\\
                \NEG &\text{otherwise},
            \end{cases}
        \]
        we see that $(f, h_1, h_2)$ is a Baire-measurable solution to $\CWplus$ on $\Free([0,1], \Z^2)$, as desired.
        \end{scproof}

        \begin{claim}
            For every free computable action $\Z^2 \acts X$ on a computable set $X \subseteq \N$, there exists a computable solution to $\CWplus$ on $X$. In other words, $\CWplus \in \COMP(\Z^2)$.
        \end{claim}
        \begin{scproof}
            Without loss of generality, assume that $X = \N$. As in the proof of Proposition~\ref{prop:compWaffle}, we describe an algorithm that enumerates a sequence $R_0$, $R_1$, \ldots{} with the following properties:
        \begin{itemize}
            \item each $R_i$ is an $N_i$-square in $X$ for some $N_i \geq i + q$,
            \item for all $i \neq j$, we have $\dist(\partial R_i, \partial R_j) > q$,
            \item for all $x \in \N$, $x \in R_x$,
            \item and, additionally, if $R_j \not\subseteq R_i$, then $\dist(\partial R_i, \partial R_j) > 2N_i$.
        \end{itemize}
        Namely, once $R_0$, \ldots, $R_{i-1}$ have been computed, we cycle through the integers $j \geq i + q$ until we find one such for all $t \in \set{0,\ldots, i-1}$,
        \[
            \dist(\partial ([-j, j]^2 \cdot i),\, \partial R_t) \,>\, \max\set{q, 2N_t},
        \]
        and, moreover, if $R_t \not\subseteq [-j,j]^2$, then
        \[
            \dist(\partial ([-j, j]^2 \cdot i),\, \partial R_t) \,>\, 4j.
        \]
        Such $j$ must exist since the set $\set{R_0, \ldots, R_{i-1}}$ is finite. We then let $N_i \defeq 2j$ and $R_i \defeq [-j, j]^n \cdot i$.
        
        As in the proof of Proposition~\ref{prop:compWaffle}, $\mathcal{R} \defeq \set{R_i \,:\, i \in \N}$ is a computable compete square $q$-toast. Hence, by Lemma~\ref{lemma:compCW}, there is a computable $\CW(q)$-labeling $f$ with $\pi \circ f = f_\mathcal{R}$. As in the proof of Claim~\ref{claim:Baireshift_big_gaps}, we now set, for each $i \in \set{1,2}$,
        \[
            h_i(x) \,\defeq\, \begin{cases}
                \POS &\text{if } x \in D_i(R_t) \text{ for some } t \in \N,\\
                \NEG &\text{otherwise}.
            \end{cases}
        \]
        Notice that the function $h_i$ is computable, since if $x \in D_i(R_t)$ for some $t \in \N$, then $R_t \subseteq R_x$, so $t \leq N_x$. Thus, $(f,h_1,h_2)$ is a computable solution to $\CWplus$, as desired.
        \end{scproof}

        We have now established all the required properties of $\CWplus$, and thus the proof of Theorem~\ref{theo:main_extra} is complete.

    \section{Proof of Theorem~\ref{theo:FFIID_Baire}}\label{sec:FFIID_Baire}

    The goal of this section is to prove Theorem~\ref{theo:FFIID_Baire}, i.e., that the inclusion $\FFIID(\G) \subseteq \BS(\G)$ is valid for every group $\G$. As mentioned in the introduction, this result follows from the general analysis of the Baire-measurable combinatorics of the shift action presented in the paper \cite{BerGeneric} by the second named author. The theory developed in \cite{BerGeneric} is more general that what we require; below we only describe the tools necessary to deduce Theorem~\ref{theo:FFIID_Baire}.

    Fix a group $\G$ and a finite set $\Lambda$. Let $\Shift_{<\infty}(\Lambda, \G)$ be the set of all partial maps $\phi \colon \G \pto \Lambda$ with finite domains. There is a natural shift action $\G \acts \Shift_{<\infty}(\Lambda, \G)$ given by
    \[
        (\gamma \cdot \phi)(\delta) \,\defeq\, \phi(\delta\gamma) \quad \text{for all $\gamma$, $\delta \in \G$ and $x \in \Shift_{<\infty}(\Lambda,\G)$},
    \]
    where if $\delta\gamma \notin \dom(\phi)$, then $(\gamma \cdot \phi)(\delta)$ is left undefined. We also let $\gamma \cdot D \defeq D\gamma^{-1}$ for every finite subset $D \subseteq \G$, so that $\dom(\gamma \cdot \phi) = \gamma \cdot \dom(\phi)$. A \emphd{$(\Lambda,\Gamma)$-ideal} is a set $\mathcal{I} \subseteq \Shift_{<\infty}(\Lambda, \Gamma)$ that is invariant under the shift action of $\Gamma$ and closed under taking restrictions of functions. A $(\Lambda,\Gamma)$-ideal $\mathcal{I}$ is \emphd{extendable} if for all $\phi \in \mathcal{I}$ and $\gamma \in \G \setminus \dom(\phi)$, there is some $\psi \in \mathcal{I}$ such that $\phi \subset \psi$ and $\gamma \in \dom(\psi)$. 

    \begin{definition}[Join property]
        A $(\Lambda,\Gamma)$-ideal $\mathcal{I}$ has the \emphd{join property} if it is possible to assign to each $\phi \in \mathcal{I}$ a finite subset $D(\phi) \in \Gamma$ so that, for all $k \in \N$, $\phi_1$, \ldots, $\phi_k \in \mathcal{I}$, and $\gamma_1$, \ldots, $\gamma_k \in \G$, if the sets $\gamma_1 \cdot D(\phi_1)$, \ldots, $\gamma_k \cdot D(\phi_k)$ are pairwise disjoint, then $(\gamma_1 \cdot \phi_1) \cup \ldots \cup (\gamma_k \cdot \phi_k) \in \mathcal{I}$. 
%
    \end{definition}

    We can now state the (somewhat weakened form of the) result from \cite{BerGeneric} that we shall use:

    \begin{theo}[{AB \cite[special case of Thm.~2.10]{BerGeneric}}]\label{theo:Generic}
        Let $\G$ be a finitely generated group and let $\Pi = (W, \Lambda, \mathcal{A})$ be an LCL problem on $\G$. Then $\Pi \in \BS(\G)$ if and only if there exists an extendable $(\Lambda, \G)$-ideal $\mathcal{I} \subseteq \Shift_{<\infty}(\Lambda, \G)$ such that:
        \begin{itemize}
            \item each $\phi \in \mathcal{I}$ can be extended to a $\Pi$-labeling of $\G$, and
            \item $\mathcal{I}$ has the join property.
        \end{itemize}
    \end{theo}

    With this, we are ready to prove Theorem~\ref{theo:FFIID_Baire}.

    \begin{scproof}[ of Theorem~\ref{theo:FFIID_Baire}]
        The argument is similar to \cite[proof of Lem.~5.1]{BerGeneric}, with measure taking the place of Baire category. Let $\G$ be a finitely generated group and let $\Pi = (W, \Lambda, \mathcal{A})$ be an LCL problem on $\G$. Suppose $\Pi \in \FFIID(\G)$ and let $f \colon \Free([0,1],\G) \to \Lambda$ be a finitary $\Pi$-labeling of $\Free([0,1],\G)$. We shall use $f$ to define a $(\Lambda,\G)$-ideal $\mathcal{I} \subseteq \Shift_{<\infty}(\Lambda,\Gamma)$ satisfying the conditions of Theorem~\ref{theo:Generic}, from which it will follow that $\Pi \in \BS(\G)$, as desired.
        
        Define a measurable $\G$-equivariant function $\pi \colon \Free([0,1],\G) \to \Shift(\Lambda,\G)$ by
        \[
            \big(\pi(x)\big)(\gamma) \,\defeq\, f(\gamma \cdot x) \quad \text{for all $x \in \Free([0,1],\G)$ and $\gamma \in \G$}.
        \]
        Since $f$ is a solution to $\Pi$, $\pi(x)$ is a $\Pi$-labeling of $\G$ for every $x \in \Free([0,1],\G)$. For $\phi \in \Shift_{< \infty} (\Lambda, \G)$, let $U_\phi \subseteq \Shift(\Lambda, \G)$ be the clopen set of all functions $\G \to \Lambda$ that extend $\phi$ and define
        \[
            \mathcal{I} \,\defeq\, \big\{\phi \in \Shift_{<\infty}(\Lambda, \G) \,:\, \pi^{-1}(U_\phi) \text{ has positive measure}\big\}.
        \]
        It is clear that, since the shift action $\G \acts \Free([0,1],\G)$ is measure-preserving, $\mathcal{I}$ is a $(\Lambda,\G)$-ideal. Moreover, since $\pi$ sends points in $\Free([0,1],\G)$ to $\Pi$-labelings of $\G$, each member of $\mathcal{I}$ can be extended to a $\Pi$-labeling of $\G$. It remains to verify that $\mathcal{I}$ is extendable and has the join property.

        Take any $\phi \in \mathcal{I}$ and $\gamma \in \G \setminus \dom(\phi)$. For each $\lambda \in \Lambda$, let $\phi_\lambda \colon \dom(\phi) \cup \set{\gamma} \to \Lambda$ be the mapping that agrees with $\phi$ on $\dom(\phi)$ and sends $\gamma$ to $\lambda$ . Note that
        \[
            U_\phi \,=\, \bigcup_{\lambda \in \Lambda} U_{\phi_\lambda}, \qquad \text{hence} \qquad  \pi^{-1}(U_\phi) \,=\, \bigcup_{\lambda \in \Lambda} \pi^{-1}(U_{\phi_\lambda}).
         \]
         Since $\phi \in \mathcal{I}$, the set $\pi^{-1}(U_\phi)$ has positive measure. It follows that for some $\lambda \in \Lambda$, $\pi^{-1}(U_{\phi_\lambda})$ must also have positive measure, i.e., $\phi_\lambda \in \mathcal{I}$. This shows that $\mathcal{I}$ is extendable.

         Finally, we need to argue that $\mathcal{I}$ has the join property. This is where we use that $f$ is finitary. Consider any $\phi \in \mathcal{I}$. 
         For each finite set $D \subseteq \G$, let
         \[
            X_{\phi, D} \,\defeq\, \big\{x \in \Free([0,1],\G) \,:\, \text{$\pi(x') \in U_\phi$ for almost all $x' \in \Free([0,1],\G)$ that agree with $x$ on $D$}\big\}.
         \]
         Let us record the following observations, which follow immediately from the above definition:
         \begin{itemize}
            \item The set $X_{\phi,D}$ is measurable by Fubini's theorem.
            
             \item The membership of a point $x$ in $X_{\phi,D}$ is determined by the restriction of $x$ to $D$.

             \item For almost all $x \in X_{\phi, D}$, we have $\pi(x) \in U_\phi$.
         \end{itemize}
         Since $f$ is finitary, almost every point $x \in \pi^{-1}(U_\phi)$ belongs to $X_{\phi,D}$ for some finite set $D \subseteq \G$. Thus, we can pick a finite set $D(\phi) \subseteq \G$ such that the measure of the set $X(\phi) \defeq X_{\phi, D(\phi)}$ is positive.
         
        Now take $k \in \N$, $\phi_1$, \ldots, $\phi_k \in \mathcal{I}$, $\gamma_1$, \ldots, $\gamma_k \in \G$ and suppose the sets $\gamma_1 \cdot D(\phi_1)$, \ldots, $\gamma_k \cdot D(\phi_k)$ are pairwise disjoint. Sample a random point $x \in \Free([0,1],\G)$. Then, for each $1 \leq i \leq k$, $x \in \gamma_i \cdot X(\phi_i)$ with positive probability. The random event $E_i \defeq \set{x \in \gamma_i \cdot X(\phi_i)}$ is determined by the restriction of $x$ to $\gamma_i \cdot D(\phi_i)$. Since the sets $\gamma_1 \cdot D(\phi_1)$, \ldots, $\gamma_k \cdot D(\phi_k)$ are disjoint, the events $E_1$, \ldots, $E_k$ are mutually independent, so with positive probability all of them occur simultaneously. In other words, the following set has positive measure:
        \[
            E \,\defeq\, \set{x \in \Free([0,1],\Gamma) \,:\, x \in \gamma_i \cdot X(\phi_i) \text{ for all } 1 \leq i \leq k}.
        \]
        Note that for almost all $x \in E$ and for all $1 \leq i \leq k$, we have $\pi(x) \in U_{\gamma_i \cdot \phi_i}$, or, equivalently,
        $
            \pi(x) \in U_{(\gamma_1 \cdot \phi_1) \cup \ldots \cup (\gamma_k \cdot \phi_k)}
        $.
        Therefore, $(\gamma_1 \cdot \phi_1) \cup \ldots \cup (\gamma_k \cdot \phi_k) \in \mathcal{I}$, and the proof is complete.
    \end{scproof}

    \section{Final remarks and open problems}\label{sec:final}

    \subsection{Unlabeled grids}\label{subsec:unlabeled}

    Throughout this paper, we have been focusing on LCL problems that are defined on finitely generated groups \ep{in particular, $\Z^n$} and their actions. In computer science and probability theory, it is also common to consider problems on ``pure'' graphs, without any additional structure. In the context of this paper, it is natural to ask whether our results extend to the framework of LCL problems on \emphd{unlabeled, undirected $n$-dimensional grids}, i.e., graphs isomorphic to the standard Cayley graph $G_n$ of $\Z^n$ \ep{without edge directions or markers indicating which generator each edge corresponds to}. The short answer is ``yes,'' with virtually the same proofs.
    
    Formally, we can define an \emphd{LCL problem on $G_n$} as an LCL problem $\Pi$ on $\Z^n$ such that the set of all $\Pi$-labelings of $\Z^n$ is invariant under the natural action of the group $\mathsf{Aut}(G_n)$ of graph-theoretic automorphisms of $G_n$. Note that $\mathsf{Aut}(G_n)$ is generated by the group of translations of $\Z^n$ \ep{which is isomorphic to $\Z^n$} and the group $\mathsf{Aut}_\bullet(G_n)$ of \emphd{root-preserving automorphisms}, i.e., the automorphisms $\phi \colon \Z^n \to \Z^n$ satisfying $\phi(\bm{0}) = \bm{0}$.  The group $\mathsf{Aut}_\bullet(G_n)$ is finite and generated by permutations and sign changes of the generators.
    
    The role of $\Z^n$-actions in this setting is taken by \emphd{$G_n$-graphs}, i.e., graphs in which every connected component is isomorphic to $G_n$. If $\Pi$ is an LCL problem on $\Z^n$ and $G$ is a graph isomorphic to $G_n$, we call $f \colon V(G) \to \Lambda$ a \emphd{solution} to $\Pi$ on $G$, or a \emphd{$\Pi$-labeling} of $G$, if $f \circ \phi$ is a solution to $\Pi$ on $\Z^n$ for any \ep{hence every} isomorphism $\phi \colon \Z^n \to V(G)$. If $G$ is a $G_n$-graph, then $f \colon V(G) \to \Lambda$ is a \emphd{solution} to $\Pi$ on $G$ if it is a $\Pi$-labeling of each component of $G$. Note that if $\Z^n \acts X$ is a free action, then the corresponding graph $G(X)$ \ep{as defined in \S\ref{subsec:qwaffleproblem} on p.~\pageref{page:G(X)}} is a $G_n$-graph. Conversely, every $G_n$-graph is generated by a free $\Z^n$-action in this way. However, not all \emph{Borel} $G_n$-graphs are generated by free \emph{Borel} $\Z^n$-actions \ep{even if $n = 1$ \cites[Rmk.~6.8]{KechrisMiller}{MillerThesis}; see also \cite[\S5]{RileyEffective} for an analysis of the $n = 2$ case}. In other words, in this framework we consider a restricted class of LCL problems, but attempt to solve them on a broader class of structures.
    
    The classes of LCL problems $\BOREL(G_n)$, $\MEAS(G_n)$, and so on are defined by modifying the definitions of $\BOREL(\Z^n)$, $\MEAS(\Z^n)$, etc.~in obvious ways. 
    We now have the following analog of Theorem~\ref{theo:main} in this setting:
    \vspace{-4pt}
    \begin{tcolorbox}
    \vspace{-5pt}
	    \begin{theorem}\label{theo:main_unlabeled}
	        For every $n \geq 2$, there exists an LCL problem $\Pi$ on $G_n$ such that:
            \begin{enumerate}[label=\ep{\normalfont{\roman*}}]
                \item $\Pi \in \MEAS(G_n)$,

                \item\label{item:unlabeled_Baire_no} $\Pi \notin \BAIRE(\Z^n)$ \ep{hence it is also not in $\BAIRE(G_n)$},

                \item\label{item:unlabeled_Baireshift_yes} $\Pi$ has a Baire-measurable solution on $G(\Free(2, \Z^n))$, 

                \item\label{item:unlabeled_finitary_no} $\Pi$ has no finitary solution on $G(\Free([0,1],\Z^n))$,

                \item $\Pi \in \COMP(G_n)$.
            \end{enumerate}
	    \end{theorem}
        \vspace{-10pt}
	\end{tcolorbox}
    To prove Theorem~\ref{theo:main_unlabeled}, we simply observe that the problem $\CW(q)$ used to witness Theorem~\ref{theo:main} is, in fact, an LCL problem on $G_n$, because the notions of a rectangular $q$-toast and a partial $2$-coloring are automorphism-invariant. 
    This immediately shows that parts \ref{item:unlabeled_Baire_no}, \ref{item:unlabeled_Baireshift_yes}, and \ref{item:unlabeled_finitary_no} of Theorem~\ref{theo:main_unlabeled} follow from the corresponding items in Theorem~\ref{theo:main}. Modifying the arguments that go into proving parts \ref{item:meas_yes} and \ref{item:comp_yes} of Theorem~\ref{theo:main} to make sure they apply to $G_n$-graphs and not just $\Z^n$-actions is routine, and we omit the details.
    
    
    We conclude this discussion with two remarks. First, by taking the quotient under the natural action $\mathsf{Aut}_\bullet(G_n) \acts \Free(2, \Z^n)$, it is 
    not difficult to arrange the Baire-measurable labeling in item \ref{item:unlabeled_Baireshift_yes}
    of Theorem~\ref{theo:main_unlabeled}
    to be $\mathsf{Aut}_\bullet(G_n)$-invariant. 
    Second, although the problem $\CWplus$ used in \S\ref{sec:extra} to prove Theorem~\ref{theo:main_extra} is not an LCL problem on $G_n$, it is not hard
    to modify it to obtain an analogous version of Theorem~\ref{theo:main_extra} for unlabeled grids. 
    
    \subsection{Open problems}

    Many fundamental questions about the complexity of LCL problems on $\Z^n$ remain open. Fig.~\ref{fig:relations} summarizes the current state of our knowledge of this matter, and any relation not shown there leads to an interesting open problem. Here we give a list of a few such problems that we find particularly intriguing. In what follows, $n$ is an integer at least $2$.

    \begin{question}
        Does the inclusion $\BAIRE(\Z^n) \subset \MEAS(\Z^n)$ hold? \ep{If it does, it must be strict.}
    \end{question}

    \begin{question}
        Do we have $\BAIRE(\Z^n) = \BOREL(\Z^n)$?
    \end{question}

    \begin{question}\label{ques:Meas_vs_FIID}
        Do we have $\MEAS(\Z^n) = \FIID(\Z^n)$?
    \end{question}

    \begin{question}
        What is the relationship between $\FFIID(\Z^n)$ and $\BOREL(\Z^n)$? 
        What about $\FFIID(\Z^n)$ and $\BAIRE(\Z^n)$? 
        No inclusion between these classes has been ruled out yet.
    \end{question}

    \begin{question}
        Does the \ep{necessarily strict} inclusion $\BAIRE(\Z^n) \subset \COMP(\Z^n)$ hold? 
    \end{question}


    \begin{question}
        Call a free Borel action $\Z^n \acts X$ on a Polish space $X$ \emphd{Baire-universal} if $\BAIRE(\Z^n)$ is precisely the class of all LCL problems that admit a Baire-measurable solution on $X$. Note that Baire-universal actions exist. Namely, for each LCL problem $\Pi \notin \BAIRE(\Z^n)$, we can fix a free Borel action $\Z^n \acts X_\Pi$ that has no Baire-measurable $\Pi$-labeling and then take the disjoint union of $X_\Pi$ over all $\Pi \notin \BAIRE(\Z^n)$. However, this is not a very satisfying example, since we do not know in general which problems are not in $\BAIRE(\Z^n)$. Is there a more ``natural,'' ``concrete'' example of a Baire-universal action? 
        The shift action $\Z^n \acts \Free(2,\Z^n)$ is an obvious candidate, but it is \emph{not} Baire-universal by Theorem~\ref{theo:not_universal}. But what about, say, the action on the circle via independent irrational rotations (see the paragraph following Proposition \ref{prop:noBaireWaffle})? Is it Baire-universal? More abstractly, does there exist a Baire-universal free continuous action $\Z^n \acts X$ on a Polish space $X$ with a dense orbit (such actions are called \emphd{generically ergodic} \cite[Prop.~20.18]{AnushDST})?
    \end{question}


    There are also plenty of interesting problems for other finitely generated groups. Even for free groups, which are by far the best understood, some questions remain unanswered. For example:
    
    \begin{question}\label{ques:free}
        Do we have $\FFIID(\mathbb{F}_n) = \FIID(\mathbb{F}_n)$ or $\FIID(\mathbb{F}_n) = \MEAS(\mathbb{F}_n)$? We conjecture the answer is ``no'' to both. 
    \end{question}

    Furthermore, there remain general questions that are of interest for all finitely generated groups. Here are two examples, the first of which is closely related to Question~\ref{ques:Meas_vs_FIID}.

    \begin{question}
        Given a finitely generated group $\G$ and a standard probability space $(X,\mu)$, let $\FIID(X, \mu, \G)$ be the class of all LCL problems on $\G$ that have a measurable solution on $\Free(X, \G)$. As a special case, $\FIID(\G) = \FIID([0,1], \Leb, \G)$, where $\Leb$ is the Lebesgue measure on $[0,1]$. Assuming $\mu$ is not concentrated on a single point, does the class $\FIID(X, \mu, \G)$ depend on the space $(X,\mu)$? It does not when $\G = \Z$ \ep{because $\FIID(\Z) = \MEAS(\Z)$} and when $\G$ is non-amenable by a result of Bowen \cite[Corl.~1.3]{LBowen}. But what about, say, $\G = \Z^n$?
    \end{question}

    \begin{question}\label{ques:pmp}
        Let $\MEAS_{\mathsf{pmp}}(\G)$ be the class of all LCL problems on $\G$ that admit a measurable solution on every free \emph{measure-preserving} action $\G \acts (X,\mu)$, where $(X,\mu)$ is a standard probability space. Clearly, $\MEAS(\G) \subseteq \MEAS_{\mathsf{pmp}}(\G)$. Can this inclusion be strict? Once again, we know $\MEAS_{\mathsf{pmp}}(\Z) = \MEAS(\Z)$ because $\FIID(\Z) = \MEAS(\Z)$, but we do not know the answer for any group that is not virtually $\Z$. 
    \end{question}

        
    \subsection{Acknowledgments}

    The fourth named author is grateful to Riley Thornton for suggesting Proposition~\ref{prop:compWaffle} as a possible tool for tackling Question~\ref{q:comp_vs_baire}, and to Andrew Marks for telling us about the source \cite{Oberwolfach}. We are also grateful to Jan Greb\'ik for insightful conversations, in particular for directing our attention to Question~\ref{ques:ffiid}, and to Alexander Kechris and Anush Tserunyan for helpful comments on an earlier version of the manuscript. Finally, we thank the organizers of the \emph{Workshop on Measurable Combinatorics} that was held at the Alfr\'ed R\'enyi Institute of Mathematics in Budapest, Hungary in June 2024, where part of this research took place.

    This material is based upon work partially supported by the National Science Foundation under grants DGE-2146752, DMS-2528522, and DMS-2402064. Any opinions, findings, and conclusions or recommendations expressed in this material are those of the authors and do not necessarily reflect the views of the National Science Foundation.

\printbibliography

@article{NaorStock,
	author = {M. Naor and L. Stockmeyer},
	title = {What can be computed locally?},
	journaltitle = {SIAM J.	Comput.},
	volume = {24},
	number = {6},
	pages = {1259--1277},
	date = {1995},
}

@book{subshifts1,
    author = {D. Lind and B. Marcus},
    title = {An Introduction to Symbolic Dynamics and Coding},
    publisher = {Cambridge University Press},
    date = {1995},
}

@book{subshifts2,
    author = {T. Ceccherini-Silberstein and M. Coornaert},
    title = {Cellular Automata and Groups},
    location = {Berlin, Heidelberg},
    publisher = {Springer},
    date = {2010},
}

@inbook{ABJ_survey,
	author = {N. Aubrun and S. Barbieri and E. Jeandel},
	title = {Chapter 9. About the domino problem for subshifts on groups},
	booktitle = {Sequences, Groups, and Number Theory},
	date = {2018},
	editor = {V. Berth{\'{e}} and M. Rigo},
	series = {Trends in Mathematics},
	pages = {331--389},
}

@article{Wang,
    author = {H. Wang},
    title = {Proving theorems by pattern recognition --- II},
    journaltitle = {Bell Syst. Tech. J.},
    volume = {40},
    number = {1},
    pages = {1--41},
    date = {1961},
}

@article{Berger,
	author = {R. Berger},
	title = {The Undecidability of the Domino Problem},
	journaltitle = {Mem. Am. Math. Soc.},
	volume = {66},
	date = {1966},
}

@article{Robinson,
    author = {R.M. Robinson},
    title = {Undecidability and nonperiodicity for tilings of the plane},
    journaltitle = {Invent. Math.},
    volume = {12},
    pages = {177--209},
    date = {1971},
}

@book{BE,
	author = {L. Barenboim and M. Elkin},
	title = {Distributed Graph Coloring: Fundamentals and Recent Developments},
	series = {Synthesis Lectures on Distributed Computing Theory},
	date = {2013},
}

@article{CS1,
	author = {S. Brandt and O. Fischer and J. Hirvonen and B. Keller and T. Lempi{\"{a}}inen and J. Rybicki and J. Suomela and J. Uitto},
	title = {A lower bound for the distributed Lov{\'{a}}sz local lemma},
	journaltitle = {ACM Symposium on Theory of Computing (STOC)},
	date = {2016},
	pages = {479--488},
	addendum = {Full version: \url{https://arxiv.org/abs/1511.00900}},
}

@article{CS2,
	title = {An exponential separation between randomized and deterministic complexity in the $\mathsf{LOCAL}$ model},
	author = {Y.-J. Chang and T. Kopelowitz and S. Pettie},
	journaltitle = {SIAM J. Comput.},
	volume = {48},
	number = {1},
	pages = {122--143},
	date = {2019},
}

@article{CS3,
	title = {A time hierarchy theorem for the $\mathsf{LOCAL}$ model},
	author = {Y.-J. Chang and S. Pettie},
	journaltitle = {SIAM J. Comput.},
	volume = {48},
	number = {1},
	pages = {33--69},
	date = {2019},
}

@article{CS4,
	author = {Y.-J. Chang and Q. He and W. Li and S. Pettie and J. Uitto},
	title = {The complexity of distributed edge coloring with small palettes},
	journaltitle = {ACM-SIAM Symposium on Discrete Algorithms \ep{SODA}},
	date = {2018},
	pages = {2633--2652},
	addendum = {Full version: \url{https://arxiv.org/abs/1708.04290}},
}

@article{CS5,
	title = {Almost global problems in the $\mathsf{LOCAL}$ model},
	author = {A. Balliu and S. Brandt and D. Olivetti and J. Suomela},
	journaltitle = {Distrib. Comput.},
	volume = {34},
	pages = {259--281},
	date = {2021},
}

@article{CS6,
	author = {Y.-J. Chang},
	title = {The complexity landscape of distributed locally checkable problems on trees},
	journaltitle = {International Symposium on Distributed Computing \ep{DISC}},
	date = {2020},
	pages = {18:1--18:17},
	addendum = {Full version: \url{https://arxiv.org/abs/2009.09645}},
}

@article{CS7,
	author = {C. Grunau and V. Rozho\v{n} and S. Brandt},
	title = {The landscape of distributed complexities on trees and beyond},
	journaltitle = {ACM Symposium on Principles of Distributed Computi \ep{PODC}},
	date = {2022},
	pages = {37--47},
	addendum = {Full version: \url{https://arxiv.org/abs/2202.04724}},
}

@article{CS8,
	author = {S. Brandt and J. Hirvonen and J.H. Korhonen and T. Lempi{\"{a}}inen and P.R.J. {\"{O}}sterg{\aa}rd and C. Purcell and J. Rybicki and J. Suomela and P. Uzna{\'{n}}ski},
	title = {$\mathsf{LCL}$ problems on grids},
	journaltitle = {ACM Symposium on Principles of Distributed Computing (PODC)},
	date = {2017},
	pages = {101--110},
	addendum = {Full version: \url{https://arxiv.org/abs/1702.05456}},
}

@article{KST,
  title={Borel chromatic numbers},
  author={A.S. Kechris and S. Solecki and S. Todorcevic},
  journal={Adv. Math.},
  volume={141},
  number={1},
  pages={1--44},
  date={1999},
}

@unpublished{KechrisMarks,
	author = {A.S. Kechris and A.S. Marks},
	title = {Descriptive Graph Combinatorics},
	date = {2020},
	howpublished = {\url{https://math.berkeley.edu/~marks/papers/combinatorics20book.pdf} (preprint)},
}

@article{Pikh_survey,
	author = {O. Pikhurko},
	title = {Borel combinatorics of locally finite graphs},
	journaltitle = {Surveys in Combinatorics, 28th British Combinatorial Conference},
	editor = {K.K. Dabrowski {et al.}},
	pages = {267--319},
	date = {2021},
}

@article{Notices,
    author = {A. Bernshteyn},
	title = {Descriptive combinatorics and distributed algorithms},
	date = {2022},
	journaltitle = {Not. Am. Math. Soc.},
	volume = {69},
    number = {9},
	pages = {1496--1507},
}

@book{KechrisDST,
	author = {A.S. Kechris},
	title = {Classical Descriptive Set Theory},
	date = {1995},
	publisher = {Springer-Verlag},
	location = {New York},
}

@unpublished{AnushDST,
	author = {A. Tserunyan},
	title = {Introduction to Descriptive Set Theory},
	date = {2022},
	howpublished = {\url{https://www.math.mcgill.ca/atserunyan/Teaching_notes/dst_lectures.pdf} (preprint)},
}

@unpublished{paths,
    author = {J. Greb{\'{i}}k and V. Rozho{\v{n}}},
    title = {Classification of local problems on paths from the perspective of descriptive combinatorics},
    date = {2021},
    howpublished = {\url{https://arxiv.org/abs/2103.14112} (preprint)},
}

@article{trees,
	author = {S. Brandt and Y.-J. Chang and J. Greb{\'{i}}k and C. Grunau and V. Rozho{\v{n}} and Z. Vidny{\'{a}}nszky},
	title = {Local problems on trees from the perspectives of distributed algorithms, finitary factors, and descriptive combinatorics},
	journaltitle = {Innovations in Theoretical Computer Science Conference (ITCS)},
    pages = {29:1--29:26},
    date = {2022},
    addendum = {Full version: \url{https://arxiv.org/abs/2106.02066}},
}

@article{Ber_cont,
	author = {A. Bernshteyn},
	title = {Probabilistic constructions in continuous combinatorics and a bridge to distributed algorithms},
	date = {2023},
    journaltitle = {Adv. Math.},
    volume = {415},
    pages = {108895},
}

@article{Marks,
	author = {A.S. Marks},
	title = {A determinacy approach to Borel combinatorics},
	journaltitle = {J. Am. Math. Soc.},
	volume = {29},
	date = {2016},
	pages = {579--600},
}

@article{RV,
	author = {M. Rahman and B. Vir{\'{a}}g},
	title = {Local algorithms for independent sets are half-optimal},
	journaltitle = {Ann. Probab.},
	volume = {45},
	number = {3},
	pages = {1543--1577},
	date = {2017},
}

@article{BerDist,
	author = {A. Bernshteyn},
	title = {Distributed algorithms, the Lovász Local Lemma, and descriptive combinatorics},
	date = {2023},
	journaltitle = {Invent. Math.},
	volume = {233},
	pages = {495--542},
}

@article{Ber_Lebesgue,
    author = {A. Bernshteyn},
    title = {Borel fractional colorings of Schreier graphs},
	date = {2022},
    volume = {5},
	journaltitle = {Ann. H. Lebesgue},
    pages = {1151--1160},
}

@article{CM16,
	author = {C.T. Conley and B.D. Miller},
	title = {A bound on measurable chromatic numbers of locally finite Borel graphs},
	date = {2016},
	journaltitle = {Math. Res. Lett.},
	volume = {23},
	number = {6},
	pages = {1633--1644},
}

@unpublished{ContinuousAbelian,
	author = {S. Gao and S. Jackson and E. Krohne and B. Seward},
	title = {Continuous combinatorics of abelian group actions},
	date = {2018},
	howpublished = {\url{https://arxiv.org/abs/1803.03872} (preprint)},
    addendum = {To appear in Mem. Am. Math. Soc.},
}

@unpublished{BorelAbelian,
	author = {S. Gao and S. Jackson and E. Krohne and B. Seward},
	title = {Borel combinatorics of abelian group actions},
	date = {2024},
	howpublished = {\url{https://arxiv.org/pdf/2401.13866} (preprint)},
}

@article{GaoJack,
	author = {S. Gao and S. Jackson},
	title = {Countable Abelian group actions and hyperfinite equivalence relations},
	journaltitle = {Invent. Math.},
	volume = {201},
	number = {1},
	pages = {309--383},
	date = {2015},
}

@article{GRgrids,
    author = {J. Greb{\'{i}}k and V. Rozho{\v{n}}},
    title = {Local problems on grids from the perspective of distributed algorithms, finitary factors, and descriptive combinatorics},
    date = {2023},
    journaltitle = {Adv. Math.},
    volume = {431},
    pages = {109241},
}

@book{Oxtoby,
    author = {J.C. Oxtoby},
    title = {Measure and Category},
    location={New York},
    publisher = {Springer-Verlag},
    date = {1980},
}

@article{STD,
	author = {B. Seward and R.D. Tucker-Drob},
	title = {Borel structurability on the $2$-shift of a countable group},
	journaltitle = {Ann. Pure Appl. Log.},
	volume = {167},
	number = {1},
	date = {2016},
	pages = {1--21},
}

@inbook{Weiss_generic,
	author = {B. Weiss},
	title = {A survey of generic dynamics},
	booktitle = {Descriptive Set Theory and Dynamical Systems},
	date = {2000},
	editor = {Foreman, M. and Kechris, A.S. and Louveau, A. and Weiss, B.},
	series = {London Mathematical Society Lecture Note Series},
	pages = {273--291},
}

@article{Holroyd,
	author = {A.E. Holroyd and O. Schramm and D.B. Wilson},
	title = {Finitary coloring},
	journaltitle = {Ann. Probab.},
	volume = {45},
	number = {5},
	date = {2017},
	pages = {2867--2898},
}

@article{Holroyd1,
	author = {A.E. Holroyd},
	title = {One-dependent coloring by finitary factors},
	journaltitle = {Ann. Inst. H. Poincar{\'{e}} Probab. Statist.},
	volume = {53},
	number = {2},
	date = {2017},
	pages = {753--765},
}

@article{Spinka,
	author = {Y. Spinka},
	title = {Finitely dependent processes are finitary},
	journaltitle = {Ann. Probab.},
	volume = {48},
	number = {4},
	date = {2020},
	pages = {2088--2117},
}

@article{FIIDSchreier,
	author = {F. Bencs and A. Hru{\v{s}}kov{\'{a}} and L.M. T{\'{o}}th},
	title = {Factor-of-iid Schreier decorations of lattices in Euclidean spaces},
	journaltitle = {Discrete Math.},
	volume = {347},
	number = {9},
	date = {2024},
	pages = {\#114056},
}

@article{Manaster2,
	author = {A.B. Manaster and J.G. Rosenstein},
	title = {Effective matchmaking and $k$-chromatic graphs},
	journaltitle = {Proc. Am. Math. Soc.},
	volume = {39},
    number = {2},
	date = {1973},
	pages = {371--378},
}

@article{Schmerl1,
	author = {J.H. Schmerl},
	title = {Recursive colorings of graphs},
	journaltitle = {Can. J. Math.},
	volume = {32},
    number = {4},
	date = {1980},
	pages = {821--830},
}

@article{Schmerl2,
	author = {J.H. Schmerl},
	title = {The effective version of Brooks' Theorem},
	journaltitle = {Can. J. Math.},
	volume = {34},
    number = {5},
	date = {1982},
	pages = {1036--1046},
}

@inbook{Gasarch,
	author = {W. Gasarch},
	title = {Chapter 16. A survey of recursive combinatorics.},
	booktitle = {Handbook of Recursive Mathematics},
    volume = {2},
	date = {1998},
	pages = {1041--1176},
}

@unpublished{ASIalgorithms,
	author = {L. Qian and F. Weilacher},
	title = {Descriptive combinatorics, computable combinatorics, and ASI algorithms},
	howpublished = {\url{https://arxiv.org/abs/2206.08426} (preprint)},
	date = {2022},
}

@unpublished{bowen2021perfect,
	author = {M. Bowen and G. Kun and M. Sabok},
	title = {Perfect matchings in hyperfinite graphings},
	date = {2021},
	howpublished = {\url{https://arxiv.org/abs/2106.01988} (preprint)},
}

@article{BPZ,
	author = {M. Bowen and A. Poulin and J. Zomback},
	title = {One-ended spanning trees and definable combinatorics},
	date = {2024},
    journaltitle = {Trans. Am. Math. Soc.},
    volume = {377},
    pages = {8411--8431},
}

@article{gao2022forcing,
  title={Forcing constructions and countable Borel equivalence relations},
  author={S. Gao and S. Jackson and E. Krohne and B. Seward},
  journal={J. Symb. Log.},
  volume={87},
  number={3},
  pages={873--893},
  year={2022},
}

@article{marks2017borel,
  title={Borel circle squaring},
  author={A.S. Marks and S.T. Unger},
  journal={Ann. Math.},
  volume={186},
  number={2},
  pages={581--605},
  year={2017},
}

@unpublished{KechrisCBER,
	author = {A.S. Kechris},
	title = {The Theory of Countable Borel Equivalence Relations},
	date = {2023},
	howpublished = {\url{https://pma.caltech.edu/documents/5608/lectures_on_CBER12book.pdf} (preprint)},
}

@unpublished{FelixTrees,
  title={Computable vs descriptive combinatorics of local problems on trees},
  author={F. Weilacher},
  addendum={To appear in J. Symb. Log.},
  year={2022},
    howpublished={\url{https://arxiv.org/abs/2208.06689} (preprint)},
}

@article{MeasurableBrooks,
	author = {C.T. Conley and A.S. Marks and R.D. Tucker-Drob},
	title = {Brooks' theorem for measurable colorings},
	date = {2016},
	volume = {4},
	journaltitle = {Forum Math. Sigma},
}

@article{Kierstead,
	author = {H.A. Kierstead},
	title = {Recursive colorings of highly recursive graphs},
	journaltitle = {Can. J. Math.},
	volume = {33},
    number = {6},
	date = {1981},
	pages = {1279--1290},
}

@article{FelixTwoEnded,
	author = {F. Weilacher},
	title = {Descriptive chromatic numbers of locally finite and everywhere two-ended graphs},
	journaltitle = {Groups Geom. Dyn.},
	volume = {16},
    number = {1},
	date = {2022},
	pages = {141--152},
}

@article{BerGeneric,
	author = {A. Bernshteyn},
	title = {On Baire measurable colorings of group actions},
	journaltitle = {Ergod. Theory Dyn. Syst.},
	volume = {41},
    number = {3},
	date = {2021},
	pages = {818--845},
}

@unpublished{BerFelixASI,
	author = {A. Bernshteyn and F. Weilacher},
	title = {Borel versions of the Local Lemma and $\mathsf{LOCAL}$ algorithms for graphs of finite asymptotic separation index},
	howpublished = {\url{https://arxiv.org/abs/2308.14941} (preprint)},
	date = {2023},
}

@article{WeissHyperfinite,
    title={Measurable dynamics},
    author={B. Weiss},
    journaltitle={Contemp. Math.},
    volume={26},
    pages={395--421},
    date={1984},
}

@inbook{SS,
	author = {T.A. Slaman and J.R. Steel},
	date = {1988},
	title = {Definable functions on degrees},
	journaltitle = {Cabal Seminar 81--85, Lect. Notes Math.},
    editor = {A.S. Kechris and D.A. Martin and J.R. Steel},
	volume = {1333},
    publisher = {Springer},
	pages = {37--55},
}

@article{DJK,
    title = {The structure of hyperfinite Borel equivalence relations},
    author = {R. Dougherty and S. Jackson and A.S. Kechris},
    journaltitle = {Trans. Am. Math. Soc.},
    volume = {341},
    number = {1},
    date = {1994},
    pages = {193--225},
}

@article{jackson2002countable,
  title={Countable Borel equivalence relations},
  author={Jackson, S. and Kechris, A.S. and Louveau, A.},
  journaltitle={J. Math. Log.},
  volume={2},
  number={01},
  pages={1--80},
  date={2002},
}

@article{conley2020borel,
  title={Borel asymptotic dimension and hyperfinite equivalence relations},
  author = {Conley, C.T. and Jackson, S. and Marks, A.S. and Seward, B. and Tucker-Drob, R.D.},
  journal={Duke Math. J.},
  volume={172},
  number={16},
  pages={3175--3226},
  year={2023},
}

@book{Coloring1,
    author = {T.R. Jensen and B. Toft},
    title = {Graph Coloring Problems},
    publisher = {Wiley},
    date = {1995},
}

@book{Coloring2,
    author = {D.W. Cranston},
    title = {Graph Coloring Methods},
    date = {2024},
    url = {https://graphcoloringmethods.com/},
}

@article{ConleyKechris,
  title={Measurable chromatic and independence numbers for ergodic graphs and group actions},
  author = {C.T. Conley and A.S. Kechris},
  journal={Groups Geom. Dyn.},
  volume={7},
  number={1},
  pages={127--180},
  year={2013},
}

@article{ornstein.weiss,
  title={Entropy and isomorphism theorems for actions of amenable groups},
  author = {D.S. Ornstein and B. Weiss},
  journal={J. Anal. Math.},
  volume={48},
  pages={1--141},
  year={1987},
}

@book{Katok,
    author = {A. Katok and B. Hasselblatt},
    title = {Introduction to the Modern Theory of Dynamical Systems},
    publisher = {Cambridge University Press },
    date = {2012},
}

@article{ornstein1980ergodic,
  title={Ergodic theory of amenable group actions, I: The Rohlin lemma},
  author={Ornstein, D.S. and Weiss, B.},
  journaltitle={Bull. Am. Math. Soc.},
  volume={2},
   number={1},
  pages={161--164},
  date={1980},
}

@article{entropy,
  title={Entropy inequalities for factors of IID},
  author={{\'{A}}. Backhausz and B. Gerencs{\'{e}}r and V. Harangi},
  journaltitle={Groups Geom. Dyn.},
  volume={13},
   number={2},
  pages={389--414},
  date={2019},
}

@article{LyonsNazarov,
  title={Perfect matchings as IID factors on non-amenable groups},
  author={R. Lyons and F. Nazarov},
  journaltitle={Eur. J. Comb.},
  volume={32},
   number={7},
  pages={1115--1125},
  date={2011},
}

@unpublished{KunHall,
	author = {G. Kun},
	title = {The measurable Hall theorem fails for treeings},
	date = {2021},
	howpublished = {\url{https://arxiv.org/abs/2106.02013} (preprint)},
}

@article{LBowen,
  title={Finitary random interlacements and the Gaboriau--Lyons problem},
  author={L. Bowen},
  journaltitle={Geom. Funct. Anal.},
  volume={29},
  pages={659--689},
  date={2019},
}

@article{Hex,
  title={The game of Hex and the Brouwer fixed-point theorem},
  author={D. Gale},
  journaltitle={Am. Math. Mon.},
  volume={86},
    number={10},
  pages={818--827},
  date={1979},
}

@book{KechrisMiller,
	author = {A.S. Kechris and B.D. Miller},
	title = {Topics in Orbit Equivalence},
	publisher = {Springer-Verlag},
	location = {Berlin/Heidelberg},
	date = {2004},
}

@thesis{MillerThesis,
	author = {B.D. Miller},
	title = {Full Groups, Classification, and Equivalence Relations},
	type = {Ph.D.~Thesis},
	institution = {University of California, Berkeley},
	date = {2004},
	location = {Berkeley, CA},
    url = {https://glimmeffros.github.io/publications/dissertation.pdf},
}

@unpublished{RileyEffective,
	author = {R. Thornton},
	title = {$\Delta_1^1$-Effectivization in Borel combinatorics},
	date = {2021},
	howpublished = {\url{https://arxiv.org/abs/2105.04063} (preprint)},
    addendum = {To appear in J. Symb. Log.},
}

@article{RileyOrient,
  title={Orienting Borel graphs},
  author={R. Thornton},
  journaltitle={Proc. Am. Math. Soc.},
  volume={150},
    number={4},
  pages={1779--1793},
  date={2022},
}

@thesis{WeilacherThesis,
	author = {F. Weilacher},
	title = {Definable Combinatorics in Descriptive Set Theory, Computability Theory, and Beyond},
	type = {Ph.D.~Thesis},
	institution = {Carnegie Mellon University},
	date = {2024},
	location = {Pittsburgh, PA},
    url = {https://kilthub.cmu.edu/articles/thesis/Definable_combinatorics_in_descriptive_set_theory_computability_theory_and_beyond/27307407},
}

@unpublished{Sequel,
	author = {K. Berlow and A. Bernshteyn and C. Lyons and F. Weilacher},
	title = {Separating complexity classes of LCL problems on grids II},
	howpublished = {in preparation},
}

@unpublished{BaireMatchings,
    author = {M. Bowen and C.T. Conley and F. Weilacher},
    title = {Measurable regular subgraphs},
    date = {2024},
    howpublished = {\url{https://arxiv.org/abs/2408.09597} (preprint)},
}

@article{Oberwolfach,
  title={Mini-Workshop: Descriptive Combinatorics, LOCAL Algorithms and Random Processes},
  editor = {Greb{\'{i}}k, J. and Pikhurko, O. and Tserunyan, A.},
  journaltitle={Oberwolfach Reports},
  volume={19},
    number={4},
  pages={429--455},
  date={2023},
}

\end{document}